\documentclass[english]{smfbook}
\usepackage[all]{xy}
\usepackage{amssymb}
\usepackage{amsmath}

\catcode`\Ž=\active \def Ž{\'e}
\catcode`\=\active \def {\`e}
\catcode`\ˆ=\active \def ˆ{\`a} 
\catcode`\=\active \def {\`u}

\catcode`\=\active \def {\^e}        
\catcode`\"=\active \def "{\^{\i}}    
\catcode`\â=\active \def â{\^a}

\catcode`\"=\active \def "{\^{\i}}
\catcode`\ô=\active \def ô{\^o}
\catcode`\ž=\active \def ž{\^u}
\catcode`\é=\active \def é{\`E} 
\catcode`\æ=\active \def æ{\^E} 

\swapnumbers
\newtheorem{thm}[subsection]{Theorem}
\newtheorem{prop}[subsection]{Proposition}
\newtheorem{lemma}[subsection]{Lemma}
\newtheorem{cor}[subsection]{Corollary}

 \renewcommand\thesection{\thechapter.\arabic{section}}

\def\d{\succ}
\def\g{\prec}

\def\oo{\omega}
\def\dd{\delta}
\def\DD{\Delta}

\def\ss{\sigma}
\def\aa{\alpha}
\def\bb{\beta}

\def\t{\otimes}

\def\Vect{{\texttt{Vect }}}

\def\alg{\textrm{-alg}}
\def\bialg{\textrm{-bialg}}

\def\pt{\textrm{-}}
\def\tto{\longrightarrow}

\newcommand{\Vt}[1]{V^{\otimes #1}}

\def\row#1#2#3{(#1_{#2},\ldots,#1_{#3})}

\def\QQ{{\mathbb{Q}}}

\def\KK{{\mathbb{K}}}

\def\AA{{\mathcal{A}}}
\def\BB{{\mathcal{B}}}
\def\CC{{\mathcal{C}}}
\def\HH{{\mathcal{H}}}
\def\PP{{\mathcal{P}}}
\def\II{{\mathcal{I}}}

\def\XXX{{\mathcal{X}}}

\def\QQQ{{\mathcal{Q}}}
\def\ZZZ{{\mathcal{Z}}} 
\def\gg{\mathfrak{g}}
\def\cc{\gamma}

\def\Prim{\mathrm{Prim\, }}
\def\Id{\mathrm{Id }}
\def\id{\mathrm{id }}
\def\Hom{\mathrm{Hom}}

\def\Ind{\mathrm{Ind}}
\def\Indec{\mathrm{Indec\, }}

\def\Bi{\mathop{{\bf B}_{\infty}}}

\def\Im{\mathop{\rm Im}}
\def\Ker{\mathop{\rm Ker}}

\def\To{\overline {T}}
\def\Vect{\mathtt{Vect}}
\def\Vtn{V^{\t n}}

\def\sm{{\bf S}-module }
\def\sms{{\bf S}-modules }

\def\bisms{${\bf S}^{op}$-{\bf S}-modules }

\def\epi{\twoheadrightarrow}
\def\mono{\rightarrowtail}
\def\proj{\textrm {proj}}

\def\ca{$\CC^c\pt\AA$ }
\def\cab{$\CC^c\pt\AA$-bialgebra }
\def\cabs{$\CC^c\pt\AA$-bialgebras }
\def\CAP{$(\mathcal{C}, \mathcal{A}, \mathcal{P})$ }

\def\KSn{\KK[S_{n}]}

\def\--{\textrm{-}}

\def\Hone{\texttt{(H1)} }
\def\Htwoiso{\texttt{(H2iso)} }
\def\Htwoepi{\texttt{(H2epi)} }

\def\row#1#2#3{(#1_{#2},\ldots,#1_{#3})}

\newenvironment{proo}{\begin{trivlist} \item{\emph{Proof.}}}
  {\hfill $\square$ \end{trivlist}}
  
\def\cpbdeuxdeuxmudd{\vcenter{\xymatrix@R=6pt@C=6pt{
*{}\ar@{-}[d]&&*{}\ar@{-}[d]\\
*{}\ar@{-}[dr]&&*{}\ar@{-}[dl]\\
&*{}\ar@{-}[dd]&\\
&&\\
&*{}\ar@{-}[dl]\ar@{-}[dr]&\\
*{}\ar@{-}[d]&&*{}\ar@{-}[d]\\
*{}&&*{}
}}}

\def\cpbtroisquatre{\vcenter{\xymatrix@R=6pt@C=6pt{
&&&&\\
&*{}\ar@{-}[u]&&*{}\ar@{-}[u]&\\
&&\mu\ar@{-}[ul]\ar@{-}[uu]\ar@{-}[ur]&&\\
&&\dd\ar@{-}[u]\ar@{-}[dll]\ar@{-}[dl] \ar@{-}[drr]\ar@{-}[dr] &&\\
*{}\ar@{-}[d]&*{}\ar@{-}[d]&&*{}\ar@{-}[d]&*{}\ar@{-}[d]\\
&&&&
}}}

\def\cpbmultiple{\vcenter{\xymatrix@R=6pt@C=6pt{
&&&&    &&&&    &&&&     \\
&&&&  *{}\ar@{-}[u]  &&&&   *{}\ar@{-}[u]   &&&&    \\
*{}\ar@{-}[rrrrrrrrrrr] &*{}\ar@{-}[uu]&&*{}\ar@{-}[ur]&  *{}\ar@{-}[uu]  &*{}\ar@{-}[ul]&*{}\ar@{-}[urr]  &*{}\ar@{-}[ur]  &  *{}\ar@{-}[u]    &&*{}\ar@{-}[ull]  &*{}&    \\
&*{}&*{}&*{}&  *{} &*{}& \omega &*{}& *{}&*{}&*{}&*{}& *{}   &*{}&*{}&*{}& *{}   \\
*{}\ar@{-}[uu] \ar@{-}[rrrrrrrrrrr] &*{}&*{}&*{}&   *{} &*{}&*{}&*{}&  *{}&*{}& *{}    &*{}\ar@{-}[uu] &  *{}  \\
&&*{}\ar@{-}[d]\ar@{-}[ul]\ar@{-}[ur]  &&  *{}\ar@{-}[d] \ar@{-}[u] 
  &&*{}\ar@{-}[d] \ar@{-}[ul] \ar@{-}[u]\ar@{-}[ur]  &&&*{}\ar@{-}[d] \ar@{-}[ul] \ar@{-}[ur] &    \\
&&&&    &&&&      &&&&    \\
}}}

\def\cpbmult{\vcenter{\xymatrix@R=6pt@C=6pt{
&&                                                &&    							&&                                              & \\
&&                                                &&  *{}\ar@{-}[u]  					&&                                               & \\
*{}\ar@{-}[rrrrrrr] &*{}\ar@{-}[uu]&&*{}\ar@{-}[ur]&  *{}\ar@{-}[uu]  &*{}\ar@{-}[ul]& &          *{}\ar@{-}[dd]    &\\
&*{}&*{}                                       &  *{}	&  		\omega               &*{}&                                               &\\
*{}\ar@{-}[uu] \ar@{-}[rrrrrrr]  &*{}&*{}&*{}&   *{} &*{}&*{} &*{} \\
&&*{}\ar@{-}[d]\ar@{-}[ul]\ar@{-}[ur]  &&   &*{} \ar@{-}[ul]\ar@{-}[ur] *{}\ar@{-}[d]  &*{}    &\\
&&                           &&    &&   & \\
}}}

\def\cpbdeuxdeux{\vcenter{\xymatrix@R=2pt@C=2pt{
&&\\
*{}\ar@{-}[u]\ar@{-}[dr]&&*{}\ar@{-}[u]\ar@{-}[dl]\\
&*{}\ar@{-}[d]&\\
&*{}\ar@{-}[dl]\ar@{-}[dr]&\\
*{}\ar@{-}[d]&&*{}\ar@{-}[d]\\
&&
}}}

\def\cpbA{\vcenter{\xymatrix@R=2pt@C=2pt{
*{}\ar@{-}[ddddd]&&*{}\ar@{-}[ddddd]\\
&&\\
&&\\
&&\\
&&\\
*{}&&*{}
}}}

\def\cpbB{\vcenter{\xymatrix@R=2pt@C=2pt{
*{}\ar@{-}[dd]&&*{}\ar@{-}[dd]\\
&&\\
*{}\ar@{-}[drr]& &*{}\ar@{-}[dll]\\
*{}\ar@{-}[dd]& &*{}\ar@{-}[dd]\\
&&\\
*{}&&*{}
}}}

\def\cpbC{\vcenter{\xymatrix@R=2pt@C=2pt{
&&&&\\
&*{}\ar@{-}[dl]\ar@{-}[dr]\ar@{-}[u]&&&\\
*{}\ar@{-}[ddd]&&*{}\ar@{-}[d]&&\\
&&*{}\ar@{-}[dr]&&*{}\ar@{-}[dl]*{}\ar@{-}[uuu]\\
&&&*{}\ar@{-}[d]&\\
&&&&
}}}

\def\cpbD{\vcenter{\xymatrix@R=2pt@C=2pt{
&&&\\
&*{}\ar@{-}[dl]\ar@{-}[dr]\ar@{-}[u]&&*{}\ar@{-}[d]*{}\ar@{-}[u]\\
*{}\ar@{-}[d]&&*{}\ar@{-}[dr]&*{}\ar@{-}[dl]\\
*{}\ar@{-}[dr]&&*{}\ar@{-}[dl]&*{}\ar@{-}[dd]\\
&*{}\ar@{-}[d]&&\\
&&&
}}}

\def\cpbE{\vcenter{\xymatrix@R=2pt@C=2pt{
&&&&\\
&&&*{}\ar@{-}[u]\ar@{-}[dl]\ar@{-}[dr]&\\
&&*{}\ar@{-}[d]&&*{}\ar@{-}[ddd]\\
*{}\ar@{-}[uuu]\ar@{-}[dr]&&*{}\ar@{-}[dl]&&\\
&*{}\ar@{-}[d]&&&\\
&&&&
}}}

\def\cpbF{\vcenter{\xymatrix@R=2pt@C=2pt{
&&&\\
&&*{}\ar@{-}[u]&\\
*{}\ar@{-}[uu]*{}\ar@{-}[dr]&*{}\ar@{-}[dl]\ar@{-}[ur]&&*{}\ar@{-}[ul]\ar@{-}[d]\\
*{}\ar@{-}[dd]&*{}\ar@{-}[dr]&&*{}\ar@{-}[dl]\\
&&*{}\ar@{-}[d]&\\
&&&
}}}

\def\cpbG{\vcenter{\xymatrix@R=2pt@C=2pt{
&&&&&\\
&*{}\ar@{-}[dl]\ar@{-}[dr]\ar@{-}[u]&&&*{}\ar@{-}[dl]\ar@{-}[dr]\ar@{-}[u]&\\
*{}\ar@{-}[d]&&*{}\ar@{-}[dr]&*{}\ar@{-}[dl]&&*{}\ar@{-}[d]\\
*{}\ar@{-}[dr]&&*{}\ar@{-}[dl]&*{}\ar@{-}[dr]&&*{}\ar@{-}[dl]\\
&*{}\ar@{-}[d]&&&*{}\ar@{-}[d]&\\
&&&&&&
}}}

 \def\arbreA{\vcenter{\xymatrix@R=3pt@C=3pt{
&& \\
&*{}\ar@{-}[ur] \ar@{-}[ul] \ar@{-}[d]     &\\
&&
}}}

\def\arbreBA{\vcenter{\xymatrix@R=2pt@C=2pt{
&&&&\\
&&&*{}\ar@{-}[ul] & \\
&&*{}\ar@{-}[uurr] \ar@{-}[uull] \ar@{-}[d]     &&\\
&&&&
}}}

\def\arbreAB{\vcenter{\xymatrix@R=2pt@C=2pt{
&&&&\\
&*{}\ar@{-}[ur] &&& \\
&&*{}\ar@{-}[uurr] \ar@{-}[uull] \ar@{-}[d]     &&\\
&&&&
}}}

\def\arbreBB{\vcenter{\xymatrix@R=2pt@C=2pt{
&&&&\\
&&&& \\
&&*{}\ar@{-}[uurr] \ar@{-}[uull] \ar@{-}[d] \ar@{-}[uu]     &&\\
&&&&
}}}

\def\arbreABC{\vcenter{\xymatrix@R=1pt@C=1pt{
&&&&&&\\
&*{}\ar@{-}[ur] &&&&& \\
&&*{}\ar@{-}[uurr] &&&&\\
&&&*{}\ar@{-}[uuurrr] \ar@{-}[uuulll] \ar@{-}[d] &&&\\
&&&&&&
}}}

\def\arbreBAC{\vcenter{\xymatrix@R=1pt@C=1pt{
&&&&&&\\
&&&*{}\ar@{-}[ul] &&& \\
&&*{}\ar@{-}[uurr] &&&&\\
&&&*{}\ar@{-}[uuurrr] \ar@{-}[uuulll] \ar@{-}[d] &&&\\
&&&&&&
}}}

\def\arbreACA{\vcenter{\xymatrix@R=1pt@C=1pt{
&&&&&&\\
&*{}\ar@{-}[ur] &&&&*{}\ar@{-}[ul] & \\
&&&&&&\\
&&&*{}\ar@{-}[uuurrr] \ar@{-}[uuulll] \ar@{-}[d] &&&\\
&&&&&&
}}}

\def\arbreCAB{\vcenter{\xymatrix@R=1pt@C=1pt{
&&&&&&\\
&&&*{}\ar@{-}[ur] &&& \\
&&&&*{}\ar@{-}[uull] &&\\
&&&*{}\ar@{-}[uuurrr] \ar@{-}[uuulll] \ar@{-}[d] &&&\\
&&&&&&
}}}

\def\arbreCBA{\vcenter{\xymatrix@R=1pt@C=1pt{
&&&&&&\\
&&&&&*{}\ar@{-}[ul] & \\
&&&&*{}\ar@{-}[uull] &&\\
&&&*{}\ar@{-}[uuurrr] \ar@{-}[uuulll] \ar@{-}[d] &&&\\
&&&&&&
}}}

\def\boitenun{\vcenter{\xymatrix@R=6pt@C=6pt{
&&          &&  		&&                     \\
&&       &    &   &   &           & \\
*{}\ar@{-}[rrrrrr] &*{}\ar@{-}[uu]&*{}\ar@{-}[uu]& &  *{}\ar@{-}[uu]  &*{}\ar@{-}[uu]  &   *{}\ar@{-}[dd]    &\\
&*{}&*{}                                       & \mu&  	              &*{}&                                               \\
*{}\ar@{-}[uu] \ar@{-}[rrrrrr]  &*{}&& *{}\ar@{-}[dd]  &*{}&*{} &*{} \\
&&  &&   &  &   \\
&&                           &&  *{}  &&    \\
}}}

\def\boiteunm{\vcenter{\xymatrix@R=6pt@C=6pt{
&&          &&  		&&                     \\
&&       &    &   &   &           & \\
*{}\ar@{-}[rrrrrr] & && *{}\ar@{-}[uu]&     &  &   *{}\ar@{-}[dd]    &\\
&*{}&*{}                                       & \dd&  	              &*{}&                                               \\
*{}\ar@{-}[uu] \ar@{-}[rrrrrr]  &*{}& *{}\ar@{-}[dd] & *{}\ar@{-}[dd]  & *{}\ar@{-}[dd] &*{} &*{} \\
&&  &&   &  &   \\
&&                           &&  *{}  &&    \\
}}}

\def\boitenm{\vcenter{\xymatrix@R=6pt@C=6pt{
&&          &&  		&&                     \\
&&       &    &   &   &           & \\
*{}\ar@{-}[rrrrrr] &*{}\ar@{-}[uu]&*{}\ar@{-}[uu]& &  *{}\ar@{-}[uu]  &*{}\ar@{-}[uu]  &   *{}\ar@{-}[dd]    &\\&*{}&*{}                                       &&  	              &*{}&                                               \\
*{}\ar@{-}[uu] \ar@{-}[rrrrrr]  &*{}& *{}\ar@{-}[dd] & *{}\ar@{-}[dd]  & *{}\ar@{-}[dd] &*{} &*{} \\
&&  &&   &  &   \\
&&                           &&  *{}  &&    \\
}}}

\def\cpbtroistrois{\vcenter{\xymatrix@R=6pt@C=6pt{
&&&\\
*{}\ar@{-}[u]&&*{}\ar@{-}[u]&\\
&*{}\ar@{-}[ul] \ar@{-}[ur] \ar@{-}[uu]&&\\
&*{}\ar@{-}[u] \ar@{-}[dl]  \ar@{-}[dr] \ar@{-}[dd]&&\\
*{}\ar@{-}[d]&&*{}\ar@{-}[d]\\
&&&&
}}}

\def\cpbtroisID{\vcenter{\xymatrix@R=6pt@C=6pt{
\ar@{-}[ddddd]&\ar@{-}[ddddd]&\ar@{-}[ddddd]&\\
&&&\\
&&&\\
&&&\\
&&&\\
&&&
}}}

\def\cpbtroisdeux{\vcenter{\xymatrix@R=6pt@C=6pt{
&&&\\
*{}\ar@{-}[u]&&*{}\ar@{-}[u]&\\
&*{}\ar@{-}[ul] \ar@{-}[ur] \ar@{-}[uu]&&\\
&*{}\ar@{-}[u] \ar@{-}[dl]  \ar@{-}[dr]&&\\
*{}\ar@{-}[d]&&*{}\ar@{-}[d]\\
&&&&
}}}

\def\cpbdeuxtrois{\vcenter{\xymatrix@R=6pt@C=6pt{
&&&\\
*{}\ar@{-}[u]&&*{}\ar@{-}[u]&\\
&*{}\ar@{-}[ul] \ar@{-}[ur] &&\\
&*{}\ar@{-}[u] \ar@{-}[dl]  \ar@{-}[dr] \ar@{-}[dd]&&\\
*{}\ar@{-}[d]&&*{}\ar@{-}[d]\\
&&&&
}}}

\def\cpbtroisAB{\vcenter{\xymatrix@R=6pt@C=6pt{
&\ar@{-}[d]&&\ar@{-}[ddd]&\ar@{-}[ddd]\\
&*{}\ar@{-}[dd]\ar@{-}[dl]\ar@{-}[dr]&&&\\
*{}\ar@{-}[d]& &*{}\ar@{-}[d]&&\\
*{}&*{}&*{}&*{}&*{}\\
&*{}&*{}\ar@{-}[ul]\ar@{-}[u]\ar@{-}[ur]&\\
\ar@{-}[uu]&&\ar@{-}[u]&&\ar@{-}[uu]
}}}

\def\cpbtroisAC{\vcenter{\xymatrix@R=6pt@C=6pt{
&\ar@{-}[d]&&\ar@{-}[ddd]&\ar@{-}[ddd]\\
&*{}\ar@{-}[dd]\ar@{-}[dl]\ar@{-}[dr]&&&\\
*{}\ar@{-}[d]& &*{}\ar@{-}[d]&&\\
*{}&*{}&*{}&*{}&*{}\\
&&*{}&*{}\ar@{-}[ul]\ar@{-}[u]\ar@{-}[ur]\\
\ar@{-}[uu]&\ar@{-}[uu]&&\ar@{-}[u]&
}}}

\def\cpbtroisBA{\vcenter{\xymatrix@R=6pt@C=6pt{
\ar@{-}[ddd]&&\ar@{-}[d]&&\ar@{-}[ddd]\\
&&*{}\ar@{-}[dd]\ar@{-}[dl]\ar@{-}[dr]&&\\
&*{}\ar@{-}[d]& &*{}\ar@{-}[d]&\\
*{}&*{}&*{}&*{}&*{}\\
*{}&*{}\ar@{-}[ul]\ar@{-}[u]\ar@{-}[ur]&\\
&\ar@{-}[u]&&\ar@{-}[uu]&\ar@{-}[uu]
}}}
\def\cpbtroisBC{\vcenter{\xymatrix@R=6pt@C=6pt{
\ar@{-}[ddd]&&\ar@{-}[d]&&\ar@{-}[ddd]\\
&&*{}\ar@{-}[dd]\ar@{-}[dl]\ar@{-}[dr]&&\\
&*{}\ar@{-}[d]& &*{}\ar@{-}[d]&\\
*{}&*{}&*{}&*{}&*{}\\
&&*{}&*{}\ar@{-}[ul]\ar@{-}[u]\ar@{-}[ur]\\
\ar@{-}[uu]&\ar@{-}[uu]&&\ar@{-}[u]&
}}}
\def\cpbtroisCA{\vcenter{\xymatrix@R=6pt@C=6pt{
\ar@{-}[ddd]&\ar@{-}[ddd]&&\ar@{-}[d]&\\
&&&*{}\ar@{-}[dd]\ar@{-}[dl]\ar@{-}[dr]&\\
&&*{}\ar@{-}[d]& &*{}\ar@{-}[d]\\
*{}&*{}&*{}&*{}&*{}\\
*{}&*{}\ar@{-}[ul]\ar@{-}[u]\ar@{-}[ur]&\\
&\ar@{-}[u]&&\ar@{-}[uu]&\ar@{-}[uu]
}}}
\def\cpbtroisCB{\vcenter{\xymatrix@R=6pt@C=6pt{
\ar@{-}[ddd]&\ar@{-}[ddd]&&\ar@{-}[d]&\\
&&&*{}\ar@{-}[dd]\ar@{-}[dl]\ar@{-}[dr]&\\
&&*{}\ar@{-}[d]& &*{}\ar@{-}[d]\\
*{}&*{}&*{}&*{}&*{}\\
&*{}&*{}\ar@{-}[ul]\ar@{-}[u]\ar@{-}[ur]&\\
\ar@{-}[uu]&&\ar@{-}[u]&&\ar@{-}[uu]
}}}

\def\cpblosg{\vcenter{\xymatrix@R=2pt@C=2pt{
&*{}\ar@{-}[d]&&*{}\ar@{-}[ddddd]&\\
&*{}\ar@{-}[dl]\ar@{-}[dr]&&&\\
*{}&&*{}&&\\
&*{}\ar@{-}[ul]\ar@{-}[ur]&&&\\
&&&&\\
&*{}\ar@{-}[uu]&&*{}&
}}}

\def\cpblosd{\vcenter{\xymatrix@R=2pt@C=2pt{
*{}\ar@{-}[ddddd]&&*{}\ar@{-}[d]&&\\
&&*{}\ar@{-}[dl]\ar@{-}[dr]&&\\
&*{}&&*{}&\\
&&*{}\ar@{-}[ul]\ar@{-}[ur]&&\\
&&&&\\
*{}&&*{}\ar@{-}[uu]&*{}&
}}}

\def\cpbloslos{\vcenter{\xymatrix@R=2pt@C=2pt{
&*{}\ar@{-}[d]& & &*{}\ar@{-}[d]&\\
&*{}\ar@{-}[dl]\ar@{-}[dr]& & &*{}\ar@{-}[dl]\ar@{-}[dr]&\\
*{}&&*{} & *{}&&*{}\\
&*{}\ar@{-}[ul]\ar@{-}[ur]& & &*{}\ar@{-}[ul]\ar@{-}[ur]&\\
&& & &&\\
&*{}\ar@{-}[uu]& && *{}\ar@{-}[uu]&*{}
}}}

\def\cpbdeuxTrois{\vcenter{\xymatrix@R=2pt@C=2pt{
&&&\\
&*{}\ar@{-}[u]\ar@{-}[dr]&&*{}\ar@{-}[u]\ar@{-}[dl]\\
&&*{}\ar@{-}[d]&\\
&&*{}\ar@{-}[dl]\ar@{-}[dr]&\\
&*{}\ar@{-}[d]&&*{}\ar@{-}[ddd]\\
&*{}\ar@{-}[dl]\ar@{-}[dr]&&\\
*{}\ar@{-}[d]&&*{}\ar@{-}[d]&\\
&&&
}}}

\def\cpbDeuxDeux#1#2{\vcenter{\xymatrix@R=2pt@C=2pt{
*{}\ar@{-}[d]&&*{}\ar@{-}[d]\\
*{}\ar@{-}[dr]&&*{}\ar@{-}[dl]\\
&#1\ar@{-}[dd]&\\
&&\\
&#2\ar@{-}[dl]\ar@{-}[dr]&\\
*{}\ar@{-}[d]&&*{}\ar@{-}[d]\\
*{}&&*{}
}}}

\def\cpbMa{\vcenter{\xymatrix@R=2pt@C=2pt{
&*{}\ar@{-}[ddddd]&&*{}\ar@{-}[ddddddd]\\
&&&\\
&&&\\
&&&\\
&&&\\
&*{}\ar@{-}[dl]\ar@{-}[dr]&&\\
*{}\ar@{-}[d]&&*{}\ar@{-}[d]&\\
&&&
}}}

\def\cpbMb{\vcenter{\xymatrix@R=2pt@C=2pt{
&&*{}\ar@{-}[dd]&&*{}\ar@{-}[dd]\\
&&&&\\
&&*{}\ar@{-}[dl]\ar@{-}[dr]&&*{}\ar@{-}[dl]\\
&*{}\ar@{-}[dd]&&*{}\ar@{-}[dddd]&\\
& &&&\\
&*{}\ar@{-}[dl]\ar@{-}[dr]&&&\\
*{}\ar@{-}[d]&&*{}\ar@{-}[d]&&\\
&&&
}}}

\def\cpbMc{\vcenter{\xymatrix@R=2pt@C=2pt{
&&*{}\ar@{-}[dd]&&*{}\ar@{-}[ddd]&\\
&&&&&\\
&&*{}\ar@{-}[dl]\ar@{-}[ddrr]&&&\\
&*{}&&&*{}\ar@{-}[dl]&\\
&*{}\ar@{-}[u]\ar@{-}[dr]&&*{}\ar@{-}[dl]&*{}\ar@{-}[dddd]&\\
&&*{}\ar@{-}[d]&&&\\
&&*{}\ar@{-}[dl]\ar@{-}[dr]&&&\\
&*{}\ar@{-}[d]&&*{}\ar@{-}[d]&&\\
&&&&&\\
}}}

\def\cpbMd{\vcenter{\xymatrix@R=2pt@C=2pt{
&&*{}\ar@{-}[dd]&&*{}\ar@{-}[ddd]&\\
&&&&&\\
&&*{}\ar@{-}[dl]\ar@{-}[ddrr]&&&\\
&*{}\ar@{-}[ddddd]&&&*{}\ar@{-}[dl]&\\
& &&*{}\ar@{-}[dl]&*{}\ar@{-}[dddd]&\\
&&*{}\ar@{-}[ddd] &&&\\
&& &&&\\
& && &&\\
&&&&&\\
}}}

\def\cpbMe{\vcenter{\xymatrix@R=2pt@C=2pt{
&&*{}\ar@{-}[dd]&&*{}\ar@{-}[ddd]\\
&&&&\\
&&*{}\ar@{-}[dl]\ar@{-}[ddrr]&&\\
&*{}\ar@{-}[dd]&&&*{}\ar@{-}[dl]\\
& &&*{}\ar@{-}[d]&*{}\ar@{-}[ddd]\\
&*{}\ar@{-}[dl]\ar@{-}[dr]&&*{}\ar@{-}[dl]\&\\
*{}\ar@{-}[d]&&*{}\ar@{-}[d]&&\\
&&&&
}}}

\def\cpbMf{\vcenter{\xymatrix@R=2pt@C=2pt{
&&*{}\ar@{-}[d]&&*{}\ar@{-}[dd]\\
&&*{}\ar@{-}[dl]\ar@{-}[ddrr]&&\\
&*{}\ar@{-}[dd]&&&*{}\ar@{-}[dl]\\
& &&*{}\ar@{-}[d]&*{}\ar@{-}[dddd]\\
&*{}\ar@{-}[dl]\ar@{-}[ddrr]&&*{}\ar@{-}[ddll] &\\
*{}\ar@{-}[dr]&&&&\\
&*{}\ar@{-}[d]&&*{}\ar@{-}[d]&\\
&&&&
}}}

\def\cpbCC#1#2{\vcenter{\xymatrix@R=2pt@C=2pt{
&*{}\ar@{-}[d]&&&*{}\ar@{-}[ddd]\\
&#1\ar@{-}[dl]\ar@{-}[dr]&&&\\
*{}\ar@{-}[ddd]&&*{}\ar@{-}[d]&&\\
&&*{}\ar@{-}[dr]&&*{}\ar@{-}[dl]\\
&&&#2\ar@{-}[d]&\\
*{}&&&*{}&
}}}

\def\cpbDD#1#2{\vcenter{\xymatrix@R=2pt@C=2pt{
&*{}\ar@{-}[d]&&*{}\ar@{-}[d]\\
&#1\ar@{-}[dl]\ar@{-}[dr]&&*{}\ar@{-}[d]\\
*{}\ar@{-}[d]&&*{}\ar@{-}[dr]&*{}\ar@{-}[dl]\\
*{}\ar@{-}[dr]&&*{}\ar@{-}[dl]&*{}\ar@{-}[dd]\\
&#2\ar@{-}[d]&&\\
&*{}&&*{}
}}}

\def\cpbEE#1#2{\vcenter{\xymatrix@R=2pt@C=2pt{
*{}\ar@{-}[ddd]&&&*{}\ar@{-}[d]&\\
&&&#1\ar@{-}[dl]\ar@{-}[dr]&\\
&&*{}\ar@{-}[d]&&*{}\ar@{-}[ddd]\\
*{}\ar@{-}[dr]&&*{}\ar@{-}[dl]&&\\
&#2\ar@{-}[d]&&&\\
&*{}&&&*{}
}}}

\def\cpbFF#1#2{\vcenter{\xymatrix@R=2pt@C=2pt{
&&*{}\ar@{-}[d]&\\
&&#1&\\
*{}\ar@{-}[uu]*{}\ar@{-}[dr]&*{}\ar@{-}[dl]\ar@{-}[ur]&&*{}\ar@{-}[ul]\ar@{-}[d]\\
*{}\ar@{-}[dd]&*{}\ar@{-}[dr]&&*{}\ar@{-}[dl]\\
&&#2\ar@{-}[d]&\\
*{}&&*{}&
}}}

\def\cpbGG#1#2#3#4{\vcenter{\xymatrix@R=2pt@C=2pt{
&*{}\ar@{-}[d]&&&*{}\ar@{-}[d]&\\
&#1\ar@{-}[dl]\ar@{-}[dr]&&&#3\ar@{-}[dl]\ar@{-}[dr]&\\
*{}\ar@{-}[d]&&*{}\ar@{-}[dr]&*{}\ar@{-}[dl]&&*{}\ar@{-}[d]\\
*{}\ar@{-}[dr]&&*{}\ar@{-}[dl]&*{}\ar@{-}[dr]&&*{}\ar@{-}[dl]\\
&#2\ar@{-}[d]&&&#4\ar@{-}[d]&\\
&*{}&&&*{}&
}}}

\makeindex

\author[J.-L. Loday]{Jean-Louis Loday}
\address{Institut de Recherche Math\'ematique Avanc\'ee\\
    CNRS et Universit\'e de Strasbourg\\
    7 rue R. Descartes\\
    67084 Strasbourg, France}
\email{loday@math.u-strasbg.fr}
\urladdr{www-irma.u-strasbg.fr/{$\sim$}loday/}

\title[Generalized bialgebras]{Generalized bialgebras and triples of operads}
\alttitle{Big\`ebres g\' en\' eralis\' ees et triples d'op\' erades}

\date{\today}
\begin{document}

\frontmatter

\begin{abstract} 
We introduce the notion of generalized bialgebra, which includes the classical notion of bialgebra (Hopf algebra) and many others, like, for instance, the tensor algebra equipped with the deconcatenation as coproduct. We prove that,  under some mild conditions, a connected generalized bialgebra is completely determined by its primitive part. This structure theorem extends the classical Poincar\'e-Birkhoff-Witt theorem and Cartier-Milnor-Moore theorem, valid for cocommutative bialgebras, to a large class of generalized bialgebras. 

Technically we work in the theory of operads which allows us to state our main results and permits us to give it a conceptual proof. A generalized bialgebra type is determined by two operads: one for the coalgebra structure ${\mathcal{C}}$, and one for the algebra structure ${\mathcal{A}}$. There is also a compatibility relation relating the two. Under some conditions, the primitive part of such a generalized bialgebra is an algebra over some sub-operad of ${\mathcal{A}}$, denoted ${\mathcal{P}}$. The structure theorem gives conditions under which a connected generalized bialgebra is cofree (as a connected ${\mathcal{C}}$-coalgebra) and can be re-constructed out of its primitive part by means of an enveloping functor from ${\mathcal{P}}$-algebras to ${\mathcal{A}}$-algebras. The classical case is $({\mathcal{C}}, {\mathcal{A}}, {\mathcal{P}}) = (Com, As, Lie)$.

This structure theorem unifies several results, generalizing the PBW and the CMM theorems, scattered in the literature. We treat many explicit examples and suggest a few conjectures.
\end{abstract}

\begin{altabstract}
On introduit la notion de bigbre gŽnŽralisŽe, qui inclut la notion de bigbre classique (algbre de Hopf) et bien d'autres, comme, par exemple, l'algbre tensorielle munie de la dŽconcatŽnation comme coproduit. On montre que, sous des hypothses raisonnables, une bigbre gŽnŽralisŽe connexe est entirement dŽterminŽe par sa partie primitive. Ce thŽorme de structure Žtend ˆ la fois le thŽorme classique de PoincarŽ-Birkhoff-Witt et le thŽorme de Cartier-Milnor-Moore valables pour les bigbres cocommutatives, ˆ une large classe de bigbres gŽnŽralisŽes. 

On travaille dans le cadre de la thŽorie des opŽrades qui nous permet d'Žnoncer le rŽsultat principal et d'en donner une dŽmonstration conceptuelle. Un type de bigbres gŽnŽralisŽes est dŽterminŽ par deux opŽrades, l'une pour la structure de cogbre, notŽe ${\mathcal{C}}$, l'autre pour la structure d'algbre, notŽe ${\mathcal{A}}$. Ces deux structures sont reliŽes par certaines relations de compatibilitŽ. Le thŽorme de structure donne des conditions sous lesquelles une bigbre gŽnŽralisŽe connexe est colibre (en tant que ${\mathcal{C}}$-cogbre connexe) et peut \^{e}tre re-construite ˆ partir de sa partie primitive gr\^{a}ce ˆ un foncteur du type ``algbre envelopante'' des ${\mathcal{P}}$-algbres dans les ${\mathcal{A}}$-algbres. Le cas classique est $({\mathcal{C}}, {\mathcal{A}}, {\mathcal{P}}) = (Com, As, Lie)$.

Ce thŽorme de structure unifie plusieurs gŽnŽralisations du thŽorme PBW et du thŽorme CMM dŽjˆ prŽsentes dans la littŽrature. On donne plusieurs exemples explicites et on formule quelques conjectures.
\end{altabstract}



\altkeywords{Bigbre, bigbre gŽnŽralisŽe, algbre de Hopf, thŽorme de Cartier-Milnor-Moore, thŽorme de Poincar\'e-Birkhoff-Witt, opŽrade, prop, triple d'opŽrades, partie primitive, dendriforme, algbre dupliciale, algbre prŽ-Lie, algbre de Zinbiel, algbre magmatique, arbre, algbre non-associative.}

\thanks{
Many thanks to Mar\' \i a Ronco and to Bruno Vallette for numerous conversations on bialgebras and operads. Thanks to E.\ Burgunder, B.\ Fresse, R.\ Holtkamp, Y.\ Lafont, and M.\ Livernet for their comments. A special thank to F.\ Goichot for his careful reading. This work has been partially supported by the  ``Agence Nationale de la Recherche".
We thank the anonymous referee for his/her helpful corrections and comments.}

\maketitle

\tableofcontents

\mainmatter

\chapter*{Introduction} \label{S:int} 

The aim of this monograph is to prove that, under some simple conditions, there is a structure theorem for generalized bialgebras.

First we introduce the notion of ``generalized bialgebras", which includes the classical notions of bialgebras, Lie bialgebras, infinitesimal bialgebras, dendriform bialgebras and many others. A type of generalized bialgebras is determined by the coalgebra structure $\CC^c$, the algebra structure $\AA$ and the compatibility relations between the operations and the cooperations. For $\CC^c=As^c$ (coassociative coalgebra) and $\AA=As$ (associative algebra) with Hopf compatibility relation, we get the classical notion of bialgebra (Hopf algebra). In the general case we make the following assumption:

\smallskip

\noindent {\texttt{(H0)}}\  \emph{the compatibility relations are distributive.}

\smallskip

\noindent It means that any composition of an operation followed by a cooperation can be rewritten as cooperations first and then operations. So is the Hopf relation:  $\DD(xy)=\DD(x)\DD(y)$. In several cases the compatibility relation is equivalent to: the cooperation is an algebra morphism, but we deal here with a far more general situation.
Then, we make an assumption on the free $\AA$-algebra with respect to the bialgebra structure:

\smallskip

\noindent \Hone\  \emph{the free $\AA$-algebra $\AA(V)$ is naturally a \cab\!.}

\smallskip

\noindent For any generalized bialgebra $\HH$ there is a notion of primitive part given by
$$\Prim \HH := \{x\in \HH\ | \ \dd(x)=0 \textrm{ for any cooperation } \dd \textrm{ of arity } \geq 2\}.$$
At this point we are able to determine a new algebra structure, denoted $\PP$, such that the primitive part of any  \cab is a $\PP$-algebra. In other words we show that the $\AA$-operations, which are well-defined on the primitive part of any bialgebra, are stable by composition. Of course we get $\PP= Lie$ in the classical case. One should observe that, even when the types $\AA$ and $\CC$ are described by explicit generators and relations, there is no obvious way to get such a presentation for the type $\PP$. Therefore one needs to work with ``abstract types of algebras", that is with \emph{algebraic operads}. In this setting we introduce the notion of \emph{connected coalgebra}, which generalizes the similar notion introduced by Quillen when the generating cooperations are binary.

The forgetful functor  $\AA\alg \to \PP\alg$ from the category of $\AA$-algebras to the category of $\PP$-algebras admits a left adjoint which we denote by $U:\PP\alg \to \AA\alg$. 
The main result unravels the algebraic structure and the coalgebraic structure of any connected \cab under the hypothesis:

\smallskip

 \noindent \Htwoepi\  \emph{the coalgebra map  $\varphi(V):\AA(V) \to \CC^c(V)$ is split surjective.}

\smallskip

\noindent More precisely we get the
\smallskip

\noindent {\bf Structure Theorem for generalized bialgebras}.   {\it  Let $\CC^c\-- \AA$ be a bialgebra type which satisfies {\texttt{(H0)}}, \Hone and  \Htwoepi.
 
 Then, for any \cab $\HH$ with primitive part $\Prim \HH$, the following are equivalent:
 
 \noindent (a) $\HH$ is connected,\\
 (b)  $\HH$ is isomorphic to $U(\Prim \HH)$,\\
(c)  $\HH$ is cofree over its primitive part, i.e.\ isomorphic to $\CC^c(\Prim \HH)$.}

\medskip

As said above, the tool to determine $\PP$, and also to prove the structure theorem, is the \emph{operad} theory. A triple of operads $(\CC, \AA, \PP)$ as above is said to be \emph{good} if the structure theorem holds.

The case of (classical) cocommutative bialgebras (i.e.\ $Com^c\pt As$-bialgebras, with Hopf compatibility relation) is well-known. Here $\PP$ is $Lie$, that is, the primitive part of a classical bialgebra is a Lie algebra. So the triple
$$(Com, As, Lie)$$
is an example of a good triple of operads. The functor $U$ is the universal enveloping functor
$$U: Lie\alg \to As\alg.$$
The isomorphism $\HH\cong U(\Prim \HH)$ is the Cartier-Milnor-Moore theorem. The isomorphism $U(\gg)\cong S^c(\gg)$, where $\gg$ is a Lie algebra, is essentially the Poincar\'e-Birkhoff-Witt theorem. It implies that, given a basis for $\gg$, we can make up a basis for $U(\gg)$ from the commutative polynomials over a basis of $\gg$ (classical PBW theorem).
\medskip
In many cases, like the one above, characteristic zero is a necessary assumption, but here is an example of  a good triple of operads which is valid over any field $\KK$. The algebra type has two associative operations denoted, $x\g y$ and $x \d y$ , which satisfy moreover the relation
$$(x\d y)\g z = x\d (y\g z ).$$
We call them \emph{duplicial algebras} and we denote by $Dup$ the associated operad. The free algebra admits an elegant description in terms of the Over and the Under operation on planar binary trees. The coalgebra type is determined by a coassociative cooperation $\dd$ and the compatibility relation is the nonunital infinitesimal relation for both pairs $( \dd, \g)$ and $(\dd, \d)$ (see below and \ref{uib}). We can show that there is a structure theorem in this case and that the primitive structure is simply a magmatic structure. The magmatic operation $x\cdot y$ is given by $x\cdot y:= x\d y - x \g y$.
In short, there is a good triple
$$(As, Dup, Mag)\ .$$
This example, which is treated in details in Chapter \ref{ch:duplicial},  is remarkable, not only because all the operads are binary, quadratic and Koszul (like in the $(Com, As, Lie)$ triple), but because they are also set-theoretic and nonsymmetric.

\medskip

When the bialgebra type satisfies the stronger hypothesis

\smallskip

\noindent  \Htwoiso \emph{the $\CC^c$-coalgebra map $\varphi(V):\AA(V)\to \CC^c(V)$ is an isomorphism,}
 
 \smallskip
 
\noindent then the operad $\PP$ is trivial, $\PP=Vect$,  that is $\PP\alg$ is the category $\Vect$ of vector spaces. The triple is $(\CC, \AA, Vect)$, and the structure theorem becomes the 
 
 \smallskip
 
\noindent {\bf Rigidity Theorem for generalized bialgebras}.   {\it  Let $\CC^c\-- \AA$ be a bialgebra type which satisfies {\texttt{(H0)}}, \Hone and \Htwoiso. Then any connected \cab $\HH$ is free and cofree:}
$$ \AA(\Prim \HH)\cong \HH \cong \CC^c(\Prim \HH)\ .$$

\medskip

There are many good triples of the form $(\AA, \AA, Vect)$, for instance
\begin{eqnarray*}
&(Com, Com, Vect),&\\
&(As, As, Vect),&\\
&(Lie, Lie, Vect).&
\end{eqnarray*}
When $\AA=Com$ the compatibility relation is the Hopf relation and the rigidity theorem is the classical Hopf-Borel theorem (originally phrased in the framework of graded vector spaces). When $\AA=As$, then the compatibility relation that we take is not the Hopf relation, but the \emph{nonunital infinitesimal relation} which reads:
 $$\dd(x y) = x \t y + x_{(1)}\t x_{(2)} y + x  y_{(1)}\t y_{(2)}\ ,$$
where $\dd$ is the comultiplication and we have put $\dd(x) = x_{(1)}\t x_{(2)}$.  When $\AA=Lie$  the compatibility relation that we take is not the cocyle relation (giving rise to the notion of Lie bialgebras), but a new one (see \ref{compatibilityLieLie}).

\bigskip

Let us say a few words about the proofs. As said before we work in the efficient and well-adapted language of operad theory. The key point in the proof of the main theorem is a conceptual construction of a universal idempotent $e_{\HH} : \HH \to \HH$ which is functorial in the bialgebra $\HH$ and which does not depend on a presentation of the operads $\CC$ and $\AA$. Its image is the primitive part of $\HH$ and so it permits us to construct a morphism: $\HH\to \CC^c(\HH)$. In the other direction the morphism $U(\Prim \HH)\to \HH$ is induced by the inclusion $\Prim \HH \mono \HH$.  

In the classical case the idempotent so obtained is precisely the \emph{Eulerian idempotent}, whose construction is based on the logarithmic series. It is an important object since it provides an explicit description of the Baker-Campbell-Hausdorff formula \cite{JLLeuler}, and it permits us to split the Hochschild and the cyclic chain complexes  \cite{JLLop, HC}. See also \cite{Burgunder07} for an application to the Kashiwara-Vergne conjecture. In the case of the triple $(As, 2as, MB)$ the idempotent is Ronco's idempotent, based on the geometric series. 
Our construction gives an analogue of the Eulerian idempotent for each triple of operads.

We give several variations of our main theorem and many examples. We show that several results on bialgebras in the literature can be interpreted in terms of triples of operads. Our conceptual proof encompasses many ad hoc proofs of particular cases.

\bigskip

Here is the {\bf content} of this monograph. The {\bf first} chapter contains elementary facts about ``types of algebras and bialgebras" from the operadic point of view. The proofs of the theorems are performed in this framework, not only because of its efficiency, but also because some of the types of algebras (namely the primitive ones) that we encounter are not defined by generators and relations, but come as the kernel of some operad morphism. We introduce the notion of ``connected coalgebra" used in the hypotheses of the main theorem. The reader who is fluent in operad theory can easily bypass this chapter.

The {\bf second} chapter contains the main results of this monograph together with their proof. First, we study the algebraic structure of the primitive part of a  generalized bialgebra of type $\CC^c\--\AA$. In general a product of two primitive elements is not primitive. However the primitive part is stable under some operations. We determine all of them under the hypotheses {\texttt{(H0)}} and \Hone and we get the ``maximal" algebraic structure for the primitive part. We call it the primitive operad and denote it by $\Prim_{\CC}\AA$ or $\PP$. 

Then we study  the generalized bialgebra types which satisfy the hypothesis \Htwoiso. Though it will become a particular case of the general theorem, we prefer to treat it independently, because of its importance and because it is the key part of the proof of the general case. The result is the  \emph{rigidity theorem} for triples of the form $(\CC, \ZZZ, Vect)$. Then we move to the structure theorem. We establish that the conditions  {\texttt{(H0)}}, \Hone and \Htwoepi  ensure that the \emph{structure theorem}, referred to above, is valid for the $\CC^c\pt \AA$-bialgebras. 

There are several consequences to the structure theorem. For instance any good triple of operads \CAP gives rise to a natural isomorphism
$$\AA(V) \cong \CC^c \circ \PP(V).$$
It generalizes the classical fact that the underlying vector space of the symmetric algebra over the free Lie algebra is isomorphic to the tensor module
$$T(V)\cong S(Lie(V)).$$ 
Equivalently, we have an isomorphism of functors and/or of \sms: $\AA \cong \CC^c \circ \PP$. Taking the Frobenius characteristic (resp.\ the generating series) gives interesting functional equations of symmetric functions (resp.\ of power series).

In the known cases the proof of the structure theorem uses an ad hoc construction of an idempotent, which depends very much on the type of bialgebras at hand. The key point of our proof is to construct an ``abstract'' idempotent which works universally.
\smallskip

{\bf Third} Chapter. In the first three sections of this chapter we give  recipes to construct good triples of operads. An important consequence of our formulation is that, starting with a good triple $(\CC, \AA, \PP)$, we can construct many others by moding out by some operadic ideal. If $J$ is an operadic ideal of $\AA$ generated by primitive operations, then $(\CC, \AA/J, \Prim_{\CC}\AA/J)$ is also a good triple. In particular any good triple \CAP determines a good triple of the form $(\CC, \ZZZ, Vect)$, where $\ZZZ = \AA/(\overline \PP)$.

Assuming that the tensor product of two $\AA$-algebras is still a $\AA$-algebra (Hopf operad, multiplicative operad, for instance), there is a natural way of constructing a notion of $As^c\--\AA$-bialgebra which verifies {\texttt{(H0)}} and \Hone.  

For quadratic operads there is a notion of Koszul dual operad. We explain how this construction should permit us to construct new triples of operads.

In order to keep the proofs into the most simple  form we treated the case of algebraic operads in vector spaces over a characteristic zero field. But the structure theorem admits several generalizations. First, if we work with regular operads, then the characteristic zero hypothesis is not necessary anymore. 
Second, the tensor category of vector spaces $\Vect$ can be replaced by any linear 
symmetric monoidal category, for example the category of sign-graded vector spaces (super vector spaces) or the category of {\bf S}-modules. The formulas are the same provided that one applies the Koszul sign rule. Third, in characteristic $p$ it is expected that similar results hold.

The structure theorem can be ``dualized'' in the sense that the role of the algebra structure and the coalgebra structure are exchanged. The role of the primitive part is played by the indecomposable part.

In the last two sections we explain the relationship with the theory of ``rewriting systems'' and we give some application to the representation theory of the symmetric groups.

\smallskip

In the {\bf fourth} chapter we study some explicit examples in details. We show how several results in the literature can be interpreted as giving rise to a good triple of operads. Any good triple \CAP gives rise to a quotient triple of the form $(\CC, \ZZZ, Vect)$ where $\ZZZ= \AA/ (\bar\PP)$. We put in the same section the triples which have the same quotient triple:

$\bullet$ $(Com, Com, Vect)$ and Hopf compatibility relation. This section deals with the classical case of $Com^c\pt As$-bialgebras (cocommutative bialgebras) and some of its variations: $(Com, Parastat, NLie)$, $(Com, Mag, Sabinin)$. Our structure theorem for the triple $(Com, As, Lie)$ is equivalent to the Poincar\'e-Birkhoff-Witt (PBW) theorem plus the Cartier-Milnor-Moore (CMM) theorem. We show that our universal idempotent identifies to the Eulerian idempotent.

$\bullet$ $(As,As,Vect)$ and nonunital infinitesimal compatibility relation. It contains the case of 2-associative bialgebras, that is the triple $(As, 2as, Brace)$, and also $(As, Dipt, MB)$, $(As, Mag, MagFine)$. It is important because it permits us to handle the structure of classical cofree Hopf algebras \cite{LRstr}.

$\bullet$ $(As,Zinb,Vect)$ and semi-Hopf compatibility relation. It contains the dendriform and dipterous bialgebras. It is interesting for its role in the study of the graph-complexes \`a la Kontsevich obtained by replacing the Lie homology by the Leibniz homology \cite{HC, Burgunder08}.

$\bullet$ $(Lie,Lie,Vect)$ and the Lily compatibility relation. This is a completely new case. It gives a criterion to show that a Lie algebra is free.

$\bullet$ $(Nap,PreLie,Vect)$ and the Livernet compatibility relation. It is due to M. Livernet \cite{Livernet06}. A variation
 $(Nap,Mag,\Prim_{NAP}PreLie)$ needs more work to find a small presentation of the primitive operad.

$\bullet$ Then we survey some examples of the form $(\AA, \AA, Vect)$. We formulate a conjecture related to a question of M. Markl and we introduce an example coming from computer sciences (interchange bialgebra). Finally we present a triple involving $k$-ary operations and cooperations.

\smallskip

In the final Chapter we treat in details the \emph{duplicial bialgebras} which give rise to the triple  $(As, Dup, Mag)$ mentioned before. We prove that this triple is good and we make explicit the analogue of the PBW isomorphism. We prove that  the operad $Dup$ is Koszul.
We treat the case $(Dup,Dup,Vect)$ and  we comment on further generalizations (operadic quantization).

In the Appendix we provide a tableau of compatibility relations and a tableau of triples summarizing the examples treated in Chapters 4 and 5.

\bigskip

\noindent {\bf Notation, convention.} In this monograph $\KK$ is a field, which is, sometimes,  supposed to be of characteristic zero.  Its unit is denoted $1_\KK$ or just 1. All vector spaces 
are over $\KK$ and the category of vector spaces is denoted by $\Vect$\index{Vect@$\Vect$}. We often say ``space'' in place of ``vector space''. An injective linear map (monomorphism) is denoted by $\mono$ and a surjective linear map (epimorphism) is denoted by $\epi$. The space spanned by the elements of a set $X$ is denoted $\KK[X]$. The tensor product of vector spaces over $\KK$
is denoted by $\t$. The tensor product of $n$ copies of the space $V$ is
denoted by $\Vt n$. For $v_i\in V$ the element $v_1\t \cdots \t v_n$ of 
$\Vt n$ is denoted by  $(v_1, \ldots , v_n)$ or
simply by $v_1 \ldots v_n$. For instance in the \emph{tensor module}\index{tensor module}
$$T(V) := \KK \oplus V \oplus \cdots \oplus V^{\t n} \oplus \cdots $$
we denote by $v_1 \ldots v_n$ an element of $V^{\t n}$, but in $T(V)^{\t k}$ we denote by 
$v_1\t \cdots \t v_k$ the element such that $v_{i}\in V \subset T(V)$ is in the $i$th factor.
 The \emph{ reduced tensor module}\index{reduced tensor module}
$$\overline T(V) :=  V \oplus \cdots \oplus V^{\t n} \oplus \cdots $$
can be considered either as a subspace of $T(V)$ or as a quotient of it.

A linear map $\Vt n \to V$ is called an \emph{$n$-ary operation}\index{$n$-ary operation}  on $V$ 
and a linear map $ V\to \Vt n$ is called an \emph{
$n$-ary cooperation} on $V$. 
The symmetric group is the automorphism group of the finite set $\{1,\ldots , n\}$ and is denoted $S_n$. It acts on $\Vtn$ on the left by $\ss \cdot \row v1n = \row v{\ss^{-1}(1)}{\ss^{-1}(n)}$. The action is extended to an action of $\KK[S_n]$ by linearity.
We denote by $\tau$ the \emph{switching map}\index{switching map} in the symmetric monoidal category $\Vect$, that is $\tau(u\t v) = v\t u$ (in the nongraded case). 

A \emph{magmatic algebra}\index{magmatic algebra} is a vector space $A$ equipped with a binary operation $A\t A \to A$, usually denoted $(a,b)\mapsto a\cdot b$\ . In the unital case it is assumed that there is an element $1$, called the unit, which satisfies $a\cdot 1 = a = 1\cdot a$. In the literature a magmatic algebra is sometimes referred to as a nonassociative algebra.

Quotienting by the associativity relation $(ab)c=a(bc)$ we get the notion of \emph{associative algebra}\index{associative algebra}. Quotienting further by the commutativity relation $ab=ba$ we get the notion of  \emph{commutative algebra}\index{commutative algebra}. So, in the terminology ``commutative algebra", associativity is understood.
\medskip

\noindent {\bf References.} References inside the monograph include the Chapter. So ``see 2.3.4" means see paragraph 3.4 in Chapter 2.


\chapter{Algebraic operads} 
  
We briefly recall the definition, notation and terminology of the operad framework (see for instance  \cite{May, MSS} or the most recent book \cite {LV}).
The reader who is familiar with algebraic operads and props can skip this first chapter and come back to it whenever needed.
 
 \section{\sm}\label{Smod}
 
 \subsection{\sm and Schur functor} An \emph{\sm}\index{\sm} $\PP$ is a family of right $S_n$-modules $\PP(n)$ for $n\geq 0$. Its associated \emph{Schur functor}\index{Schur functor} $\PP : \Vect \to \Vect $ is defined as
 $$\PP (V) := \bigoplus _{n\geq 0 }\PP(n) \t_{S_n} V^{\t n} $$
 where $S_n$ acts on the left on $\Vtn$ by permuting the factors. We also use the notation $\PP(V)_{n} := 
 \PP(n) \t_{S_n} V^{\t n} $ so that $\PP (V) := \bigoplus _{n\geq 0 }\PP(V)_{n} $.
 
 In this monograph we always assume that $\PP(n)$ is finite dimensional. In general we assume that $\PP(0) =0$ and $\PP(1)=\KK\, \id$ (\emph{simply connected operad}). In a few cases we assume instead that $\PP(0)=\KK$, so that  the algebras can be equipped with a unit. The natural projection map which sends $\PP(V)_{n}$ to 0 when $n>1$ and $\PP(V)_{1}$ to
  itself, that is $\KK\ \id \t V=V$ is denoted
  $$\proj : \PP(V) \to V.$$
  
  A morphism of {\bf S}-modules $f : \PP \to \PP' $ is a family of $S_{n}$-morphisms $f(n) : \PP(n) \to \PP'(n) $. They induce a morphism of Schur functors: 
  $$f(V): \PP(V) \to \PP'(V).$$
  
  \subsection{Composition of {\bf S}-modules}\label{compositionSmod} Let $\PP$ and $\QQQ$ be two {\bf S}-modules. It can be shown that the composite $\QQQ \circ \PP$ of the Schur functors (as endomorphisms of $\Vect$) is again the Schur functor of an \sm , also denoted $\QQQ \circ \PP$. The explicit value of $(\QQQ \circ \PP)(n)$ involves sums, tensor products and induced representations of the representations $\QQQ(i)$ and $\PP(i)$ for all $i\leq n$. This composition makes the category of {\bf S}-modules into a monoidal category whose neutral element is the identity functor. Observe that this is not a symmetric monoidal category since the composition of functors is far from being symmetric.

  \subsection{Generating series}\label{genseries} The \emph{generating series}\index{generating series} of an \sm $\PP$ is defined as
  $$f^{\PP}(t)  := \sum_{n\geq 1} \frac{ \dim \PP(n) }{n!} t^n\ .$$
  It is immediate to check that the generating series of a composite is the composite of the generating series:
  $$f^{\QQQ\circ \PP}(t) =f^{\QQQ}(f^{\PP}(t) ).$$
  
\section{Algebraic operad} 

\subsection{Definition}\label{def:operad}\index{operad} By definition an \emph{algebraic operad}\index{algebraic operad}, or \emph{operad} for short, is a Schur functor $\PP$ equipped with two transformations of functors $\iota : \Id_{\Vect} \to  \PP$ and $\cc : \PP\circ \PP \to \PP$ which make it into a monoid. In other words we assume that $\cc$ is associative and that $\iota$ is a unit for $\cc$.
 Such an object is also called a \emph{monad} in $\Vect$.
 
 The identity functor $\Id_{\Vect}$ is itself an operad that we denote by $Vect$ (instead of $\Id_{\Vect}$)\index{Vect@$Vect$} when we consider it as an operad. We call it the \emph{identity operad}.

We usually assume that the operad is \emph{connected}, that is $\PP(0)=0$.

\subsection{Algebra over an operad}\label{algoverop} By definition an \emph{algebra over the operad}\index{algebra over an operad} $\PP$, or a $\PP$-{\it
algebra} for short, is a vector space $A$ equipped with a linear
map $\cc_A : \PP (A) \to A$ such that the following diagrams are
commutative:
$$\xymatrix{
\PP \circ \PP (A) \ar[r]^-{\PP(\cc_A)}
\ar[d]_{\cc(A)}&\PP(A)\ar[d]^{\cc_A}&\qquad &
A \ar[r]^{\iota (A)} \ar[dr]_{=}&\PP(A)\ar[d]^{\cc_A}\\
\PP(A) \ar[r]^{\cc_A} & A&\qquad&&  A }$$ 
There is an obvious notion of morphism of $\PP$-algebras. The category of $\PP$-algebras is denoted $\PP$-alg. Since we made the assumption $\PP(0)=0$, these algebras are nonunital.

If we want a unit for the $\PP$-algebras, then we take $\PP(0)=\KK$ and the image of $1\in \KK$ by $\cc_{A}:\PP(0)\to A$ is the unit of $A$.

The operation $\id\in \PP(1)$ is the identity operation: $\id(a)=a$ for any $a\in A$. For $\mu\in \PP(k)$ and $\mu_1\in \PP(n_1),\ldots , \mu_k\in \PP(n_k)$ the composite $\cc(\mu; \mu_1,\cdots , \mu_k)\in \PP(n_1 +\cdots + n_k)$ is denoted $\mu\circ (\mu_1,\cdots , \mu_k)$ or $\mu (\mu_1,\cdots , \mu_k)$  if no confusion can arise.

$$
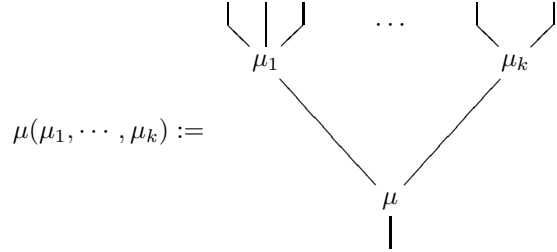
{\xymatrix@R=3pt@C=3pt{
&\ar@{-}[d]&\ar@{-}[dd]&\ar@{-}[d]& &&&& &\ar@{-}[d]&&\ar@{-}[d]& \\
&*{}\ar@{-}[dr]&&*{}\ar@{-}[dl]&                                                 & &\cdots&&& *{}\ar@{-}[dr]&&*{}\ar@{-}[dl]& \\
&&\mu_1\ar@{-}[ddddrrrr]&& &&&& &&\mu_k\ar@{-}[ddddllll]&& &&& \\
&&&& &&&& &&&& \\
\mu (\mu_1,\cdots , \mu_k):=&&&& &&&& &&&& \\
&&&& &&&& &&&& \\
&&&& &&\mu\ar@{-}[dd]&&& &&& \\
&&&& &&&& &&&& \\
&&&& &&&& &&&& \\ }}
$$

So, an operad can also be described as a family of linear maps
$$\cc : \PP(k)\t \PP(i_{1}) \t \cdots \t \PP(i_{k})\to \PP(i_{1}+\cdots + i_{k})$$
which assemble to give a monad $(\PP, \cc, \iota)$.

The restriction of $\cc_{A}$ to $\PP(n)\t_{S_{n}}A^{\t n}$ is denoted $\cc_{n}:\PP(n)\t_{S_{n}}A^{\t n}\to A$ if no confusion can arise.

Given an element $\mu\in \PP(n)$ and an $n$-tuple $(a_{1}, \ldots , a_{n})$ of elements of $A$, we can construct 
$$\mu(a_{1}, \ldots , a_{n}):= \cc_{n}(\mu\t (a_{1}, \ldots , a_{n}))\in A\ .$$
Hence $\PP(n)$ is referred to as the ``space of $n$-ary operations" for $\PP$-algebras. The integer $n$ is called the ``arity"\index{arity} of the operation $\mu$.

The category of algebras over the operad $Vect$ is simply the category of vector spaces $\Vect$. Hence we have $Vect(1)=\KK$ and $Vect(n)=0$ for $n\neq 1$.

\subsection{Free $\PP$-algebra}\label{freePalgebra}\index{free algebra} By definition a $\PP$-algebra $A_0$ is \emph{free} over the vector space $V$ if it is equipped with a linear map $i : V\to A_0$ and if it satisfies the following universal property:

any map $f : V\to A$, where $A$ is a $\PP$-algebra, extends
uniquely into a $\PP$-algebra morphism $\tilde f : A_0 \to A$:
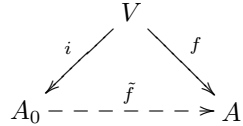
$$\xymatrix{
&V \ar[dl]_{i} \ar[dr]^{f}&\\
A_0 \ar@{-->}[rr]^{\tilde f}&&A }$$ 
Observe that the free algebra over $V$ is well-defined up
to a unique isomorphism. 

For any  vector space $V$ one can equip $\PP(V)$ with a structure of
$\PP$-algebra by setting $\cc_{\PP(V)} := \cc (V) : \PP (\PP (V)) \to \PP (V).$
The axioms defining the operad $\PP$ show that $(\PP(V),\cc (V))$
is  the free $\PP$-algebra over  $V$.

 Categorically, the functor $\PP:  \Vect \to\PP\alg, V\mapsto \PP(V)$ is  left adjoint to the forgetful functor
which assigns, to a $\PP$-algebra $A$, its underlying vector space:
\[ {\Hom}_{\PP\alg}(\PP(V), A) \cong {\Hom}_{\Vect}(V, A).\]

\subsection{Operadic ideal}\label{opideal}\index{operadic ideal} For a given operad $\PP$ and a family of operations $\{\nu\}$ in $\PP$ the ideal $\II$, generated by this family, is the sub-\sm $\II$ linearly generated by all the compositions 
$\mu\circ (\mu_1,\cdots , \mu_k)$ where at least one of the operations is in the family. The quotient $\PP/\II$, defined as $(\PP/\II )(n) = \PP(n)/\II(n)$, is an operad.

If $\QQQ$ is a suboperad of $\PP$, then we denote by $\overline \QQQ$ the sub-\sm of $\QQQ$ such that $\bar \QQQ(1)=\QQQ(1)/\KK \id$ and $\bar \QQQ(n) =  \QQQ(n)$ for $n\geq 2$. We denote by $(\bar \QQQ)$ the operadic ideal generated by $\bar \QQQ$ in $\PP$. So the quotient $\PP/(\bar \QQQ)$ is an operad.

\subsection{Type of algebras and presentation of an operad}\index{type of algebras} For a given type of algebras defined by generators and relations (supposed to be multilinear), the associated operad is obtained as follows. Let $\PP(V)$ be the free algebra of the given type over $V$. Let $V= \KK x_1 \oplus \cdots \oplus  \KK x_n$ be a based $n$-dimensional vector space. The multilinear part of $\PP(V)$ of degree $n$ ( i.e.\ linear in each variable) is a subspace of $\PP(V)_{n}$ denoted $\PP(n)$. It is clear that $\PP(n)$ inherits an action of the symmetric group. The universal property of the free algebra $\PP(V)$ permits us to give a structure of operad on the Schur functor $\PP$. The category of $\PP$-algebras is precisely the category of algebras we started with.

The operad $\PP$ can also be constructed by taking the free operad over the generating operations and quotienting by the ideal (in the operadic sense) generated by the relators.

For instance the free operad on one binary operation $\mu$ (with no symmetry) is the magmatic  operad\index{magmatic} $Mag$. In degree $n$ we get $Mag(n)=\KK[PBT_{n}]\t \KSn$ where $PBT_{n}$ is the set of planar binary rooted trees with $n-1$ internal vertices (and $n$ leaves), cf.~\ref{pbtree}. The tree $|\in PBT_{1}$ codes for  $\id\in Mag(1)$ and the tree  $\arbreA\in PBT_{2}$ codes for the generating operation $\mu\in Mag(2)$. 

The operad $As$ of (nonunital) associative algebras is the quotient of $Mag$ by the ideal generated by the relator 
$$\mu\circ (\mu \t\id ) - \mu\circ (\id \t \mu)\in Mag(3) .$$

Observe that a morphism of operads $\PP\to \QQQ$ gives rise to a functor between the corresponding categories of algebras in the other direction:
$$\QQQ\alg \tto \PP\alg.$$

\subsection{Binary and quadratic operad}\index{binary}\index{quadratic} An element $\mu\in \PP(2)$ defines a map
$$\mu : A^{\t 2} \to A, \quad a\t b \mapsto \mu(a,b)\ ,$$
called a \emph{binary operation}.
Sometimes such an operation is denoted by a symbol, for instance $*$, and we write $a*b$ instead of  $\mu(a,b)$. We allow ourselves to talk about ``the operation $a*b$ ".

An operad is said to be \emph{binary}, resp. $\langle k\rangle$-ary, if it is generated by binary, resp. $\langle k\rangle$-ary, operations (elements in $\PP(2)$, resp.~$\PP(k)$). An operad is said to be \emph{quadratic}\index{quadratic} if the relations are made of monomials involving only the composition of two operations. In the binary case, it means that the relations are of the form
$$ \sum_i \mu_i (\nu_i\t \Id) = \sum_j \mu_j(\Id \t\nu_j)$$
where the elements $\mu_i , \nu_i , \mu_j , \nu_j$  are binary operations (not necessarily the generating ones). Sometimes, in the literature, the adjective quadratic is used in place of binary and quadratic (see for instance \cite{GK}). Most classical types of algebras are defined by binary quadratic operads: associative, commutative, Lie, Poisson, pre-Lie,  Leibniz, dendriform, 2-associative, alternative, magmatic, etc. Some are generated by $n$-ary operations, but are still quadratic: Lie triples, Jordan triples, $A_{\infty}, C_{\infty}, L_{\infty},Brace, MB, \mathcal{M}ag^{\infty}$, etc.

\subsection{Nonsymmetric operad}\label{nsoperad}\index{nonsymmetric operad} Let $\PP$ be an operad whose associated type of algebras has the following property. The generating operations do not satisfy any symmetry property and, in the relations, the variables stay in the same order. Then, it is easy to show that $\PP(n)= \PP_n\t \KSn$ for some vector space $\PP_n$. Here $\KSn$ stands for the regular representation. Moreover the operadic structure is completely determined by composition maps
$$\cc_{i_1\cdots  i_n}:\PP_n \t \PP_{i_1} \t \cdots  \t \PP_{i_n} \longrightarrow \PP_{i_1+ \cdots + i_n}\ .$$
Such operads are called \emph{regular operads}\index{regular operad}. The object $(\PP_{n},\cc_{i_1\cdots  i_n})_{n\geq 1}$ is called a \emph{nonsymmetric operad}\index{nonsymmetric operad}. The operads as defined in \ref{def:operad} are also called \emph{symmetric operad} when there is a risk of confusion. The operads 
$$As, Dend, Dipt, 2as, Mag, Dup, A_{\infty}, Mag^{\infty}$$ are regular and so are determined by nonsymmetric operads that we denote by the same symbol.

\subsection{Set-theoretic operad} So far we have defined an operad in the monoidal category of vector spaces (and tensor product $\t$), but we could choose the category of sets (and cartesian product $\times$). This would give us the notion of \emph{set-operad}\index{set-operad}. Since the functor $\KK[-]: Set \to \Vect$ is strong monoidal, any set-operad gives rise to an algebraic operad. Such an operad is said to be \emph{set-theoretic}\index{set-theoretic}.

\subsection{Classical examples: the three graces} The classical examples of algebraic operads are the operad $As$ of associative algebras, the operad $Com$ of commutative algebras (understood to be associative) and the operad $Lie$ of Lie algebras. In each case the free algebra is well-known, so the operad is easy to describe: $As(V)$ is the (nonunital) algebra of noncommutative polynomials over $V$ (i.e.\ reduced tensor algebra $\To(V)$),  $Com(V)$ is the (nonunital) algebra of  polynomials over $V$ (i.e.\ reduced symmetric algebra $\overline S(V)$),  $Lie(V)$ is the subspace of $As(V)$ generated by $V$ under the bracket operation $[x,y]=xy-yx$. It follows that in the associative case we get $As(n)=\KK[S_n]$ (regular representation).
 In the commutative case we get $Com(n)=\KK$ (trivial representation).  In the Lie case we get $Lie(n)=\Ind^{S_n}_{C_n}(\KK)$ (induced representation from the cyclic representation over the cyclic group $C_n$ when $\KK=\mathbb{C}$). 
 
 Observe that $As$ is nonsymmetric, $As$ and $Com$ are set-theoretic, but  $Lie$ is not regular nor set-theoretic.
 
\section{Coalgebra and cooperad}

 \subsection{Coalgebra over an operad}\label{coalgop}\index{coalgebra} By definition a \emph{coalgebra over the operad} $\CC$ is a vector space $C$ equipped with $S_n$-equivariant maps 
$$\CC(n)\t C \to C^{\t n}, \quad \dd\t c \mapsto \dd(c)= c_{1}^{\dd}\t \cdots \t c_{n}^{\dd}$$
 (we omit the summation symbol on the right hand side) which are compatible with the operad structure of $\CC$
 . In particular, there is a commutative diagram
 $$\xymatrix{
 \CC(k)\t \CC(i_{1})\t \cdots \t \CC(i_{k})\t C\ar[rr]\ar[d]&  &\CC(i_{1})\t  \cdots \t \CC(i_{k})\t C^{\t k} \ar@{=}[d]\\
\CC(i_{1}+\cdots +i_{k})\t C \ar@{=}[d]         &  & \CC(i_{1})\t C\t  \cdots \t \CC(i_{k})\t C\ar[d]\\
 \CC(n)\t C \ar[rr]                                                         &   &C^{\t n} =  C^{\t i_{1}}\t \cdots \t C^{\t i_{k}}\\
 }$$
 
 In this framework the elements of $\CC(n)$ are called \emph{$n$-ary cooperations}.

\subsection{Primitive part, connectedness (conilpotency)}\label{connectedness}\index{connected}\index{primitive}\index{conilpotent} Let $\CC$ be an operad such that $\CC(0)=0$ and $\CC(1)=\KK\, \id$. We suppose that there is only a finite number of generating cooperations in each arity so that $\CC(n)$ is finite dimensional.  The identity operation $\id$ is not considered as a generating cooperation. Let $C$ be a $\CC$-coalgebra. We define a filtration on $C$ as follows:
$$F_1C= \Prim C := \{x\in C \ \vert \ \dd(x) = 0 \textrm{ for any generating cooperation } \dd \}.$$
The space $\Prim C$ is called the \emph{primitive part} of $C$, and its elements are said to be \emph{primitive}. Then we define the filtration by:
$$F_rC := \{x\in C \ \vert \ \dd(x)=0 \textrm{ for any } \dd\in \CC(n), n>r \}.$$
We adopt Quillen's terminology (cf.\  \cite{Quillen} Appendix B) and say that the coalgebra $C$ is \emph{connected}, or \emph{conilpotent},  if $C= \bigcup_{r\geq 1}F_{r}C$.

If the operad $\CC$ is binary, then this definition of connectedness is equivalent to the definition which uses the filtration 
$$F'_rC := \{x\in C \ \vert \ \dd(x)\in (F'_{r-1} C)^{\t |\dd |} \textrm{ for any generating cooperation } \dd \}.$$
If $C$ is a connected graded coalgebra, that is $C=\oplus_{n\geq 1}C_{n}$, then $C$ is connected in the aforementioned sense.

\subsection{Cofree coalgebra}\label{cofree}\index{cofree} By definition a $\CC$-coalgebra $C_{0}$ is said to be \emph{cofree} over the vector space  $V$ if it is connected, equipped with a map $s:C_{0} \to V$ and if it satisfies the following universal property:

any map $p: C\to V$, where $C$ is a connected $\CC$-coalgebra, extends uniquely into a $\CC$-coalgebra morphism  $\tilde p : C \to C_{0}$:
$$\xymatrix{
C \ar@{-->}[rr]^{\tilde p} \ar[dr]_{p}&&C_{0}\ar[dl]^{s}\\
&V &
 }$$ 
The cofree coalgebra over $V$ is well-defined up
to a unique isomorphism. Observe that we are working within the category of ``connected" coalgebras. If we were working in the whole category of coalgebras, the notion of cofree object would be different.

Let $\CC=As$. The cofree coassociative coalgebra over $V$ is the reduced tensor module $\To(V)$ equipped with the \emph{deconcatenation}\index{deconcatenation} operation:
$$\dd(v_{1}\ldots v_{n}) = \sum_{1\leq i \leq n-1} v_{1}\ldots v_{i}\t v_{i+1}\ldots v_{n}\ .$$

\subsection{Cooperad and coalgebra over a cooperad}\label{cooperad}\index{cooperad} Taking the linear dual of an operad $\CC$ gives a \emph{cooperad} denoted $\CC^c$. Let us recall that a cooperad is a comonoid structure on a Schur functor. As a vector space $\CC^c(n)= \CC(n)^*=\Hom (\CC(n), \KK)$. We equip this space with the following right $S_n$-module structure: 
$$f^{\ss} (\mu) := f(\mu^{\ss^{-1}}),$$
for $f\in \CC(n)^*, \mu\in \CC(n)$ and $\ss\in S_n$.


The cooperadic composition is denoted $\theta: \CC^c \to \CC^c\circ \CC^c$.  There is an obvious notion of \emph{coalgebra} $C$ over a cooperad $\CC^c$ given by maps:
$$\theta(C) : C \tto \CC^c(C)^{\wedge} := \prod_{{n\geq 1} } \CC^c(n) \t_{S_n} C^{\t n}.$$

It coincides with the notion of coalgebra over an operad when the cooperad is the linear dual of the operad. Here we use the characteristic zero hypothesis to identify invariants and coinvariants under the symmetric group action. We always assume that $\CC^c(0)=0$, $\CC^c(1)=\KK\id$ and that $\CC(n)$ is finite dimensional. The elements of $\CC^c(1)$ are called the trivial cooperations and an element $f\in\CC^c(n), n\geq 2$, is called a \emph{nontrivial cooperation}\index{nontrivial cooperation}.


The projection of $\theta(C) $ to the $n$-th component is denoted 
$$\theta_{n} :C \to \CC^c(n) \t_{S_n} C^{\t n}$$
 if no confusion can arise. Let $\langle-, -\rangle : \CC(n) \t \CC^c(n) \to \KK$ be the evaluation pairing. The relationship with the notation introduced in \ref{coalgop} is given by the commutative diagram
 
 $$\xymatrix{
        &\CC(n) \t \CC^c(n)\t_{S_{n}} C^{\t n }\ar[dr]^{\langle-, -\rangle\t \id} & \\
\CC(n)\t C \ar[ru]^{\id\t \theta_{n}}\ar [rr] && C^{\t n }\\
\dd\t c &\mapsto &\dd(c)=c_{1}^{\dd}\t \cdots \t c_{n}^{\dd}\\
}$$

\begin{lemma} If the coalgebra $C$ is connected (cf.~\ref{connectedness}), then the map $\theta(C)$ factors through the direct sum:
$$\theta(C) : C \tto \CC^c(C) := \bigoplus_{{n\geq 1} } \CC^c(n) \t_{S_n} C^{\t n}.$$
\end{lemma}
\begin{proo} Since for any $x\in C$ there is an integer $r$ such that $\dd(x)=0$ for any cooperation $\dd$ such that $|\dd | > r$, it follows that there is only finitely many nonzero components in $\theta(C)(x)$.
\end{proo}

\subsection{Cofree coalgebra and cooperad}\label{cofreecoalg}
From the axioms of a cooperad it follows that $\CC^c(V)$ is the cofree $\CC^c$-coalgebra over $V$.

Explicitly, for any connected coalgebra $C$,  the universal lifting $\tilde p : C \to \CC^c(V)$ induced by a linear map $p:C\to V$ is obtained as the composite 
$$C\stackrel{\theta(C)}{\tto}\CC^c(C) \stackrel{\CC^c(p)}{\tto}\CC^c(V).$$

Suppose that  the operad $\CC$ is binary quadratic, generated by operations $\mu_1, \mu_2 , \ldots$ and relations of the form
$$ \sum_{i,j} \aa_{ij} \mu_i (\mu_j\t \id) = \sum \bb_{ij} \mu_i(\id \t\mu_j), \quad  \aa_{ij}, \bb_{ij}\in \KK .$$
Then a coalgebra $C$ over $\CC$ is defined by the cooperations $\mu_i^* : C \to C\t C$ satisfying the relations:
$$ \sum_{i,j} \aa_{ij} (\mu_j^*\t \id) \mu_i^* = \sum_{i,j} \bb_{ij} (\id \t\mu_j^*)\mu_i^*.$$

\subsection{Invariants versus coinvariants} Saying that a binary cooperation $\dd :C\to C\t C$ is symmetric means that its image lies in the invariant subspace $(C\t C)^{S_{2}}$. In characteristic zero the natural map from invariants to coinvariants is an isomorphism $(C\t C)^{S_{2}}\cong (C\t C)_{S_{2}}$. Therefore, if $\CC(2)= \KK$ (trivial representation), then $\dd$ defines a map 
$C\to \CC(2)\t_{S_{2}}C^{\t 2}$.
\subsection{Nonsymmetric cooperad}
If the operad $\CC$ is nonsymmetric (cf. \ref{nsoperad}), then the equivalence between the two notions of coalgebra does not need the characteristic zero hypothesis. Indeed, since $\CC(n)=\CC_{n}\t \KSn$, we simply take $\CC^c(n):=\CC^*_{n}\t\KSn$.

\section{Prop}\label{prop}
A type of algebras is governed by an operad. Similarly there is an algebraic device which governs a type of bialgebras: it is called a prop.

\subsection{Definition} In the algebra framework (i.e.\ operad framework) an operation $\mu$ can be seen as a box with $m$ inputs and one output:

$$\boitenun$$

In the coalgebra framework (i.e.\ cooperad framework) a cooperation $\dd$ can be seen as a box with $1$ input and $n$ outputs:

$$\boiteunm$$

If we want to deal with  bialgebras, then we need boxes with multiple inputs and multiple outputs:

$$\boitenm$$

\noindent called \emph{multivalued operations} or {bioperations}.\index{multivalued operation}
Hence we have to replace the {\bf S}-modules by the \bisms,  i.e.\ families $\PP(m,n)$ of $S_{m}^{op}\times S_{n}$-bimodules. A composition of multivalued operations is prescribed by a bipartite graph. This composition is supposed to be compatible with the symmetric group actions and to satisfy some obvious associativity and unitality axioms. The whole structure is called a \emph{prop} (also denoted previously PROP), cf.\ \cite{McLane, FM, Vallette04}. An algebra (or better a \emph{gebra}) over a prop is a vector space $\HH$ equipped with maps
$$\PP(m,n) \t \HH^{\t m} \tto  \HH^{\t n}$$
which are compatible with the symmetric group actions and with the composition in the prop.

A type of gebra can be defined by generators and relations. When the generators are either operations (i.e.\ elements in $\PP(m,1)$) and/or cooperations (i.e.\ elements in $\PP(1,n)$), the gebras are called \emph{generalized bialgebras}, or simply bialgebras if there is no ambiguity with the classical notion of bialgebras.
An explicit description of the props corresponding to the classical bialgebras can be found in \cite {Pira02}.

Given a generating set of multivalued operations one can construct the free prop over them (cf.\ for instance \cite{FM} and \cite {Vallette04} section2). Moding out by  relations gives rise to a quotient prop. 

Unlike the case of algebras and operads, this prop cannot be interpreted as a ``free gebra'' because the forgetful functor from gebras (of a given type) to vector spaces does not admit a left adjoint in general.


 \chapter{Generalized bialgebra and triple of operads}\label{Ch:genebialg}
 
 We introduce the notion of \emph{generalized bialgebra} and its primitive part.  We denote by $\CC$ the operad driving the coalgebra structure and by $\AA$ the operad driving the algebra structure. We prove that, under some hypotheses {\texttt{(H0)}} and \Hone, a generalized bialgebra type determines an operad called the \emph{primitive operad}. 
 Primitive elements in a generalized bialgebra do not, in general, give a primitive element under an operation. However they do when this operation is primitive.  So, the primitive part of a generalized bialgebra has the property of being an algebra over the primitive operad. 
 
 Then we treat the case where the primitive operad is trivial (i.e.\ $\Vect$). We prove that under the hypotheses {\texttt{(H0)}} (distributive compatibility condition), \Hone (the free algebra is a bialgebra), and \Htwoiso (free isomorphic to cofree), any connected \cab is both free and cofree. This is the \emph{rigidity theorem} \ref{thm:rigidity}. The key of the proof is the construction of a universal idempotent $e_{\HH} : \HH \to \HH$ whose image is the space of primitive elements $\Prim \HH$. 
 
 For a given prop $\CC^c\--\AA$ (satisfying {\texttt{(H0)}} and \Hone) whose primitive operad is $\PP$, we call \CAP a \emph{triple of operads}. Our aim is to find simple conditions under which the ``structure theorem" holds for \CAP. This structure theorem says that any connected $\CC^c\--\AA$-bialgebra is isomorphic to $U(\Prim \HH)$ as an algebra and is cofree over $\Prim \HH$ as a coalgebra. These simple conditions are {\texttt{(H0)}}, \Hone and \Htwoepi . The latter condition says that the coalgebra map $\AA(V)\to \CC^c(V)$ is surjective and admits a coalgebra splitting.
 
 Then we give some immediate consequences of the main theorem.
 
 In this chapter we suppose that the ground field is of characteristic zero. We indicate in the next chapter how to avoid this hypothesis in certain cases.

\section{Generalized bialgebra}\label{genbialg} 

We consider a certain type of prop generated by operations and cooperations. A gebra over such prop is called a generalized bialgebra.

  \subsection{Compatibility relation and generalized bialgebra}\label{comprel} Let $\AA$ and $\CC$ be two algebraic operads. We always assume that there is a finite number of generating operations in each arity. As a consequence $\CC(n)$ and $\AA(n)$ are  finite dimensional vector spaces.
  
  By definition  a \emph{$( \CC^c, \between, \AA)$-bialgebra}, or \emph{$ \CC^c\pt\AA$-bialgebra} for short, also called \emph{generalized bialgebra}\index{generalized bialgebra},  is a vector space $\HH$ which is an $\AA$-algebra, a $\CC$-coalgebra, and such that the operations of $\AA$ and the cooperations of $\CC$ acting on $\HH$ satisfy  some \emph {compatibility relations}\index{compatibility relation}, denoted $\between$, read ``between" (some equalities involving composition of operations and cooperations valid for any elements of $\HH$). This set of relations is, of course, part of the structure. A category of generalized bialgebras is governed by an \emph{algebraic prop} (we simply say a prop) as mentioned in \ref{prop} (cf.~for instance \cite{Vallette04}). Starting with any presentation of the operad $\AA$ and of the cooperad $\CC$, this prop is obtained as the quotient of the free prop generated by the generators of $\AA$ and $\CC$ (considered as multivalued operations), modulo the relations  between the operations, the relations between the cooperations and the relations entwining operations and cooperations. The gebras over this quotient are the generalized bialgebras.

A \emph{distributive compatibility relation}\index{distributive compatibility relation} between the operation $\mu$  and the cooperation $\dd$ is a relation of the form
 $$\dd\circ \mu = \sum_{i} (\mu_1^i\t\cdots \t  \mu_m^i)\circ\omega^{i} \circ (\dd_1^{i}\t \cdots \t \dd_n^{i})
\qquad \qquad ( \between)$$
 where
 \begin{displaymath}
\left\{\begin{array}{l}
\mu\in \AA(n),\  \mu_1^i\in \AA(k_{1}), \ldots , \mu_m^i \in \AA(k_{m}), \\
\dd\in \CC(m),\  \dd_1^i\in \CC(l_{1}), \ldots , \dd_n^i \in \CC(l_{n}), \\
k_{1}+\cdots + k_{m}= l_{1}+\cdots + l_{n}= r_{i} , \\
\omega^{i} \in \KK[S_{r_{i}}].
\end{array}
\right.
\end{displaymath}

Hence, in a generalized bialgebra, the composite of an operation and a cooperation can be re-written as cooperations first and then operations. Observe that the identity is both an operation and a cooperation.

\subsection{Hypothesis {\texttt{(H0)}}}\index{{\texttt{(H0)}}}  There is  a distributive compatibility relation for any pair $(\dd,\mu)$ where $\mu$ is an operation and $\dd$ is a cooperation.

\medskip

Of course, it suffices to check this hypothesis for $\mu$ a generating operation and $\dd$ a generating cooperation.

For a given relation $\between$ we denote by $\Phi$ the right-hand side term.

The distributive compatibility relations induce a \emph{mixed distributive law} in the sense of Fox and Markl \cite {FM}, that is a map
$$\PP(m)\t \CC(n) \to \bigoplus \CC(i_{1})\t \cdots \t \CC(i_{m})\t _{S_{i}}\KK[S_{N}]\t_{S_{j}}\PP(j_{1})\t \cdots \t \PP(j_{n}).$$
Here we used the multi-index notation for $i=(i_{1}, \ldots , i_{m})$ and for $j$, and $S_{i}:= S_{i_{1}}\times\cdots \times S_{i_{m}}$, $N= i_{1}+ \cdots + i_{m}= j_{1}+\cdots +j_{n}$.

\subsection{Diagrams}\label{diagram} It will prove helpful to write the compatibility relations as diagrams instead of long algebraic expressions. For instance, for a binary operation $\mu $ and  a  binary cooperation $\dd $ we draw

$\vcenter
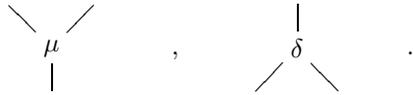
{\xymatrix@R=2pt@C=2pt{
\ar@{-}[ddrr]&&&&\ar@{-}[ddll]\\
&&&&\\
&&\mu\ar@{-}[dd]&&\\
&&&&\\
&&&&
}}$
\qquad , \qquad $\vcenter{\xymatrix@R=2pt@C=2pt{
&&&&\\
&&&&\\
&&\dd\ar@{-}[uu]&&\\
&&&&\\
\ar@{-}[uurr]&&&&\ar@{-}[uull]
}}$\qquad  .

\medskip

The associativity property of $\mu$, which is written $\mu(\mu(x,y),z)= \mu(x,\mu(y,z))$ algebraically, becomes 

\medskip

$\vcenter
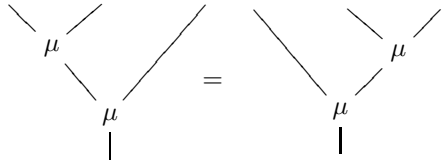
{\xymatrix@R=2pt@C=2pt{
\ar@{-}[ddrr]&&&&\ar@{-}[ddll]&&&&\\
&&&& &&&&\\
&&\mu\ar@{-}[ddrr]&& &&&&\\
&&&& &&&& =\\
&&&&\mu\ar@{-}[dd]\ar@{-}[uuuurrrr] &&&&\\
&&&& &&&&\\
&&&& &&&&
}}$
$\vcenter{\xymatrix@R=2pt@C=2pt{
\ar@{-}[ddddrrrr]&&&& &&&&\\
&&&& &&&&\\
&&&& &&\mu\ar@{-}[uurr]\ar@{-}[uull]&&\\
&&&& &&&&\\
&&&&\mu\ar@{-}[dd]\ar@{-}[uurr] &&&&\\
&&&& &&&&\\
&&&& &&&&
}}$

\noindent pictorially.

Example of a compatibility relation  for the pair $(\dd,\mu)$ with $n=3, m=4$ and $r=8$~:

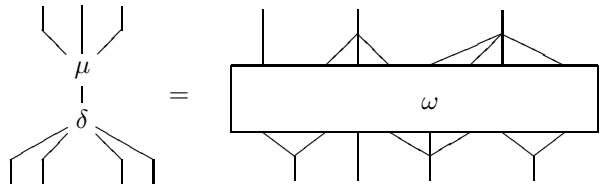
$$\cpbtroisquatre =\quad  \cpbmultiple$$

Here we have $l_{1}=1, l_{2}=3, l_{3}=4 ; k_{1}=2, k_{2}=1, k_{3}=3, k_{4}=2$ and so 
$r=1+3+4=2+1+3+2$. 
Observe that, in the general case, the right-hand side term $\Phi$ is a sum of such compositions. We split $\Phi$ into two summands
$$\Phi = \Phi_{1} + \Phi_{2},$$
as follows. The summand
  $\Phi_{1}$ contains all the terms for which $r=n$ and  $\Phi_{2}$ contains all the terms for which $r> n$, see section \ref{comprel} for the meaning of $r$ and $n$. There is no term with $r<n$ since we assume that $\CC(0)=0$. The important point of this splitting is the following: for each summand 
$$ (\mu_1^i\t\cdots \t  \mu_m^i)\circ\omega^{i} \circ (\dd_1^{i}\t \cdots \t \dd_n^{i})$$
of  $\Phi_{2}$ at least one of the cooperations $\dd_{k}^{i}$ is  nontrivial ( i.e.\ of arity $\geq 2$). In $\Phi_{1}$ the only cooperation which pops up is the identity.

When both operads $\AA$ and $\CC$ are nonsymmetric and, in the compatibility relations, there is no crossing (in particular the only permutations $\oo$ are the identity), then we say that this is a \emph{nonsymmetric case} and that $\CC^c\pt \AA$ is a \emph{nonsymmetric prop}.
\medskip

\subsection{Examples of distributive compatibility relations}
 {}
 \smallskip
 
\subsubsection{Example 1}\label{ex1} Hopf algebra (classical bialgebra, Hopf relation):

A classical bialgebra is a unital associative algebra equipped with a counital coproduct $\DD$ which satisfies the \emph{Hopf compatibility relation}\index{Hopf compatibility relation}
$$\DD(xy) =\DD(x)\DD(y).$$
Since, here, we want to work without unit nor co-unit, we work over the augmentation ideal and with the \emph{reduced comultiplication}\index{reduced comultiplication} $\dd$ defined by
$$\dd(x):= \DD(x) -x\t 1 - 1 \t x.$$

The classical Hopf compatibility relation becomes $\between_{Hopf}$ :

\medskip

\noindent $\dd(xy) =$
$$ x\t y + y\t x  + x_{(1)}\t x_{(2)}y + x_{(1)}y\t x_{(2)}  + xy_{(1)}\t y_{(2)}+ y_{(1)}\t xy_{(2)}+ x_{(1)}y_{(1)}\t x_{(2)}y_{(2)}$$

\noindent under the notation $\dd(x):= x_{(1)}\t x_{(2)}$ (Sweedler notation with summation sign understood).

Pictorially the relation  $\between_{Hopf}$ reads:

$$\cpbdeuxdeux = \underbrace{\cpbA +\cpbB}_{\Phi_{1}} +\underbrace{\cpbC +\cpbD +\cpbE +\cpbF +\cpbG}_{\Phi_2}$$

\medskip

\subsubsection{Example 2}\label{ex2} Nonunital infinitesimal bialgebra ($As^c\pt As$-bialgebra, n.u.i. relation).

The motivation for this case is the tensor algebra that we equip with the deconcatenation coproduct (instead of the shuffle coproduct). In the nonunital framework the compatibility relation satisfied by the concatenation product and the (reduced) deconcatenation coproduct  is

$$\dd(xy)= x\t y  + x_{(1)}y\t x_{(2)} + xy_{(1)}\t y_{(2)}$$

under the notation $\dd(x):= x_{(1)}\t x_{(2)}$.

Pictorially we get $\between_{nui}$ :

$$\cpbdeuxdeux = \underbrace{\cpbA }_{\Phi_{1}}+\underbrace{\cpbC +\cpbE}_{\Phi_{2}}$$

See \cite{LRstr}, where a unital version is handled, and \ref{uib} for more details. The prop defined by this type of generalized bialgebras is nonsymmetric.

\medskip

\subsubsection{Example 3}\label{ex3} Bimagmatic bialgebra ($Mag^c\pt Mag$-bialgebra, magmatic relation). The motivation  for this case is the magmatic algebra $Mag(V)$. We equip it with the magmatic coproduct obtained by identifying the classical basis of $Mag(\KK)$ (planar binary trees) with its dual. In the nonunital framework the compatibility relation is $\between_{mag}$ :

$$\cpbdeuxdeux =   \underbrace{\cpbA}_{\Phi_{1}} +  \underbrace{{\vcenter{\xymatrix@R=2pt@C=2pt{
\\
\\
0\\
\\
\\
}}}}_{\Phi_{2}} $$

This is a nonsymmetric prop.

\medskip

\subsubsection{Example 4}\label{ex4} Frobenius algebra ($As^c\pt As$-bialgebra, $Com^c\pt Com$-bialgebra, Frobenius relation):

$$\cpbdeuxdeux = \underbrace{{\vcenter{\xymatrix@R=2pt@C=2pt{
\\
\\
0\\
\\
\\
}}}}_{\Phi_{1}}  + \underbrace{\cpbC}_{\Phi_{2}}   \qquad\mathrm{and} \qquad  \cpbdeuxdeux = \underbrace{{\vcenter{\xymatrix@R=2pt@C=2pt{
\\
\\
0\\
\\
\\
}}}}_{\Phi_{1}}  + \underbrace{\cpbE}_{\Phi_{2}} $$
This is a nonsymmetric case if the algebra and the coalgebra are not supposed to be commutative.

 These examples and many more will be treated in Chapters \ref{ch:examples} and \ref{ch:duplicial}. See \ref{tableauCR} for a list of some compatibility relations.
 
 \section{The primitive operad}\index{primitive operad}
 
 The primitive part of a generalized bialgebra is, in general, not stable under the operations of the operad $\AA$. However it may be stable under \emph{some} operations. In this section we describe the maximal suboperad $\PP$ of $\AA$ such that the primitive part $\Prim \HH$ of the \cab $\HH$ is a $\PP$-algebra. Both operads $\AA$ and $\CC$ are supposed to be finitely generated and simply-connected, that is $\AA(0)=0=\CC(0)$ and $\AA(1)=\KK\, \id =\CC(1)$.
 \subsection{The hypothesis {\texttt{(H0)}}} Given a type of bialgebras, for instance a set of generating operations, a set of generating cooperations and a set of relations, it may happen that the only bialgebra that can exist is $0$. See for instance the discussion in \cite{FM}, section 11. The following hypothesis asserts that this is not the case, that is, the prop is nontrivial.
 
\noindent  \Hone\index{\Hone} The free $\AA$-algebra $\AA(V)$ is  equipped with a \cab structure which is functorial in $V$.
 
 One could translate this hypothesis as a condition on the associated prop.
 
  \subsection{The primitive part of a bialgebra} Let   $\HH$ be a $\CC^c\pt \AA$-bialgebra. By definition the \emph{primitive part} of the $\CC^c\pt \AA$-bialgebra $\HH$, denoted $\Prim \HH$, is
 $$\Prim \HH := \{ x\in \HH \, | \,  \dd (x) = 0 \textrm{ for all } \dd \in \CC^c(n), n\geq 2\}.$$
 Hence, if $\CC$ is generated by $\dd_{1},\ldots , \dd_{k}, \ldots $, then we have 
$$\Prim \HH= \Ker \dd_{1} \cap \ldots\cap \Ker \dd_{k}\cap \ldots\ .$$
 Let us suppose that the hypotheses {\texttt{(H0)}} and \Hone are true.
  By definition an element $\mu\in \AA(n)$ is called a \emph{primitive operation} if,  for any independent variables $x_1,\ldots , x_n$, the element $\mu(x_1,\ldots , x_n)\in \AA(\KK x_1\oplus \cdots \oplus \KK x_n)$ is primitive. In terms of the prop $\CC^c\pt \AA$, $\mu$ being primitive is equivalent to the following: for any generating cooperation $\dd$ the $\Phi_{1}$-part of the compatibility relation for $(\dd, \mu)$ is $0$.

  Let $(\Prim_{\CC}\AA)(n)\subset \AA(n)$ be the \emph{space of primitive operations} for $n\geq 1$:
  $$(\Prim_{\CC}\AA)(n) := \{\mu\in \AA(n) \ | \ \mu \textrm{ is primitive}\}.$$
   By functoriality of the hypothesis, $(\Prim_{\CC}\AA)(n)$ is a sub-$S_n$-module of $\AA(n)$ and so we obtain an inclusion of Schur functors 
$$\Prim_{\CC} \AA \mono \AA.$$
As a result we have $(\Prim_{\CC}\AA)(V) = \Prim (\AA(V))$.

  \begin{thm}[The primitive operad]\label{thm:primitive} Let $(\CC^c,\between ,\AA)$ be a type of generalized bialgebras
 over a  characteristic zero
  field $\KK$. We suppose that the following hypotheses are fulfilled:\\
  
  \noindent {\texttt{(H0)}} any pair $(\dd,\mu)$ satisfies a distributive compatibility relation,\\
\noindent \Hone\index{\Hone} the free $\AA$-algebra $\AA(V)$ is  equipped with a \cab structure which is functorial in $V$.\\
 Then the Schur functor $\PP$, given by $\PP(V):= (\Prim_{\CC} \AA)(V)$, is a suboperad of $\AA$. For any \cab $\HH$ the space $\Prim \HH$ is a $\PP$-algebra and the inclusion $\Prim \HH \mono \HH$ is a morphism of $\PP$-algebras.
 \end{thm}

 \begin{proo} First we remark that the elements of $V\subset \AA(V)$ are primitive, hence $V\subset \PP(V)$,  and $\id$ is a primitive operation. Indeed, since the bialgebra structure of $\AA(V)$ is functorial in $V$, any cooperation $\delta$ on $\AA(V)$ respects the degree. For $n\geq 2$ the degree one part of $\AA(V)^{\t n}$ is trivial. Since $V$ is of degree one in $\AA(V)$, we get $\dd(V)=0$. Hence any element of $V$ is primitive and the functor $\iota:\Id \to \AA$ factors through $\PP$.
 
 To prove that the Schur functor $\PP:\Vect \to \Vect$ is an operad, it suffices to show that it inherits a monoid structure $\PP \circ \PP \to \PP$ from the monoid structure of $\AA$. In other words it suffices to show that composition of primitive operations, under the composition in $\AA$, provides a primitive operation:
 $$\xymatrix{
\PP\circ \PP \ar@{>->}[d]  \ar@{.>}[r] & \PP\ar@{>->}[d]\\
\AA \circ \AA  \ar[r]^-{\cc} &\AA\\
}$$ 
 
 We use the hypothesis of distributivity of the compatibility relation between operations and cooperations, cf. \ref{genbialg}. 
  
 Let $\mu, \mu_1,\ldots , \mu_n$ be operations, where $\mu\in \PP(n)$. We want to prove that the composite $\mu\circ ( \mu_1,\ldots , \mu_n)$ is primitive when all the operations are primitive. It suffices to show that $\dd \circ \mu\circ  ( \mu_1,\ldots , \mu_n)$ applied to the generic element $\row x1s$ is 0 for any  nontrivial cooperation $\dd$.
  
  By \ref{genbialg} we know that $\dd \circ\mu  = \Phi_1 + \Phi_2$, where $\Phi_1$ involves only operations, and  $\Phi_2$ is of the form
  $$ \Phi_2= \sum_{i} (\mu_1^i\t\cdots \t  \mu_m^i)\circ\omega^{i} \circ (\dd_1^{i}\t \cdots \t \dd_n^{i})$$
   where, for any $i$,  at least one of the cooperations $\dd_k^{i},\  k=1,\ldots, n$, is nontrivial.
   We evaluate this expression on a generic element $\row x1n$. On the left-hand side $\mu\row x1n$ is primitive by hypothesis, so 
 $( \dd \circ\mu )\row x1n=0$. On the right-hand side $\Phi_2 \row x1n=0$ because the evaluation of a nontrivial cooperation on a generic element (which is primitive) is 0. Hence we deduce that $\Phi_1\row x1n=0$. Therefore the operation 
$\Phi_1$ is $0$. 

 Let us now suppose that, not only $\mu$ is primitive, but $\row {\mu}1n$ are also primitive operations. By the preceding argument we get
 $$\Phi_2\circ \row {\mu}1n = \sum_{i} \Psi\circ (\dd_1^{i}\mu_{1}\t \cdots \t \dd_n^{i}\mu_{n}) $$
where,  for any $i$,  at least one of the cooperations $\dd_k^{i},\  k=1,\ldots, n$ is nontrivial.
 Hence, summarizing our arguments, the evaluation on $\row x1s$ gives:
 
 \begin{displaymath}
\begin{array}{rcl}
 ( \dd \circ \mu \circ \row {\mu}1k )\row x1s &=& \Phi_2\circ \row {\mu}1n \row x1s \\
 &=&  \Psi\circ 
 (\dd_1^{i}\mu_{1}^{i}(x_1\ldots ),\ldots , \dd_n^{i}\mu_{n}^{i}(\ldots x_s))\\
 &=& 0 ,
\end{array}
\end{displaymath}

\noindent because  a nontrivial cooperation applied to a primitive element gives 0.

In conclusion we have shown that, when $\mu, \mu_1, \ldots, \mu_n$ are primitive, then
$$\dd \circ \mu \circ \row {\mu}1n  \row x1s =0 .$$ 
Hence the operation $ \mu \circ \row {\mu}1n$ is primitive. As a consequence the image of the composite 
$$\PP\circ \PP \mono \AA\circ \AA \stackrel{\cc}{\to} \AA$$
lies in $\PP$ as expected, and so $\PP$ is a suboperad of $\AA$.

\medskip

From the definition of the primitive part of the \cab $\HH$ it follows that $\Prim \HH$ is a $\PP$-algebra. Since $\PP$ is a suboperad of $\AA$, $\HH$ is also a $\PP$-algebra and the inclusion $\Prim \HH \to \HH$ is a $\PP$-algebra morphism.
 \end{proo}
 
\subsection{Examples} Theorem \ref{thm:primitive} proves the existence of an operad structure on $\PP = \Prim_{\CC} \AA$, however, even when $\AA$ and $\CC$ are described by generators and relations, it is often a challenge to find a small presentation of $\PP$ and then to find explicit formulas for the functor $F: \AA\alg \to \PP\alg$.

 In the case of the classical bialgebras, the primitive operad is $Lie$ and the functor $F:As\alg \to Lie\alg$ is the classical Liezation functor: $F(A)$ is $A$ as a vector space and the bracket operation is given by $[x,y]=xy-yx$, cf.~\ref{ComAsLie}. 

In the case of u.i. bialgebras the primitive operad is $Vect$ and the functor $F$ is simply the forgetful functor, cf.~\ref{AsAsVect}.  

In the magmatic bialgebras case  the primitive operad is $Vect$ and the functor $F$ is simply the forgetful functor, cf.~\ref{bimag}.  

In the case of Frobenius bialgebras the primitive operad is $As$ (that is the whole operad) and the functor $F$ is the identity.

We end this section with a result which will prove helpful in the sequel.

\begin{lemma}\label{fi} Let \ca be a generalized bialgebra type verifying the hypotheses {\texttt{(H0)}} and \Hone  of Theorem \ref{thm:primitive}. Let $\varphi(V): \AA(V) \to \CC^c(V)$ be the unique coalgebra map induced by the projection map ${\proj} : \AA(V) \to V$. Denote by
$$\langle -, - \rangle : \CC(n) \times \CC^c(n) \to \KK$$
the pairing between the operad and the cooperad. 

For any  cooperation $\dd\in \CC(n)$ the image of $( \mu; x_{1}\cdots x_{n})\in \AA(V)$ is
$$\dd( \mu; x_{1}\cdots x_{n})= \langle \dd, \varphi_n (\mu) \rangle  x_{1}\t \cdots \t x_{n}  \in \Vtn\subset \AA(V)^{\t n}.$$
\end{lemma}
\begin{proo} Let us recall from \ref{cooperad} that the map $\varphi_n :\AA(n)\to  \CC^c(n)$ is given by the composite

$$\AA(V)\stackrel{\theta_{n}}{\tto}\CC^c(n)\t_{S_{n}}\AA(V)^{\t n} \xrightarrow{\Id\t{\proj}^{\t n}}\CC^c(n)\t_{S_{n}}V^{\t n}.$$

By assumption the bialgebra structure of $\AA(V)$ is functorial in $V$. Therefore $\theta_{n}(\mu\t\row x1n)$ is linear in each variable $x_{i}$. Hence it lies in $\CC^c(n)\t_{S_{n}}\AA(V)_{1}{}^{\t n}=\CC^c(n)\t_{S_{n}}V^{\t n}$. So we have proved that $\theta_{n}(\mu\t\row x1n)= \varphi_n (\mu)\t \row x1n$.

By definition, the coalgebra structure of $\AA(V)$
$$\CC(n)\t \AA(V) \to \AA(V)^{\t n}$$
 is dual (cf.~\ref{cooperad}) to
 $$\theta_{n}:\AA(V) \to \CC^c(n)\t  \AA(V)^{\t n}$$
 via the pairing $\langle-,-\rangle$.
 Hence we get
 $$\dd(\mu\t\row x1n)=  \langle \dd, \varphi_n (\mu) \rangle x_{1}\t \cdots \t x_{n}.$$
\end{proo}

\section{Rigidity theorem}\label{s:rigidity} We first study the generalized bialgebra types for which the primitive operad is trivial. The paradigm is  the case of cocommutative commutative bialgebras (over a characteristic zero field). The classical theorem of Hopf and Borel \cite{Borel}, can be phrased as follows:\\

\noindent {\bf Theorem} (Hopf-Borel)\label{Hopf-Borel}\index{Hopf-Borel theorem}. \emph{In characteristic zero any connected cocommutative commutative bialgebra is both free and cofree over its primitive part}.\\

In other words such a bialgebra $\HH$ is isomorphic to $S(\Prim \HH)$ (symmetric algebra over the primitive part), see \ref{comcom} for more details. Recall that, here, we are working in the monoidal category of vector spaces. The classical Hopf-Borel theorem  was originally phrased in the  monoidal category of sign-graded vector spaces, cf.\ \ref{graded}. 

Our aim is to generalize this theorem to the \cab types for which $\PP = Vect$.

\medskip

\subsection{Hypotheses}\label{hypo}

 In this section we make the following  assumptions on the given \cab type:

\medskip

\noindent {\texttt{(H0)}}\index{{\texttt{(H0)}}}  for any pair $(\dd,\mu)$ of generating operation $\mu$ and generating cooperation $\dd$ there is a distributive compatibility relation,

\medskip

\noindent  \Hone\index{\Hone} the free $\AA$-algebra $\AA(V)$ is naturally equipped with a \cab structure.

\medskip

\noindent   \Htwoiso \index{\Htwoiso} the natural coalgebra map $\varphi(V): \AA(V)\to \CC^c(V)$ induced by the projection $\proj:\AA(V)\epi V$ is an isomorphism of {\bf S}-modules $\varphi: \AA \cong \CC^c$.

\begin{prop}\label{proptechnique} Let \ca be a type of bialgebras which verifies hypotheses {\texttt{(H0)}} , \texttt {(H1)} and
 \texttt {(H2iso)}. Then the  primitive operad is the identity operad: $\PP = Vect$.
\end{prop}
\begin{proo} It follows from {\texttt{(H0)}}, \Hone  and Theorem \ref{thm:primitive} that there exists a primitive operad $\PP$. Let $\mu\in \PP(n)$ be a nonzero $n$-ary operation for $n\geq 2$. Since $\varphi : \AA(n) \cong \CC^c(n)$ is an isomorphism by  \texttt {(H2iso)}, there exists a cooperation $\dd\in \CC(n)$  such that $\langle \dd , \varphi \mu\rangle = 1$. Let $V= \KK x_1 \oplus \cdots \oplus \KK x_n$. It follows  from Lemma \ref{fi} that
$$\dd \circ \mu (x_1, \ldots , x_n) = x_1\t  \cdots \t  x_n \in \Vtn \subset \AA(V)^{\t n}. $$
Therefore $\dd\circ \mu \neq 0$ and there is a contradiction. Hence we have $\PP(n)=0$ for any $n\geq 2$.
\end{proo}

\begin{prop}\label{bialgebrafunct} Let \ca be a type of bialgebras which verifies hypotheses {\texttt{(H0)}}, \Hone  and \texttt {(H2iso)}.  Let $\HH$ be a $\CC^c\pt \AA$-bialgebra and let $V\to \Prim \HH$ be a linear map. Then the unique algebra lifting $\tilde \aa : \AA(V)\to \HH$ of the composite $\aa : V \to \Prim \HH \mono \HH$  is a bialgebra map.
\end{prop}
\begin{proo} First let us observe that, by Proposition \ref{proptechnique}, we have $\Prim \AA(V) = V$.  Since $\tilde \aa$ is an algebra map by construction, we need only to prove that it is a coalgebra map. We work by induction on the filtration of $\AA(V)$ given by 
$$F_{n}\AA(V) := \bigoplus _{k\leq n}\AA(k)\t _{S_{k}}V^{\t k}.$$
 When $x$ is primitive, that is $x$ lies in $F_{1}\AA(V)=V$, then  $\tilde \aa(x)= \aa (x)$ is primitive by hypothesis. Let $x\in \AA(V)$ be an obstruction of minimal filtration degree $m$, that is an element such that $x\in F_{m}\AA(V)$ and $\DD_{\HH}\circ \tilde \aa(x)\neq (\tilde \aa\t \tilde \aa)\circ \DD_{\AA(V)}(x)$. From the definition of the filtration by the cooperations there exists some cooperation which provides an obstruction of minimal filtration degree $m-1$. But, since for $m=1$ there is no obstruction, we get a contradiction and $\tilde \aa$ is a coalgebra morphism.
\end{proo}

\subsection{The universal idempotent $e$}\label{universalidempotent}\index{universal idempotent}   Let $\HH$ be a \cab. We define a linear map $\oo^{[n]}:\HH \to \HH$ for each $n\geq 2$ as the following composite
$$ \oo^{[n]} : \HH \stackrel{\theta_{n}}{\tto} \CC^c(n)\t_{S_n}\HH^{\t n} \stackrel{\varphi^{-1}\t\Id}{\cong} \AA(n)\t_{S_n} \HH^{\t n} \stackrel{\cc_{n}}{\tto}\HH .$$

We define a linear map $e:\HH \to \HH$ by the formula: 
$$ e = e_{\HH}:= (\Id - \oo^{[2]})(\Id - \oo^{[3]})\cdots (\Id - \oo^{[n]})\cdots \ .$$
We will show that, though $e$ is given by an infinite product, it is well-defined.
We also denote by $e_{\HH}$, or simply $e$, the surjective map $\HH \epi \Im (e)$.

Before stating and proving the main theorem  of this section, we prove some technical results on the universal idempotent $e$.

We denote by $\tilde \iota_{\HH}: \AA(\Prim \HH)\to \HH$ the unique algebra lifting induced by the inclusion map $\iota_{\HH} : \Prim \HH \mono \HH$.

\begin{prop}\label{idempotent} If the  \cab $\HH$ is connected, then the map $e= e_{\HH}:\HH \to \HH$ is well-defined and satisfies the following properties:

\medskip

a) $e_{\AA(V)}= \proj_{V}: \AA(V) \to V$,

b) the image of $e_{\HH}$ is $\Prim \HH$,

c) $e$ is an idempotent: $e^2=e$.

\end{prop}
\begin{proo} Let $(\dd_{i_{1}}^n, \ldots, \dd_{i_{k}}^n)$ be a linear basis of $\CC(n)$ and let 
$(\bar \dd_{i_{1}}^n, \ldots, \bar \dd_{i_{k}}^n)$ be the dual basis (of $\CC^c(n)$). We know by Lemma \ref{fi} that for any $x\in C$ we have
$$\theta_{n}(x) = \sum_{j} \bar \dd_{i_{j}}^n\t \dd_{i_{j}}^n(x)\ .$$
From the connectedness assumption on $C$ (cf.~\ref{connectedness}) it follows that, for any $x\in C$,  there exists an integer $r$ such that $x\in F_r\HH$. Hence we have $\theta_{n}(x)=0$ whenever $n>r$ and therefore $\oo^{[n]}=0$ on $F_r\HH$ whenever $n>r$.
 As a consequence 
$$e(x) = \big((\Id - \oo^{[2]})(\Id - \oo^{[3]})\cdots (\Id - \oo^{[r]})\big)(x),$$
and so $e(x)$ is well-defined.

\medskip

Proof of (a). We consider the following diagram (where $\underline {\t}$ means $\t _{S_n}$):

$$\xymatrix{
 \AA(n)\underline {\t}V^{\t n}\ar[r]^{\varphi\t \Id}\ar@{>->}[d] & \CC^c(n)\underline {\t} V^{\t n}\ar[r]^{\varphi^{-1}\t \Id} \ar@{>->}[d]& \AA(n)\underline {\t} V^{\t n}\ar[r]^{=}\ar@{>->}[d] & \AA(n)\underline {\t} V^{\t n}\ar@{>->}[d] \\
  \AA(V)\ar[r]^-{\theta_{n}} & \CC^c(n)\underline {\t} \AA(V)^{\t n}\ar[r]^{\varphi^{-1}\t \Id} & \AA(n)\underline {\t} \AA(V)^{\t n}\ar[r]^-{\cc_{n}} & \AA(V) \\
 }$$ 
  where the composition in the last line is $\oo^{[n]}$. The left hand side square is commutative by definition of $\varphi$, cf. \ref{cooperad}. The middle square is commutative by construction. The right hand side square is commutative by definition of the $\AA$-structure of $\AA(V)$, cf. \ref{freePalgebra}. As a consequence the whole diagram is commutative. Since, in the diagram,  the lower composite is $\oo^{[n]}$ and the upper composite is the identity, we deduce that the restriction of $\oo^{[n]}$ on the $n$-th component $\AA(V)_n$ is the inclusion into $\AA(V)$. As a consequence $\Id - \oo^{[n]}$ is 0 on the $n$-th component for any $n\geq 2$. So $e$ is the projection on $V= \AA(V)_{1}$ parallel to the higher components, since $e(x)=x$ for any primitive element.
  
 \medskip

Proof of (b). First we remark that the statement is true for $\HH = \AA(V)$ by virtue of (a). Since $\aa$ is a bialgebra morphism by Proposition \ref{bialgebrafunct}, there is a commutative diagram:

$$\xymatrix{
\AA(V) \ar[r]^{\tilde \aa}\ar[d]_{e_{\AA(V)}} & \HH\ar[d]^{e_{\HH}} \\
V\ar[r]^-{\aa} & \Prim \HH \\
}$$

where $V= \Prim \HH$. Statement (a) implies that ${e_{\HH}}$ is surjective.

 \medskip

Proof of (c). From the definition of $e$ we observe that $e(x)=x$ for any $x\in  \Prim \HH $ because $\oo^{[n]}(x) = 0$ for any $n\geq 2$. Since $e(x)$ is primitive by (b) we get $e^2=e$.

\end{proo}

\begin{cor}\label{surjectivity} Let \ca be a type of generalized bialgebras which verifies hypotheses  {\texttt{(H0)}} and \Hone. For any connected bialgebra $\HH$ the natural algebra map $\tilde \iota : \AA(\Prim \HH)\to  \HH$ induced by the inclusion $\iota : \Prim \HH\to  \HH$ is surjective.
\end{cor}
\begin{proo} If $x\in \Prim \HH = \Im \iota$, then clearly $x\in \Im \iota$. Let us now work by induction on the filtration of $\HH$. Assume that $F_{m-1}\HH \subset \Im \tilde \iota$ and let $x\in 
F_{m}\HH$. In the formula 
$$x =e(x) + \big( \sum\oo^{[i]}(x) - \sum\oo^{[i]}\circ \oo^{[j]}(x) + \cdots \big)$$
the first summand $e(x)$ is in $\Prim \HH\subset \Im \tilde \iota$ by Proposition \ref{idempotent}. The second summand is also in $\Im \tilde \iota$ because it is the sum of elements which are products of elements in $\Prim \HH$ by induction. Therefore we proved $x\in \Im \tilde \iota$ for any $x\in \HH$, so $ \tilde \iota$ is surjective.
\end{proo}

\begin{thm}[Rigidity theorem]\label{thm:rigidity} Let $\CC^c$-$\AA$ be a type of generalized bialgebras (over a characteristic zero field) verifying the following hypotheses:

\noindent {\texttt{(H0)}}  the operad $\CC$ is finitely generated and for any pair $(\dd,\mu)$ of generating operation $\mu$ and generating cooperation $\dd$ there is a distributive compatibility relation,

\noindent \Hone the free $\AA$-algebra $\AA(V)$ is naturally equipped with a \cab structure,

 \noindent \Htwoiso\index{\Htwoiso}  the natural coalgebra map  $\varphi(V):\AA(V) \to \CC^c(V)$ is an isomorphism.
 
Then any \cab $\HH$ is free and cofree over its primitive part:
$$\AA(\Prim \HH)\cong \HH \cong \CC^c(\Prim \HH).$$
\end{thm}

\begin{proof} By Proposition \ref{bialgebrafunct} the map $\tilde \iota : \AA(\Prim \HH)\to \HH$ is a bialgebra morphism. On the other hand the projection $e: \HH\to  \Prim \HH$ induces a coalgebra map $\tilde e: \HH \to \CC^c(\Prim \HH)$ by universality (cf.~\ref{cofree}). We will prove that both morphisms are isomorphisms and that the composite
$$ \AA(\Prim \HH)\stackrel{\tilde \iota}{\to} \HH \stackrel{\tilde e}{\to} \CC^c(\Prim \HH)$$
 is $\varphi$.

By Proposition \ref{bialgebrafunct} and the fact that the idempotent is functorial in the bialgebra, there is a commutative diagram

$$\xymatrix{
 & \HH\ar[dr]^{e} & \\
\AA(\Prim \HH) \ar[rr]^{\tilde \aa}\ar[ur]^{\tilde \iota} &  & \Prim \HH\\
}$$

 which induces, by universality of the cofree coalgebra, the commutative diagram:

$$\xymatrix{
 & \HH\ar[dr]^{\tilde e} & \\
\AA(\Prim \HH) \ar[rr]^{\varphi}\ar[ur]^{\tilde \iota} &  &\CC^c( \Prim \HH)\\
}$$

Since $\varphi$ is an isomorphism by \texttt{(H2iso)}, it follows that ${\tilde \iota}$ is injective.

In Proposition \ref{surjectivity} we proved that ${\tilde \iota}$ is surjective, therefore  $e=\tilde \iota : \AA(\Prim \HH)\tto \HH$ is a bialgebra isomorphism and, as a consequence,  $\tilde e$ is also an isomorphism.
\end{proof}

\begin{cor} Let $\HH$ be a connected \cab and let $\HH^2$ denote the image in $\HH$ of $\bigoplus_{n\geq 2}\AA(n)_{S_{n}}\HH^{\t n}$. Then one has $e_{\HH}(\HH^2)=0.$
\end{cor}
\begin{proo} By the rigidity theorem it suffices to show that this assertion is valid when $\HH$ if free. By definition of $e$ we have $e_{\AA(V)}= \proj_{V}$, whose kernel is precisely $\AA(V)^2$.
\end{proo}

\subsection{Explicit universal idempotent}\label{explidemp} Let us suppose that $\AA$ and $\CC$ are given by generators and relations, and that one knows how to describe $\AA(n)$ and $\CC(n)$ explicitly in terms of these generators. Then it makes sense to look for an explicit description of $e$ in terms of the elements of $\AA(n)$ and $\CC(n)$. In the cases already treated in the literature (cf.~for instance \cite{Quillen, LRstr, Livernet06, Holtkamp03, Holtkamp05, Burgunder06, Foissy}, the first step of the proof of the rigidity theorem consists  always in writing down such an explicit idempotent. In the case at hand ( i.e.\ under \Htwoiso) the universal idempotent and the explicit idempotent coincide because, on $\AA(V)$, it is the projection onto $V$ parallel to the other components $\oplus_{n>1} \AA(n)\t_{S_{n}}V^{\t n}$.

The exact form of the compatibility relation(s) depends on the choice of the presentation of $\AA$ and of $\CC$. Let us suppose that hypotheses {\texttt{(H0)}}, \Hone and \Htwoiso hold. Once a linear basis $(\mu_{1}^{[n]}, \ldots, \mu_{k}^{[n]})$ of $\AA(n)$ is chosen, then we can choose, for basis of $\CC(n)$, its dual $(\dd_{1}^{[n]}, \ldots, \dd_{k}^{[n]})$ under the isomorphism $\varphi$: $\langle\varphi(\mu_{i}^{[n]}), \dd_{j}^{[n]}\rangle = 1$ if $i=j$ and $0$ otherwise. Then the compatibility relation of the pair $(\dd_{j}^{[n]}, \mu_{i}^{[n]})$ is such that $\Phi_{1}= \id_{n}$  if $i=j$ and $0$ otherwise.

\section{Triple of operads} \label{s:triple}

We introduce the notion of \emph{triple of operads}\index{triple of operads}

$$(\CC, \AA,  \PP)= (\CC, \between, \AA, F, \PP)$$

\noindent deduced from the prop $\CC^c\--\AA$, that is from a notion of $\CC^c\--\AA$-bialgebra. We construct and study the universal enveloping functor

$$U: \PP\alg \to \AA\alg\ .$$

\subsection{Triple  of operads}\label{triple} Let $(\CC, \between,  \AA)$ be a type of generalized bialgebras. Suppose that the hypotheses {\texttt{(H0)}} and \Hone are  fulfilled, cf.~\ref{hypo}. Then it determines an operad $\PP:= \Prim_{\CC}\AA$ and a  functor $F: \AA\pt alg \to \PP\pt alg$.  Observe that the operad $\PP$  is the largest suboperad of $\AA$ such that any $\PP$-operation applied on primitive elements gives a primitive element. For any \cab $\HH$ the inclusion $\Prim \HH \mono \HH$ becomes a morphism of $\PP$-algebras.  We call this whole structure a \emph{triple of operads} and we denote it by 

\smallskip

$(\CC, \between, \AA, F, \PP)$, or $(\CC, \between, \AA, \PP)$, or more simply $(\CC, \AA,  \PP)$.


\subsection{The map $\varphi$ and the hypothesis \Htwoepi}\label{phi} Since, by hypothesis \Hone, $\AA(V)$ is a \cab, the projection map $\proj_V:\AA(V)\epi V$ determines a unique coalgebra map (cf.~\ref{fi}):
$$\varphi(V): \AA(V)\to \CC^c(V).$$
We recall from \ref{cofreecoalg} that $\varphi(V)$ is the composite
$$\AA(V)\stackrel{\theta(\AA(V))}{\tto}\CC^c(\AA(V)) \stackrel{\CC^c(\proj)}{\tto}\CC^c(V).$$

We denote by $\varphi : \AA \to \CC^c$ the underlying functor of {\bf S}-modules and by $\varphi _{n}:\AA(n) \to \CC^c(n)$ its arity $n$ component.

  We make the following assumption:\\
 

 \noindent \Htwoepi\index{\Htwoepi} the natural coalgebra map $\varphi(V)$ is surjective and admits a natural coalgebra map splitting $s(V):\CC^c(V) \to \AA(V)$,  i.e.\ $\varphi(V)\circ s(V)=\Id_{\CC^{c}(V)}$.
\medskip

 

\subsection{Universal enveloping functor}\label{uef} The functor
$$F: \AA\alg \tto \PP\alg$$
is a \emph{forgetful functor}\index{forgetful functor} in the sense that the composition 
$$\AA\alg \stackrel{F}{\tto} \PP\alg\tto \Vect$$
is the forgetful functor $\AA\alg \to \Vect$. In other words, in passing from an $\AA$-algebra to a $\PP$-algebra we keep the same underlying vector space. Hence this forgetful functor has a left adjoint denoted by 
$$U:  \PP\alg \tto \AA\alg$$
and called the \emph{universal enveloping algebra functor}\index{universal enveloping algebra functor}  (by analogy with the classical case $U: Lie\alg \to As\alg$). Let us recall that adjointness means the following: for any $\PP$-algebra $L$ and any $\AA$-algebra $A$ there is a binatural isomorphism
$$\Hom_{\AA\alg}(U(L), A) = \Hom_{\PP\alg}(L, F(A) ).$$

\begin{prop}\label{uea} Let $L$ be a $\PP$-algebra. The universal enveloping algebra of $L$ is given by 
$$U(L) = \AA(L)/ \sim$$
where the equivalence relation $\sim$  is generated, for any $x_1,\ldots, x_n$ in $L\subset \AA(L)$, by 
$$\mu^{\PP}(x_1,\ldots, x_n) \sim 
(\mu^{\AA};x_1,\ldots, x_n),\quad \mu^{\PP}\in \PP(n),$$
 where $\mu^{\PP}\mapsto \mu^{\AA}$ under the inclusion $\PP(n)\subset \AA(n)$.
\end{prop}
\begin{proo} We have $\mu^{\PP}(x_1,\ldots, x_n)\in L=\AA(1)\t L$ and $(\mu^{\AA};x_1,\ldots, x_n)\in \AA(n)\t_{S_{n}}L^{\t n}$. So the equivalence relation does not respect the graduation. However, it respects the filtration given by 
$$F_{n}\AA(V) := \bigoplus_{j\leq n} \AA(V)_{j}.$$

Let us show that the functor $L\mapsto U(L):=\AA(L)/ \sim$ is left adjoint to the forgetful functor $F$. Let $A$ be an $\AA$-algebra and let $f:L\to F(A)$ be a $\PP$-morphism. There is a unique $\AA$-algebra extension of $f$ to $\AA(L)$ since $\AA(L)$ is free. It is clear that this map passes to the quotient by the equivalence relation and so defines an $\AA$-morphism $U(L)\to A$. 

In the other direction, let $g: U(L)\to A$ be a $\AA$-morphism. Then its restriction to $L$ is a $\PP$-morphism $L\to F(A)$ by Theorem \ref{thm:primitive}. It is immediate to verify that these two constructions are inverse to each other. Therefore we have an isomorphism
$$\Hom_{\AA\alg}(U(L), A) \cong \Hom_{\PP\alg}(L,F( A)),$$
which proves that $U$ is left adjoint to $F$.
\end{proo}

\begin{prop}\label{ueabialg} Under the hypotheses  {\texttt{(H0)}}, \Hone and \texttt{(H2epi)} the universal enveloping algebra $U(L)$ of the $\PP$-algebra $L$  is a \cab.
\end{prop}
\begin{proo} Since we mod out by an ideal, the quotient is an $\AA$-algebra. By hypothesis the free algebra $\AA(L)$ is a \cab. The coalgebra structure of $U(L)$ is induced by the coalgebra structure of $\AA(L)$. For any nontrivial cooperation $\dd$ we have 
$$\dd(\mu^{\PP}(x_1,\ldots, x_n))=0$$
since $\mu^{\PP}(x_1,\ldots, x_n)$ lies in $L$, and we have 
$$\dd(\mu^{\AA}(x_1,\ldots, x_n))=0$$
because $\dd(\mu^{\AA})$ is a primitive operation. Hence for any cooperation $\dd$ we have $\dd(\textrm{relator})=0$. Then $\dd$ is also 0 on $\Ker(U(L) \to \AA(L))$ by the distributivity property of the compatibility relation $\between$.
\end{proo}

\section{Structure Theorem for generalized bialgebras} In this section we show that any triple of operads $(\CC,\AA,\PP)$, which satisfies \Htwoepi, gives rise to a structure theorem analogous to the classical CMM+PBW theorem valid for the triple $(Com, As, Lie)$ (cf.~\ref{thm:ComAs}). It says that any connected \cab is cofree over its primitive part as a  coalgebra and that, as an algebra, it is the universal enveloping algebra over its primitive part.

\begin{thm}[Structure Theorem for generalized bialgebras]\label{thm:structure} Let $\CC^c$-$\AA$ be a type of generalized bialgebras over a field of characteristic zero. Suppose that the following hypotheses are fulfilled:

\noindent {\texttt{(H0)}} for any pair $(\dd,\mu)$ of cooperation $\dd$ and operation $\mu$ there is a distributive compatibility relation,

\noindent \Hone the free $\AA$-algebra $\AA(V)$ is naturally equipped with a \cab structure,

 \noindent \Htwoepi  the natural coalgebra map  $\varphi(V):\AA(V) \to \CC^c(V)$ is surjective and admits a natural coalgebra map splitting $s(V):\CC^c(V) \to \AA(V)$.



\medskip

Then for any \cab $\HH$ the following are equivalent:

\noindent a) the \cab $\HH$ is connected,

\noindent b) there is an isomorphism of bialgebras $\HH \cong U(\Prim \HH)$,

\noindent c) there is an isomorphism of connected coalgebras $\HH \cong \CC^c(\Prim \HH)$.
\end{thm}

We need a construction and two Lemmas before entering the proof of the structure Theorem.
We first introduce a useful terminology.

\subsection{The versal idempotent $e$}\label{versalidempotent} The choice of a coalgebra splitting $s$ permits us to construct a functorial idempotent $e=e_{\HH}:\HH\to \HH$ as follows. 
First we define $\oo^{[n]}:\HH \to \HH$ as the composite 
$$ \oo^{[n]} : \HH \stackrel{\theta_{n}}{\tto} \CC^c(n)\t_{S_n}\HH^{\t n} \stackrel{s(n)\t \Id}{\longrightarrow} \AA(n)\t_{S_n} \HH^{\t n} \stackrel{\cc_{n}}{\tto}\HH .$$

We define a linear map $e:\HH \to \HH$ by the formula: 
$$ e = (\Id - \oo^{[2]})(\Id - \oo^{[3]})\cdots (\Id - \oo^{[n]})\cdots \ .$$

By the very same argument as in the proof of Proposition \ref{idempotent} we show that $e$ is well-defined though it is given by an infinite product.

We will show below that $e$ is an idempotent ($e^2=e$). We call it the \emph{versal idempotent} \index{versal idempotent} (and not universal since it depends on the choice of a splitting). Different choices of splitting lead to different idempotents.

\begin{lemma}\label{basis} We assume hypotheses {\texttt{(H0)}}, \Hone and   \texttt{(H2epi)}. Let $\dd_{1}^{c}, \ldots, \dd_{k}^{c}$ be a basis of $\CC^c(n)$. Let $\mu_{i}:= s(\dd_{i}^{c})\in \AA(n)$ and complete it into a basis $\mu_{1}\ldots , \mu_{k},\mu_{k+1},\ldots , \mu_{l}$ of $\AA(n)$.

Then one has:
\begin{displaymath}\left\{
\begin{array}{ccl}
\dd_{i}\circ \mu_i& = & \id +  \textrm{higher\ terms},\\
\dd_{j}\circ \mu_i& = &0 + \textrm{ higher\ terms,\  when\ } j\neq i ,\\
\end{array}\right .
\end{displaymath}
where ``higher terms" means a sum of some multivalued operations which begin with at least one nontrivial cooperation ($\Phi_{2}$-type multivalued operations).
\end{lemma}
\begin{proo} Graphically, for $n=2$ the statement that we want to prove is

\medskip

$\cpbdeuxdeux$ =\Big( $\cpbA$ or 0\Big) \quad +\quad  $\sum\quad \cpbmult$

\medskip

Since we are interested only in the $\Phi_{1}$-type part, it is sufficient to compute the value of $\dd_{j}\circ \mu_{i}$ on the generic element $x_{1}\cdots x_{n}$ of $\AA(\KK x_{1}\oplus \cdots \oplus \KK x_{n})$. By Proposition \ref{fi} we get 

\begin{displaymath}
\begin{array}{ccl}
\dd_{j}\circ \mu_i(x_{1}\cdots x_{n})& = & \langle \dd_j , \varphi_n (\mu_i) \rangle x_{1}\t \cdots \t x_{n},\\
                                                        & = & \langle \dd_j , \dd_i^c \rangle x_{1}\t \cdots \t x_{n},\\
 & = & (\id\ \textrm{or}\ 0 )x_{1}\t \cdots \t x_{n},\\
\end{array}
\end{displaymath}
depending on $j=i$ or $j\neq i$.
\end{proo}

\begin{lemma}\label{generalidempotent}  If the  \cab $\HH$ is connected, then the map $e= e_{\HH}:\HH \to \HH$ is well-defined, functorial in $\HH$,  and satisfies the following properties:

\medskip

a)  the image of $e$ is $\Prim \HH$,

b) $e$ is an idempotent.

\end{lemma}
\begin{proo} First we observe that, if $x$ is primitive, then $e(x)=x$. Indeed, it is clear that $\oo^{[n]}(x) =0$ for any $n\geq 2$ since $\oo^{[n]}$ begins with a nontrivial cooperation. Hence we have $e(x)=\Id (x)= x$.

Proof of a).  We will prove by induction on $n$ that the image of $F_{n}\HH$ by $e$ lies in $\Prim \HH$.
It is true for $n=1$, since $F_{1}\HH = \Prim \HH$. We use the notation of the previous Lemma.

For $n=2$ we have 
$$\theta_{2}(x) = \sum_{i=1}^{i=k} \mu_{i} \circ \dd_i (x).$$
On $F_{2}\HH$ we have $e=\Id -\oo^{[2]}$. We want to prove that for any $x\in F_{2}\HH$ we have $(\Id -\oo^{[2]})(x)\in F_1\HH = \Prim \HH$, that is,  for any $\dd_j\in \CC(2)$,   $\dd_j(x) = \dd_j\oo^{[2]}(x)$.

We have 
\begin{displaymath}\left\{
\begin{array}{ccl}
\dd_{j}\oo^{[2]}(x)& = & \sum_{i=1}^{k}\dd_j \mu_i \dd_i (x),\\
                       & = &\dd_i (x) + \textrm{ higher\ terms}, \\
\end{array}\right .
\end{displaymath}
by Lemma \ref{basis}. So we have 
$$\dd_j (\Id -\oo^{[2]})(x)= \dd_j (x) - \dd_j (x) + \sum \dd_i \circ  \textrm{ higher\ terms}(x).$$
Since $x\in F_2\HH$, we have $ \dd_i \circ \textrm{higher\ terms}(x)=0$ and therefore $\dd_j (\Id -\oo^{[2]})(x)=0$ as expected.

A similar proof shows that $x\in F_{n}\HH$ implies $(\Id -\oo^{[n]})(x) \in F_{n-1}\HH$. Hence, putting all pieces together, we have shown that   $x\in F_{n}\HH$  implies $e(x) \in F_{1}\HH= \Prim \HH$.
The expected assertion follows from the connectedness of $\HH$.

Proof of b)  Since $e(x) = x$ when $x$ is primitive and since $e(x)\in \Prim \HH$ for any $x\in \HH$, it is clear that $e\circ e = e$.
\end{proo}

\subsection{Proof of the structure Theorem} 

\noindent $(a)\Rightarrow (b)$. Since the functor $U:\PP\alg \to \AA\alg$ is left adjoint to the forgetful functor $F: \AA\alg \to \PP\alg$, the adjoint to the inclusion map $\iota :\Prim \HH\to \HH$ is an algebra map $\aa: U(\Prim \HH) \to \HH$. 

If $x\in \HH$ is in $\Prim \HH$, then obviously it belongs to the image of $\aa$. For any $x\in \HH$ there is an integer $m$ such that $x\in F_{m}\HH$ by the connectedness hypothesis. We now work by induction and suppose that $\aa$ is surjective on $F_{m-1}\HH$. From the formula
$$x = e(x) + \sum_{i }\oo^{[i]}(x) - \sum_{i,j}\oo^{[i]}\oo^{[j]}(x)+ \cdots $$
and the fact that $\oo^{[n]}$ consists in applying a cooperation first and then an operation, it follows that $x-e(x)$ is the sum of products of elements in $F_{m-1}\HH$. From the inductive hypothesis we deduce that $x-e(x)$ is in the image of $\aa$. Since $e(x)\in \Prim \HH$ by Proposition \ref{universalidempotent} we have proved that any $x\in \HH$ belongs to the image of $\aa$ and so $\aa $ is surjective.

The inductive argument as in the proof of Theorem \ref{thm:rigidity} shows the injectivity of $\aa$.

In conclusion the algebra map $\aa: U(\Prim \HH) \to \HH$ is surjective and injective, so it is an isomorphism. It is also a coalgebra map by \ref{bialgebrafunct} and \ref{ueabialg}, so it is a bialgebra isomorphism.
\medskip

\noindent $(b)\Rightarrow (c)$. Let $L$ be a $\PP$-algebra. Since $U$ is left adjoint to the functor $F:\AA\alg \to \PP\alg$ the map $\varphi(L):\AA(L)\to \CC^c(L)$ factors through $U(L)$:
$$\xymatrix{
\AA(L) \ar[rr]^-{\varphi(L))}\ar@{->>}[rd]&&\CC^c(L)\\
& U(L)\ar[ru]&}$$ 

We first show that the map $U(L) \to \CC^c(L)$ is an isomorphism when $L$ is a free $\PP$-algebra. In this case we know that $\Prim U(L)=L$ since $U(\PP(V)) = \AA(V)$ (left-adjointness), and $\Prim \AA(V)= \PP(V)$ by definition of $\PP$.
 Since $U(L)$ is a \cab and $\Prim U(L)=L$ by Proposition \ref{ueabialg}, there is a surjection $e_{U(L)} : U(L) \to L$. By the universality of the cofree coalgebra there is a lifting $\tilde e_{U(L)}$ and so a commutative diagram:

$$\xymatrix{
U(L) \ar[rr]^-{\tilde e_{U(L)}}\ar@{->>}[rd]&&\CC^c(L)\ar@{->>}[ld]\\
& L&}$$ 

Observe that the composite $\AA(L)\epi U(L) \to \CC^c(L)$ is the map $\varphi(L)$, hence $\tilde e_{U(L)}$ is compatible with the filtration.

We claim that the associated morphism $gr (U(L))\to gr \CC^c(L)= \CC^c(L)$ on the graded objects is an isomorphism. Indeed, the quotient $F_{n}U(L) / F_{n-1}U(L)$ consists in moding out $\AA(L)_{1}\oplus \cdots \oplus \AA(L)_{n} $ by the subspace $J_{n}$ generated by the elements $\mu^{\PP}(x_{1}\ldots x_{n}) - \mu^{\AA}(x_{1}\ldots x_{n})$ and by $\AA(L)_{1}\oplus \cdots \oplus \AA(L)_{n-1} $. On the other hand, the quotient  $F_{n} \CC^c(L) / F_{n-1} \CC^c(L)$ is simply $ \CC^c(L)_{n}$, which is the quotient of $\AA(L)_{n}$ by the homogeneous degree $n$ part of $J_{n}$. These two quotients are the same because the relations\\
$$u+v\sim 0 \textrm{ and } v\sim 0,$$
and
$$u\sim 0 \textrm{ and } v\sim 0,$$
are equivalent (recall that $\mu^{\PP}(x_{1}\ldots x_{n})$ is in degree 1).

Let us now prove that $U(L) \to \CC^c(L)$ is an isomorphim for any Lie algebra $L$. Let 
$$L_{1} \to L_{0}  \epi L \to 0$$
be a free resolution of $L$. Since the morphism $U(V)\to \CC^c(V)$ is natural in $V$, the isomorphisms for $L_{0}$ and $L_{1}$ imply the isomorphism for $L$.
\medskip

The implication $(c) \Rightarrow (a)$ is a tautology. \hfill $\square$

\subsection{Good triple of operads}  If a triple of operads $(\CC, \AA, \PP)$ satisfies the structure theorem, then we call it a \emph{good triple of operads}\index{good triple}. So Theorem \ref{thm:structure} shows that, if the triple of operads $(\CC, \AA, \PP)$ satisfies the hypothesis \Htwoepi, then it is a good triple. Conversely, if the triple is good, then the coalgebra isomorphism $\AA(V)\cong \CC^c(\PP(V))$ composed with the projection induced by $\proj: \PP(V) \to V$ defines $\AA(V)\to \CC^c(V)$, which is a splitting of $\varphi(V)$. Hence the hypothesis \Htwoepi is fulfilled. Therefore the triple \CAP is good if and only if the hypothesis \Htwoepi is fulfilled.

\subsection{About the verification of the hypotheses}\label{verification} 
\subsubsection{{\texttt{(H0)}}} The hypothesis {\texttt{(H0)}} (distributivity of the compatibility relation) is, in general, immediate to check by direct inspection. Observe that, when the operads $\CC$ and $\AA$ are given by generators and relations, it suffices to check the compatibility relations on the pairs $(\dd, \mu)$ when they are both generators. 

\subsubsection{\Hone} In order to verify hypothesis \Hone there are, essentially, three strategies. 

\noindent (1). When the free algebra $\AA(V)$ is known explicitly (for instance a basis of $\AA(n)$ is identified with some explicit combinatorial objects), then one can usually construct explicitly the generating cooperations and check that they satisfy the relations of $\CC^c$ and the compatibility relations. 

\noindent (2). Another strategy consists in taking advantage of the distributivity of the compatibility relations. One constructs inductively the cooperations on $\AA(V) = \oplus_{n\geq 1}\AA(V)_{n}$ by sending $V$ to $0$ and then  by using the compatibility relations. Then one checks, again inductively, that they satisfy the relations of $\CC^c$ and the compatibility relations. This is very close to the techniques used in ``rewriting systems'', cf.~\cite{Lafont97} and \ref{IC}.

\noindent (3). The third strategy consists in viewing a given triple as a quotient of a good triple. It is given in Proposition \ref{prop:quotient} below.

\subsubsection{\Htwoepi} The map $\varphi$ is in fact a map of {\bf S}-modules. Therefore it sends the degree $n$ part of $\AA$ to the degree $n$ part of $\CC^c$. So in order to check surjectivity, it is sufficient to compute the composite $\dd \circ \mu$ for any pair $(\dd,\mu)$ where $\mu$ is a linear generator of $\AA(n)$ and $\dd$ is a linear generator of $\CC^c(n)$. From this functorial property of $\varphi$ we deduce that the element $\dd \circ \mu\row x1n$ is of the form $\sum_{\ss}a_{\ss} (x_{\ss(1)}, \ldots , x_{\ss(n)})$. Let $a(\mu, \dd)$ be the coefficient $a_{\id}$ of this sum. The map $\varphi$ is
$$\varphi_{n}(\mu) = \sum_{\dd} a(\mu, \dd) \dd$$
where the sum is over a basis of $\CC^c(n)$.

\section{A few consequences of the structure theorem} We derive a  few consequences of the structure theorem, namely by applying it to the free algebra $\AA(V)$. It gives some criterion to check if a given triple of operads has some chances to be good.

\subsection{From the structure theorem to the rigidity theorem} The rigidity theorem is a Corollary of the structure theorem. Indeed, if the hypothesis \texttt{(H2iso)} holds, then  \texttt{(H2epi)} holds (unique choice for the splitting) and the primitive operad is $Vect$ by Proposition \ref{proptechnique}. Hence the triple $(\CC, \AA, Vect)$ is a good triple and the functor $F$ is simply the forgetful functor to $\Vect$. The left adjoint functor of $F$ is the free $\AA$-algebra functor, so item (b) in the structure theorem becomes $\HH\cong \AA(\Prim \HH)$. So $\HH$ is free and cofree over its primitive part, as claimed in Theorem \ref{thm:rigidity}.

\subsection{Dualization} Observe that if $(\CC, \AA, Vect)$ is a good triple of operads, then so is the triple $(\AA,\CC,Vect)$. The compatibility relation(s) is obtained by dualization, i.e.\ reading $\between$ upsidedown. The new map $\varphi$ is simply the linear dual of the former one.

\begin{thm} If $(\CC, \AA, \PP)$ is a good triple of operads over the field $\KK$, then there is an equivalence of categories between the category of connected (i.e.\ conilpotent) \cabs and the category of $\PP$-algebras:

\begin{displaymath} \{\textrm{con.\ } \CC^c\pt\AA\bialg \} 
{\begin{array}{c} U\\
\leftrightarrows\\
\Prim 
\end{array}}
\{\PP\alg \} \ .
\end{displaymath}
\end{thm}
\begin{proo} We already know that if $L$ is a free $\PP$-algebra,  i.e.\ $L=\PP(V)$, then $\Prim U(L) = \Prim U(\PP(V))= \Prim \AA(V) = \PP(V)=L$. By the same argument as in the proof of $(b) \Rightarrow (c)$ it is true for any $\PP$-algebra $L$.

In the other direction, let $\HH$ be a connected \cab. By item (b) in Theorem \ref{thm:structure} we have an isomorphism $\HH \cong U(\Prim \HH)$.
\end{proo}

\begin{prop}\label{isosms}  If \CAP is a good triple of operads, then there is an isomorphism of Schur functors:
$$\AA \cong \CC^c \circ \PP\ .$$
\end{prop}       
\begin{proo} It suffices to apply the structure theorem to the free algebra $\AA(V)$, which is a \cab by hypothesis. It is connected because $\AA(1)= \KK$.

Since the composite of left adjoint functors is still left adjoint, we have $U(\PP(V)) = \AA(V)$. Hence, by the structure Theorem \ref{thm:structure}, we get the expected isomorphism.
\end{proo}
\begin{cor}\label{cor:genser}  If \CAP is a good triple of operads, then there is an identity of formal power series:
$$f^{\AA}(t) = f^{\CC}( f^{\PP}(t)).$$
\end{cor}    
\begin{proo} Since $\CC^c(n)$ is finite dimensional, the two Schur functors $\CC$ and $\CC^c$ have well-defined generating series which are equal. The formula follows from Proposition \ref{isosms} and the computation of the generating series of a composite of Schur functors, cf. \ref{genseries}.
\end{proo}

\subsection{Searching for good triples} Observe that this relationship entertwining the generating series gives a criterion to the possible existence of a good triple. Indeed, let us suppose that we start with a forgetful functor $\AA\alg \stackrel{F}{\to} \PP\alg$ and we would like to know if it can be part of a good triple $(\CC, \between, \AA, F, \PP)$. Then, there should exist a power series $c(t)= \sum_{n\geq 1} \frac{c(n)}{n!} t^n$ where the coefficients $c(n)$ are integers (and $c(1)=1$), such that $f^{\AA}(t) =c( f^{\PP}(t)) $. 

For instance, if $(Com, \AA, \PP)$ is a good triple of operads and if $\BB\alg \to \AA\alg $ is a forgetful functor, then there is no  good triple for the composite  $\BB\alg \to \PP\alg $ (unless $\BB = \AA$).

\subsection{Frobenius characteristic} There is an invariant which is finer than the generating series. It consists in taking the Frobenius characteristic of the Schur functor in the symmetric functions. Indeed the isomorphism
$$\AA(n)\cong \sum_{i_{1}+\cdots +i_{k}=n}\CC^{c}(k) \t_{S_{k}}\Ind ^{S_{n}}_{S_{i_{1}}+\cdots +S_{i_{k}}}\big(\PP(i_{1})\t \cdots \t \PP(i_{k})\big)$$
implies that the composite of the Frobenius characteristic of $\CC^c$ and $\PP$ is the Frobenius characteristic of $\AA$. See \cite{LPo} for an explicit example.



\chapter{Applications and variations}

In this chapter we give a few applications of the structure theorem, some generalizations and we give some general constructions to obtain a good triple of operads. Concrete examples will be given in the next chapters.

One of the most easy ways of constructing a good triple from an existing triple is to mod out by primitive relators. It gives rise to many examples.

 There are some techniques to obtain triples of the form  $(As, \AA, \Prim_{As}\AA)$. For instance one can assume that $\AA$ is a Hopf operad, that is, the Schur functor is a functor to coalgebras. Another assumption is to suppose that there exists an associative operation verifying some good properties (multiplicative operad).

 In the nonsymmetric case, we do not need the characteristic zero hypothesis to get a good triple, so, under this hypothesis, the structure theorem is valid in a characteristic free context.
   
  We show how Koszul duality should help to construct good triples out of existing ones. 
  
  Our basic category is the category $\Vect$ of vector spaces. It is a linear symmetric monoidal category and this  is exactly the structure that we used. So there is an immediate extension of our main theorem to any linear symmetric monoidal category, for instance the category of sign-graded vector spaces and the category of {\bf S}-modules.
  
  We can reverse the roles of algebraic and coalgebraic structures. Then the primitives are replaced by the indecomposables and we obtain a ``dual" result.
  
  The classical result (PBW+CMM) admits a characteristic $p$ variant. We expect similar generalizations in characteristic $p$ and we explain how to modify the operad framework to do so.
  
  We  mention briefly the relationship with rewriting systems in computer sciences.
  
   Finally, we give an application to a natural problem in representation theory of the symmetric groups.

\section{Quotient triple}\label{quotient}

Let $\CC^c\pt \AA$ be a type of bialgebras, which satisfies the hypotheses {\texttt{(H0)}}, \Hone and \Htwoepi. So we start with a good triple of operads $(\CC, \AA, \PP)$. We will show that moding out by primitive operations gives rise to other good triples of operads.

\begin{prop}\label{prop:quotient} Let $\CC^c\pt \AA$ be a type of bialgebras, which satisfies the hypotheses {\texttt{(H0)}}, \Hone and \Htwoepi. Let $J$ be an operadic ideal of $\AA$ generated by (nontrivial) primitive operations. 
Then $(\CC, \AA/J, \Prim_{\CC}\AA/J)$ is a good triple of operads.
\end{prop}
\begin{proo} Let $\PP= \Prim_{\CC}\AA$. We have $\PP=\KK\id\oplus \bar \PP$. In this proof a primitive operation is an element in $\bar \PP(n)$ for some $n$ (so we exclude $\id$). 

Since {\texttt{(H0)}} is fulfilled for $\CC^c\pt \AA$, it is also fulfilled for $\CC^c\pt \AA/J$.

By hypothesis the ideal $J$ is linearly generated by the composites $\mu \circ \row {\mu}1n$ where at least one of the operations is primitive. If $\mu$ is primitive, then the $\Phi_{1}$-part of $\dd\circ \mu$ is $0$. If $\mu_{l}$ is primitive for some index $l$, then the $\Phi_{1}$-part of $\dd\circ \mu\circ \row {\mu}1n$ is of the form $\sum \ss \circ \row {\mu}1n$ for some permutations $\ss$. Since $\mu_{l}$ is primitive, the value of this operation on a generic element is $0$ in the quotient $\AA/J$. Therefore $(\AA/J)(V)$ is a $\CC^c \pt \AA/J$-bialgebra. Hence hypothesis \Hone is fulfilled.

The above calculation shows that the map $\varphi: \AA \to \CC^c$ factors through $\AA/J$ and the resulting map $ \AA/J \to \CC^c$ is the map $\varphi$ for   $\CC^c\pt \AA/J$. Hence the composite
$$ \CC^c\buildrel{s}\over{\longrightarrow}\AA \epi \AA/J $$
is a splitting of $\varphi$. So hypothesis \Htwoepi 
is fulfilled.

Therefore, by Theorem \ref{thm:structure}, $(\CC, \AA/J, \Prim_{\CC}\AA/J)$ is a good triple of operads.
\end{proo}

\begin{cor}\label{cor:quotient}
Any good triple of operads $(\CC,\AA,\PP)$ determines a good triple of operads $(\CC, \ZZZ,Vect)$, called its quotient triple\index{quotient triple},  where the operad $\ZZZ$ is the quotient of $\AA$ by the operadic ideal $(\bar \PP)$ generated by the nontrivial primitive operations: $\ZZZ = \AA/(\bar \PP)$.
\end{cor}

\subsection{Remark} Observe that in many cases, including the classical case, the \sm  isomorphism $\CC\cong \AA/(\bar \PP)$ is in fact an operad isomorphism. But this is not always true, see for instance examples \ref{Zinbiel} and \ref{NAP}.

\begin{thm}[Analogue of the classical PBW Theorem]\label{classicalPBW} Let $(\CC, \AA, \PP)$ be a good triple of operads, and let $\ZZZ := \AA/(\bar \PP)$ be the quotient operad of $\AA$ by the ideal generated by the (nontrivial) primitive operations. Then, for any $\PP$-algebra $L$ there is an isomorphism of $\ZZZ$-algebras:
$$ \ZZZ(L) \to gr\ U(L)\ .$$
\end{thm}
\begin{proo}  First we observe that $gr\ U(L)$ is a $\ZZZ$-algebra by direct inspection of the structure of $U(L)$, cf. \ref{uea}. The composite map 
$$L\to \AA(L) \to U(L) \to gr_{1}U(L)\subset U(L)$$
 induces a $\ZZZ$-algebra map $ \ZZZ(L) \to gr\ U(L)$. The commutativity of the diagram

$$\xymatrix{
\ZZZ(L)  \ar[r]\ar@{->>}[rd]&gr\ U(L)\ar@{->>}[d] \ar[r]&\CC^c(L)\ar@{->>}[ld]\\
& L&}$$ 
shows that the composite of the horizontal arrows is the isomorphism $\ZZZ(L) \to \CC^c(L)$ coming from the good triple $(\CC, \ZZZ, Vect)$. Since $gr\ U(L)\to \CC^c(L)$ is an isomorphism (cf. the proof of $(b)\Rightarrow (c)$ in Theorem \ref{thm:structure}), we are done.
\end{proo}

\subsection{Remark on PBW}\label{rk:PBW} In the classical case $(Com, As, Lie)$ (see \ref{ex:comas} for details) the isomorphism of commutative algebras
$$ S(\gg) \to gr\ U(\gg)$$
is often called the Poincar\'e-Birkhoff-Witt theorem (cf.~\cite{Grivel}).

\subsection{Split triple of operads}\label{splittriple} Let \CAP be a triple of operads and let $\ZZZ:= \AA/(\bar \PP)$ be the quotient operad. We say that \CAP is a \emph{split triple}\index{split triple} if there is a morphism of operads $s:\ZZZ\mono \AA$ such that the composite $\ZZZ \mono \AA \epi \ZZZ$ is the identity and the map $s(V):\ZZZ(V) \to \AA(V)$ is a $\CC^c\--\ZZZ$-bialgebra morphism. For instance the triple $(Com, As, Lie)$ is not split, and the triple $(As, Dup, Mag)$ (cf.~\ref{AsDupepi}) admits two different splittings.

\begin{prop} Let \CAP be a split triple of operads. Then any \cab $\HH$  is also a $\CC^c\--\ZZZ$-bialgebra and the idempotent $e_{\HH}$ is the same in both cases.
\end{prop}

\begin{proo} Since, by hypothesis, the splitting $s$ induces a morphism of bialgebras on the free algebras, it induces a morphism of props $\CC^c\--\ZZZ\to \CC^c\--\AA$, or, equivalently, a functor between the category of bialgebras $\CC^c\--\AA\--bialg$ and $\CC^c\--\ZZZ\--bialg$. In the construction of the idempotent $e_{\HH}$ for $\CC^c\--\AA\--$bialgebras we need a coalgebra splitting $\CC^c\to \AA$. We can take the composite $\CC^c\cong \ZZZ\buildrel{s}\over {\to} \AA$. From the construction of $e_{\HH}$ (cf. \ref{universalidempotent}) it follows that we get precisely the universal idempotent for $\CC^c\--\ZZZ\--$bialgebras.
\end{proo}

\section{Hopf operad, multiplicative operad}\label{Hopfoperad}

Under some reasonable assumptions on an operad $\PP$ we can show that the tensor product of two $\PP$-algebras is still a $\PP$-algebra. As a consequence one can equip the free $\PP$-algebra with a coassociative cooperation. It gives rise to a notion of $As^c\--\PP$-bialgebras. Quite often the assumptions are easy to verify and show immediately that hypothesis \Hone is fulfilled. We present two cases: Hopf operad and  multiplicative operad.

\subsection{Hopf operad} 
By definition a \emph{Hopf operad}\index{Hopf operad} is an operad $\PP$ in the category of coalgebras (cf. \cite{Moerdijk} for instance). Moreover we assume that $\PP(0)=\KK$, so $\PP$-algebras have a unit (the image of $1\in \PP(0)$). Explicitly the spaces $\PP(n)$ are coalgebras,  i.e.\ they are equipped with a coassociative map $\DD : \PP(n) \to \PP(n)\t \PP(n)$, compatible with the operad structure. For instance a set-theoretic operad gives rise to a Hopf operad by using the diagonal on sets. As a consequence the tensor product of two $\PP$-algebras $A$ and $B$ is a $\PP\t \PP$-algebra (where $\PP\t \PP$ is the Hadamard product, i.e.\ $(\PP\t \PP)(n)=\PP(n)\t \PP(n$
) and the map $\DD$ makes it into a $\PP$-algebra.

An algebra over $\PP$ in this framework is a coalgebra equipped with a coalgebra map $\PP(A) \to A$. Hence we get a notion of $As^c\--\PP$-bialgebra. In particular the free $\PP$-algebra is a $As^c\--\PP$-bialgebra since there is a unique $\PP$-algebra morphism
$$\PP(V)\to \PP(V)\t \PP(V)$$
which extends $v\mapsto v\t 1 + 1 \t v$. Several examples of triples of operads $(As, \PP, \Prim_{As} \PP)$ are of this type (see Chapter 4). 

\subsection{Multiplicative operad \cite{JLLsci}}\label{split}
Let $\PP$ be a binary quadratic nonsymmetric operad which contains an associative operation, denoted $*$ (this hypothesis is sometimes called ``split associativity'' as $*$ comes, often,  as the sum of the generating operations). In other words we suppose that there is given a morphism of operads $As\to \PP$. We call it a \emph{multiplicative operad}\index{multiplicative operad}. We suppose that there is a partial unit $1$ in the following sense. We give ourselves two maps
$$\aa:\PP_{2} \to \PP_{1} =\KK\quad \textrm{and} \quad \beta:\PP_{2}  \to \PP_{1} =\KK$$
which give a meaning to $x\circ 1$ and $1\circ x$ for any $\circ \in \PP_{2} $ and any $x\in A$ (where $A$ is a $\PP$-algebra):
$$x\circ 1 = \aa(\circ) x \quad \textrm \quad 1\circ x = \bb (\circ) x.$$
We always assume that $1$ is a two-sided unit for $*$ ( i.e.\ $\aa(*)=1=\bb(*)$). Observe that we do not require  $1\circ 1$ to be defined. Let $A_{+}= A\oplus \KK 1_{A}$ be the augmented algebra. For two $\PP$-algebras $A$ and $B$ the augmentation ideal of $A_{+}\t B_{+}$ is
$$A\t \KK 1_{B} \ \oplus \ \KK 1_{A} \t B \ \oplus \ A\t B\ $$
The \emph{Ronco's trick}\index{Ronco's trick} (cf.\ \cite{Ronco00}) consists in constructing an operation $\circ$ on the augmentation ideal as follows:
$$(a\t b) \circ (a'\t b') = (a*a')\t (b\circ b')$$
whenever all the terms are defined, and (when $b=1_{B}=b'$)
$$(a\t 1_{B}) \circ (a'\t 1_{B}) = (a\circ a')\t (1_{A\t B}).$$
Observe that the relations of $\PP$ are verified for any $a,a'\in A$ and $b,b'\in B$. If they are also verified in all the other cases, then the choice of $\aa$ and $\bb$ is said to be \emph{coherent}\index{coherent} with $\PP$.

It was proved in \cite {JLLsci} that, under this coherence assumption, the free $\PP$-algebra $\PP(V)$ is equipped with a coassociative coproduct $\dd$. It is constructed as follows. By hypothesis there is a $\PP$-algebra structure on
$$\PP(V)\t \KK 1 \ \oplus \ \KK 1 \t \PP(V) \ \oplus \ \PP(V)\t \PP(V)\ $$
We define $\DD$ from $\PP(V)$  to the $\PP$-algebra above as the unique $\PP$-algebra map  which sends $v\in V$ to $v\t 1 + 1 \t v$. The projection to $\PP(V)\t \PP(V)$ gives the expected map $\dd$. 

From this construction we get a well-defined notion of $As^c\--\PP$-bialgebra for which the hypotheses {\texttt{(H0)}} and \Hone are fulfilled.
Of course, this construction can be re-written in the nonunital context. But, then, the formulas are more complicated to handle (see \ref{Hopf}).

In many cases we get a good triple of operads $(As, \PP, \Prim_{As}\PP)$, see \cite {JLLsci}. Several cases will be described in Chapter 4. The interesting point about these examples is that $\PP(V)$ is a Hopf algebra in the classical sense. In fact many \emph{combinatorial Hopf algebras}\index{combinatorial Hopf algebra} (cf.\ \cite{LRcha}) can be constructed this way.

A refinement of this method gives triples of the form $(\AA,\AA, Vect)$, cf.~\ref{AAVect}.

\section{The nonsymmetric case}\label{nonsymmetriccase} In the preceding chapter we always made the hypothesis: $\KK$ is a characteristic zero field. The reason was the following. In the interplay between operad and cooperad we had to identify invariants and coinvariants, cf.~\ref{cooperad}. There is an environment for which this hypothesis is not necessary, it is the nonsymmetric case (for the characteristic $p$ case, see \ref{charp}).

We suppose that $\CC$ and $\AA$ are nonsymmetric operads (cf.~\ref{nsoperad}) and that, in the compatibility relations, the only permutation, which is involved, is the identity (cf~\ref{diagram}). Such a type is called a \emph{nonsymmetric type} of generalized bialgebras (nonsymmetric prop). In hypothesis \texttt{(H2epi)} we suppose that there is a cooperad splitting of the form $\CC^c_{n} \to \AA_{n}$,  i.e.\ not involving the symmetric group. This is called the \emph{nonsymmetric version} of \texttt{(H2epi)}. Then the very same proof as in the structure Theorem can be performed and we get the following result.

\begin{thm} [Structure Theorem for ns generalized bialgebras]\label{thm:structurens}

 Let $\CC^c$-$\AA$ be a nonsymmetric type of generalized bialgebras over a field $\KK$. Suppose that the hypotheses \texttt{(H1)}  and nonsymmetric \texttt{(H2epi)} are fulfilled. 

Then the good triple of operads \CAP has the following property.

\medskip

For any \cab $\HH$ the following are equivalent:

\noindent a) the \cab $\HH$ is connected,

\noindent b) there is an isomorphism of bialgebras $\HH \cong U(\Prim \HH)$,

\noindent c) there is an isomorphism of connected coalgebras $\HH \cong \CC^c(\Prim \HH)$.
\end{thm}

The following Rigidity Theorem is a Corollary of the Structure Theorem in the nonsymmetric case.

\begin{thm}[Rigidity Theorem, nonsymmetric case]  Let $\CC^c$-$\AA$ be a nonsymmetric type of generalized bialgebras over a field $\KK$. Suppose that the hypotheses \texttt{(H1)}  and \texttt{(H2iso)} are fulfilled. 

Then any \cab $\HH$ is free and cofree over its primitive part:
$$\AA(\Prim \HH)\cong \HH \cong \CC^c(\Prim \HH).$$
\end{thm}

Explicit examples will be given in the next section.


\section{Koszul duality and triples}\label{Koszul}

We provide a method to construct good triples of operads by using the Koszul duality for operads.

\subsection{Koszul duality of quadratic operads} Let us recall briefly from \cite{GK} and \cite{Fresse04} (see also \cite{LV} for details) that any quadratic operad $\PP$ gives rise to a dual operad $\PP^!$. It is also quadratic and $(\PP^!)^!=\PP$. For instance $As^!=As, Com^!=Lie, Mag^!=Nil$. When the operad $\PP$ is binary, then the generating series of $\PP$ and of $\PP^!$ are related by the formula:
$$f^{\PP^!}(-f^{\PP}(-t))=t,$$
(in the Appendix of \cite{JLLdig} there is a brief account of Koszul duality of operads). When the operad $\PP$ is $k$-ary, one needs to introduce the skew-generating series
 $$g^{\PP}(t):=\sum_{n\geq 1}(-1)^k\frac{\dim \PP((k-1)n+1)}{n!}\ t^{((k-1)n+1)} .$$
 If the operad $\PP$ is Koszul, then Vallette proved in  \cite{Vallette05}) the formula:
 $$f^{\PP^!}(-g^{\PP}(-t))=t.$$
 In the most general case (generating operations of any arity), it is best to work with a series in two variables, cf. \cite{Vallette05} for details.
 
 \subsection{Extension of operads}\label{extension} Let us say that the sequence of operads
 $$\PP \mono \AA \epi \ZZZ$$
 is an \emph{extension of operads}\index{extension of operads} if $\mono$ is a monomorphism, $\epi$ is an epimorphism, and if there is an isomorphism of {\bf S}-modules (which is part of the structure) $\AA\cong \ZZZ \circ \PP$ such that $\Id_{\Vect}\mono \ZZZ$ induces  $\PP \mono \AA$, and  $\PP\epi \Id_{\Vect}$ induces  $\AA\epi \PP$. Under this hypothesis we say that $\AA$ is an \emph{extension} of $\ZZZ$ by $\PP$. For instance $As$ and $Pois$ are both extensions of $Com$ by $Lie$.
 
 In many cases where all the operads are quadratic we can check that 
  $$\ZZZ^! \mono \AA^! \epi \PP^!$$
is also an extension of operads (Exercise: show that it works at the level of generating functions). For instance the classical extension
$$Lie \mono As \epi Com$$
is self-dual.

Suppose that \CAP is a good triple and that $\PP^!$ is part of a good triple $(\QQQ^!, \PP^!, Vect)$. We recall that we have also a good triple $(\CC, \ZZZ, Vect)$. Then comparing the generating functions we can expect the existence of a good triple of the form:
$$(\QQQ^!, \AA^!, \ZZZ^!).$$
In the case where $\QQQ=\PP$ and $\CC=\ZZZ$, we would get the triple
$$(\PP^!, \AA^!, \CC^!).$$
Similar structures have been studied in \cite{FM}.

\section{Other symmetric monoidal categories}\label{graded}

Until now our ground category was the category of vector spaces over a field. The only property of $\Vect$ that we used is that it is a symmetric monoidal category. Hence we can replace it by any other symmetric monoidal category.

\subsection{Graded vector space} The category of \emph{graded vector spaces} (more accurately we should say the \emph{sign-graded vector spaces}\index{sign-graded vector spaces}) is a symmetric monoidal category.  Recall that the objects are the graded vector spaces $\{V_{n}\}_{n\geq 0}$ and the symmetric isomorphism (twisting map) is given by 
$$\tau(x\t y) = (-1)^{|x||y|}y \t x\ ,$$
where $x$ and $y$ are homogeneous elements of degree $|x|$ and $|y|$ respectively.

\subsection{Structure theorem in the graded  case} The main result (cf.~\ref{thm:structure}) holds in this more general setting because, as already said, our proofs use only the symmetric monoidal properties of $\Vect$. The result is not significantly different when the operads are nonsymmetric (since the symmetric group does not play any role). However it is different for general operads since, for instance, the free commutative algebra over an odd-degree vector space is the exterior algebra instead of the symmetric algebra. In particular when the vector space is finite dimensional, the exterior algebra is finite dimensional, while the symmetric algebra is not.

In algebraic topology it is the sign-graded framework which is relevant.

\subsection{Structure theorem for twisted bialgebras} The category of {\bf S}-modules can be equipped with a symmetric  product as follows. Let $M$ and $N$ be two {\bf S}-modules. We define their tensor product $M\t N$ by
$$(M\t N) (V) := M(V) \t N(V)\ .$$
Here we use the interpretation of {\bf S}-modules in terms of endofunctors of $\Vect$ (Schur functors, see \ref{Smod}).  For any operad $\AA$ one can define a notion of $\AA$-algebra in the category of {\bf S}-modules. They are sometimes called \emph{twisted $\AA$-algebras}. Similarly, given a prop $\CC^c\pt \AA$, there is an obvious notion of \emph{twisted $\CC^c\pt \AA$-bialgebra} over this prop. The structure theorem admits an obvious extension to twisted bialgebras.


\subsection{Generalization to colored operads} A category is a generalization of a monoid in the sense that the composition of two elements is defined only if certain conditions are fulfilled (source of one = target of the other). Similarly there is a generalization of the notion of operads in which the composition of operations is defined only if some conditions are fulfilled among the operations. This is called a \emph{colored operad} or a \emph{multicategory}. One should be able to write a structure theorem in this framework. See \cite{Ronco06} for an example.

\subsection{Generalization of the nonsymmetric case} In the nonsymmetric case we don't even use the symmetry properties of the monoidal category $\Vect$. Hence we can extend our theorem to other monoidal categories. For instance we can replace $(\Vect, \t)$ by the category of {\bf S}-modules equipped with the composition product $\circ $ , cf.~\ref{compositionSmod}. Then, there are notions of (generalized) algebras, coalgebras and bialgebras in this context and we can extend our main theorem. For instance the analogue of associative algebra (resp. associative coalgebra) is the notion of operad (resp. cooperad). So a \emph{unital infinitesimal bioperad}\index{unital infinitesimal bioperad} $\BB$ is an \sm $\BB$ equipped with an operad structure $\cc :\BB \circ \BB \to \BB$ and a cooperad structure $\theta : \BB \to \BB\circ \BB$ satisfying the following compatibility relation 
$$\theta \cc = - \id_{\BB\circ\BB} + (\id_{\BB}\circ \theta)(\cc\circ \id_{\BB}) + (\theta\circ \id_{\BB})(\id_{\BB}\circ \cc)$$
which is nothing but the unital infinitesimal relation \ref{ex2}. Observe that we denoted the composition of functors by concatenation (for instance $\theta \cc$) to avoid confusion with the composition of {\bf S}-modules. The generalization of our theorem says that the only example of unital infinitesimal bioperad is the free operad.

Observe that there is no such object as Hopf bioperad since the monoidal category $({\bf S}\pt mod, \circ)$ is not symmetric.

We plan to come back to this notion of generalized bioperads in a subsequent paper.

\section{Coalgebraic version}\label{coalgebraic}

Let $\CC^c\--\AA$ be a bialgebra type. The notion of ``indecomposable'' is dual to the notion of ``primitive''. By definition, the \emph{indecomposable space}\index{indecomposable space} of a $\CC^c\--\AA$-bialgebra $\HH$ is the quotient
$$\Indec \HH := \HH/\HH^2$$
where $\HH^2$ is the image of $\bigoplus_{n\geq 2}\AA(n)\t_{S_{n}}\HH^{\t n} $ in $\HH$ under $\cc$. Observe that it depends only on the $\AA$-algebra structure of $\HH$. In general $\Indec \HH$ is not a $\CC^c$-coalgebra, but we will construct a quotient cooperad of $\CC^c$ on which $\Indec \HH$ is a coalgebra.

\begin{prop} Suppose that the bialgebra type $\CC^c\--\AA$ satisfies {\texttt{(H0)}} and

\texttt{(H1${}^c$)} the cofree $\CC^c$-coalgebra $\CC^c(V)$ is equipped with a natural $\CC^c\--\AA$-bialgebra structure.

Then the \sm $\QQQ^c(V)=\Indec_{\AA}\CC^c(V) := \CC^c(V)/ \CC^c(V)^2$ inherits a cooperad structure from $\CC^c$. 

Moreover for any $\CC^c\--\AA$-bialgebra $\HH$ the indecomposable space $\Indec \HH$ is a $\QQQ^c$-coalgebra, and the surjection $\HH \to \Indec \HH$ is a $\QQQ^c$-coalgebra morphism.
\end{prop}

\begin{proo} It suffices to dualize the proof of \ref{thm:structure}.
\end{proo}

\noindent{\bf Example}. Let $As^c\--Com$ be the type ``commutative (classical) bialgebras". Then $\QQQ=Lie$ and the surjection $As^c\to Lie^c$ is simply the dual of the inclusion $Lie \to As$. Explicitly the coLie structure of the coassociative coalgebra $(C, \dd)$ is given by $(\id - \tau) \dd$ where $\tau$ is the twisting map.

Since $\QQQ^c$ is a quotient of $\CC^c$, there is a forgetful functor
$$F^c : \CC^c\--coalg \to \QQQ^c\--coalg $$
which admits a \emph{right} adjoint, that we denote by
$$F^c : \QQQ^c\--coalg \to \CC^c\--coalg .$$

So now we have all the ingredients to write a structure theorem in the dual case.

Observe that $U^c(C)$ acquires a $\CC^c\--\AA$-bialgebra structure from the $\CC^c\--\AA$-bialgebra structure of $\CC^c(V)$. 

The dual PBW has been proved in \cite{Michaelis85}. See also \cite{Fresse06} \S 4.2 for the dual PBW and dual CMM theorems. The Eulerian idempotent has been worked out in this context by M.~Hoffman \cite{Hoffman00}.

\section{Generalized bialgebras in characteristic $p$}\label{charp}

First, observe that in the nonsymmetric framework (cf.~\ref{nscase}) there is no characteristic assumption, therefore the structure theorem holds in characteristic $p$. In \cite{Cartier} and \cite{MM} the authors give a characteristic $p$ version  of the PBW theorem and of the CMM theorem. So there is a characteristic $p$ version of the structure theorem for cocommutative (classical) bialgebras. But the notion of Lie algebra has to be replaced by the notion of $p$-restricted Lie algebras. 

\subsection{$p$-restricted Lie algebras} By definition a  \emph{$p$-restricted Lie algebra} is a Lie algebra over a characteristic $p$ field which is equipped with a unary operation $x\mapsto x^{[p]}$ called the Frobenius operation. It is supposed to satisfy all the formal properties of the iterated bracket 
$$\underbrace{[x, [x, [\  \cdots , [x}_{p\ \textrm {times}},-]\cdots ]]]$$
in an associative algebra (cf.~loc.cit.).

In the PBW and CMM theorems the forgetful functor is replaced by the functor $As\alg\to p\--Lie\alg$ where the bracket is as usual and the Frobenius is the iterated bracket as above. Since this functor admits a left adjoint $U$ all the ingredients are in place for a structure theorem in that case 
 (cf.~loc.cit.).

\subsection{Operads in characteristic $p$} Since the Frobenius is not a linear operation (it is polynomial of degree $p$), a $p$-restricted Lie algebra is not an algebra over some operad in the sense of \ref{algoverop}. Note that, in our definition of operad, we defined the Schur functor $\PP(V) := \bigoplus_{n}\PP(n)\t_{S_{n}}V^{\t n}$ by using the coinvariants. If, instead, we had taken the invariants, then there would be no difference in characteristic 0, but it would be different in characteristic $p$. In short, there is a way to handle $p$-restricted Lie algebras in the operad framework by playing with the two different kinds of Schur functors. In fact B.~Fresse showed in \cite{Fresse00} how to work with any operad in characteristic $p$ along this line. It gives rise to the notion of $\PP$-algebra with \emph{divided symmetries}\index{divided symmetries}. For instance, in the commutative case, it gives rise to the divided power operation. If the operad is nonsymmetric the notions of $\PP$-algebra and  $\PP$-algebra with divided symmetries are equivalent.

\subsection{Structure theorem in characteristic $p$} Now we have all the ingredients to write down a structure theorem for generalized bialgebras in characteristic $p$, including a toy-model. I conjecture that such a theorem exists. In fact some cases, with $\CC=As$ or $Com$, have already been proved, see  \cite{Patras99} and the references in this paper.

Observe that there are two levels of difficulty. First write the general theorem and its proof, second handle explicit cases. Recall for instance that, for Poisson algebras, we have to work with two  divided operations: the Frobenius operation and the divided power operation. The relationship between these two are quite complicated formulas (cf. \cite{Fresse06}).

\section{Relationship with rewriting systems}\label{rewriting} The rewriting theory aims at computing a monoid (or a group) starting from a presentation. The idea is to write any relation under the form $u_{1} \ldots u_{k}= v_{1}\ldots v_{l}$ and to think of it as a ``rewriting procedure" $u_{1} \ldots u_{k}\to  v_{1}\ldots v_{l}$. In this setting one can define the notions of  noetherianity, confluence, critical peak and convergence. There is a way to extend the rewriting theory to operads and even props, see Y.~Lafont \cite{Lafont03}. For instance a distributive compatibility relation like $\dd \circ \mu = \Phi$ (cf.~\ref{comprel}) can be thought of as a rewriting procedure $\dd \circ \mu \to \Phi$. The aim is to find a \emph{reduced form} for the multivalued operations. In this setting Hypothesis \Hone can be proved by verifying that a rewriting system is convergent. See \ref{IC} for an example taken out of \cite{Lafont97}.

\section{Application to representation theory}\label{representationtheory} Given an \sm $\AA$ and a sub \sm $\PP$ it is usually difficult to decide whether there exists an \sm $\ZZZ$ such that $\AA = \ZZZ \circ \PP$ (recall that in this framework the composition $\circ$ is called the \emph{plethysm}\index{plethysm}). We will show that, in certain cases, we can give a positive answer to this question.

\begin{prop} Let $\AA$ be an operad and let $\PP$ be a suboperad of $\AA$. The following condition is sufficient to ensure that there is an isomorphism of {\bf S}-modules $\AA\cong \ZZZ \circ \PP$, where $\ZZZ:= \AA/(\overline \PP)$.

Condition: there exists an operad $\CC$ and a good triple of operads $(\CC, \AA, \PP)$ giving rise to the inclusion $\PP \subset \AA$.
\end{prop}
\begin{proo} If $(\CC, \AA, \PP)$ is good, then so is $(\CC, \AA / (\overline \PP), Vect)= (\CC, \ZZZ, Vect)$ by Proposition \ref{quotient}. So we get isomorphisms of {\bf S}-modules $\CC^c \cong \ZZZ$ and $\AA\cong \CC^c \circ \PP$ (Proposition \ref{isosms}), which imply $\AA\cong \ZZZ \circ \PP$.
\end{proo}



\chapter{Examples}\label{ch:examples}

The problem of determining if a triple of operads \CAP\!, or more accurately $(\CC, \between, \AA, F, \PP)$, is good may crop up in different guises. Most of the time the starting data is the prop $(\CC, \between, \AA)$, that is the type of bialgebras. Verifying {\texttt{(H0)}} is, most of the time, immediate by direct inspection. The first problem is to verify the hypotheses \Hone and \Htwoepi. The second problem (and often the most difficult) is to find a small presentation of the operad $\PP= \Prim_{\CC} \AA$ and make explicit the functor $F: \AA\alg \to \PP\alg$.

Another kind of problem is to start with a forgetful functor $F: \AA\alg \to \PP\alg$ (i.e.~$\PP$ is a suboperad of $\AA$) and to try to find $\CC$ and $\between$ so that the data 
 $(\CC, \between,\AA, F, \PP)$ is a good triple.
 
In both problems Corollary \ref{cor:genser} relating the generating series of $\CC, \AA$ and $\PP$ is a good criterion since the knowledge of two of the operads determines uniquely the generating series of the third.

As said in the introduction the (uni)versal idempotent $e$ is a very powerful tool. In Chapter 3 it is constructed abstractly. To get it explicitly in a given example is often a challenge. 

In this chapter we present several concrete cases. For many of them, existing results in the literature permit us to prove the hypotheses and to find a small presentation of the primitive operad. In some cases the technique is very close to the \emph{rewriting techniques} in computer sciences. Proposition \ref{prop:quotient} is quite helpful in proving that $(\CC, \AA, \Prim_{\CC}\AA)$ is a good triple since it reduces several cases to one case. For instance when $\AA$ is generated by one operation, it suffices to prove the hypotheses for the prop $\CC^c\--Mag$.

We have seen that any good triple \CAP gives rise to a triple of the form $(\CC, \ZZZ, Vect)$ (the quotient triple) by moding out by the primitive operad $\PP$. We put in the same section the triples which have the same quotient triple (with a few exceptions).

We give only the proofs of the statements which are not already in the literature.
\bigskip

In section 1 we treat $(Com, Com, Vect)$ with Hopf compatibility relation. It includes the classical case $(Com, As, Lie)$ as well as $(Com, Parastat, Nil)$ and $(Com, Mag, Sabinin)$. The rigidity theorem for $(Com, Com, Vect)$ is the Hopf-Borel theorem and the structure theorem for $(Com, As, Lie)$ is the union of the CMM theorem and the PBW theorem. We prove that, in this classical case, the universal idempotent is precisely the well-known Eulerian idempotent.

In section 2 we treat $(As,As,Vect)$ with nonunital infinitesimal compatibility relation. It includes the case $(As, Mag, MagFine)$ and the case $(As, 2as, MB )$ which is important because the category of cofree Hopf algebras is equivalent to the category of $MB $-algebras. The triple $(As, Dup, Mag)$ where a $Dup$-algebra is a space equipped with two associative operations satisfying further the relation 
$$(x\d y)\g z = x\d (y\g z),$$
should be in this section. It will be treated in full detail in the next Chapter.

In section 3 we treat $(As, Zinb, Vect)$ with semi-Hopf compatibility relation. It includes the case $(As, Dipt, MB )$ and the case $(As, Dend, brace)$ (due to Mar\'\i a Ronco \cite{Ronco02}) which is important since it permits us to unravel the structure of a free brace algebra.

In section 4 we treat  $(Lie, Lie, Vect)$. It should be noted that the notion of $Lie^c\--Lie$-bialgebra is NOT what is commonly called Lie bialgebras because the compatibility relation is different. In particular there is a nontrivial $\Phi_{1}$ term in our case (cf.~\ref {diagram}). 

In section 5 we treat $(NAP, PreLie, Vect)$ due to M. Livernet \cite{Livernet06} and the triple $(NAP, Mag, \Prim_{NAP}Mag)$.

In section 6 we describe several cases of the form $(\AA, \AA, Vect)$. 

In section 7 we describe the interchange bialgebra case. Here the operads are no more quadratic but cubic. 

In section 8 we treat a case where the generating operations and cooperations are of arity $k$. 

When there is no $\Phi_{1}$ term in the compatibility relation(s) (see \ref{diagram}), every operation is a primitive operation and there is nothing to prove. This is why we do not treat the Frobenius case.

 \section{Hopf algebras: the classical case}\label{ex:comas} In this section we treat several triples admitting the triple $(Com,Com, Vect)$ (Hopf-Borel) as a quotient triple. It includes the classical case $(Com, As, Lie)$. The compatibility relation in these cases is the Hopf relation. More examples can be found in the subsequent paper \cite[LRcha].
 
 \subsection{The Hopf compatibility relation}\label{Hopf} First, let us recall some elementary facts about unital associative algebras. The tensor product $A\t B$ of the two unital associative algebras $A$ and $B$ is itself a unital associative algebra with product given by  $(a\t b) (a'\t b') = aa' \t bb'$ and with unit $1_{A}\t 1_{B}$. The free unital associative algebra over $V$ is the tensor algebra 
$$T(V) = \KK \oplus V \oplus \cdots \oplus V^{\t n} \oplus \cdots $$
whose product is the concatenation:
$$(v_{1}\cdots v_{p})(v_{p+1}\cdots v_{n}) = v_{1}\cdots v_{n}.$$
 Let $V\to T(V)\t T(V)$ be the map given by $v\mapsto v\t 1 + 1 \t v$. Since $T(V)$ is free, there is a unique extension as algebra homomorphism denoted
$$\DD : T(V) \longrightarrow T(V) \t T(V)\ .$$
It is easy to show, from the universal property of the free algebra,  that $\DD$ is coassociative and cocommutative. The fact that $\DD$ is an algebra morphism reads
$$\DD(xy) = \DD(x) \DD(y)\ ,$$
which is the \emph{classical Hopf compatibility relation}. Hence the tensor algebra, equipped with this comultiplication, is a classical cocommutative bialgebra. The map $\Delta$ can be made explicit in terms of shuffles, cf. \cite{JLLeuler}.

In order to work in a non-unital framework, we need to restrict ourselves to the augmentation ideal of the bialgebra and to introduce the \emph{reduced coproduct} $\dd$
$$\dd(x):= \DD(x)-x\t 1 - 1\t x\ .$$
As already mentioned (cf.~\ref{ex1}) the compatibility relation between the product $\mu$ and the coproduct  $\dd$ becomes $\between_{Hopf}$:
$$\dd(xy) =x\t y + y \t x + \dd(x) (y\t 1 + 1 \t y) + (x\t 1 + 1 \t x) \dd(y) + \dd(x) \dd (y),$$
where $xy= \mu(x,y)$.
Diagrammatically it reads:
$$\cpbdeuxdeux = \cpbA +\cpbB+\cpbC +\cpbD +\cpbE +\cpbF +\cpbG .$$
Observe that this is a distributive compatibility relation.

\subsection{The triple $(Com, As, Lie)$}\label{ComAsLie} By definition a $Com^c\--As$-bialgebra is (in the non-unital framework) a vector space $\HH$ equipped with a (nonunital) associative operation  $\mu$, a commutative associative comultiplication $\dd$, satisfying the Hopf compatibility relation $\between_{Hopf}$. Obviously hypothesis {\texttt{(H0)}} is fulfilled.
As we already mentioned in \ref{ex1} a $Com^c\--As$-bialgebra is equivalent to a classical bialgebra\footnote{Here we only deal with connected bialgebras for which an antipode always exists. So there is an equivalence between connected bialgebras and Hopf algebras.} by the map $\HH \mapsto \HH_{+} := \KK 1\oplus \HH$.

The free $As$-algebra\index{free associative algebra} over $V$ is the reduced tensor algebra 
$$As(V)=\To(V)= V \oplus \cdots \oplus V^{\t n} \oplus \cdots $$
 equipped with the concatenation product. From the property of the tensor algebra recalled above, it follows that $\To(V)$ is a cocommutative bialgebra, in other words it satisfies the hypothesis \texttt{(H1)}. We claim that the operad $\PP$ deduced from Theorem \ref{thm:structure} is the operad $Lie$ of Lie algebras. Indeed it is well-known that $\Prim \To(V) = Lie(V)$, cf. for instance \cite{Wigner} for a short proof. 

The map $\varphi : As \to Com^c$ is given by $x_{1}\cdots x_{n} \mapsto x_{1}\cdots x_{n}$ in degree $n$, where on the left side we have a noncommutative polynomial, and on the right side we have a commutative polynomial. In other words $\varphi_{n}:As(n)\to Com^c(n)$ is the map $\KK[S_{n}] \to \KK, \ss \mapsto 1_{\KK}$.
This map has a splitting in characteristic zero, given by $x_{1}\cdots x_{n} \mapsto \frac{1}{n!}\sum_{\ss\in S_{n}} \ss(x_{1}\cdots x_{n})$. It is a coalgebra morphism for the coalgebra structure of $\To(V) = As(V)$ constructed above. Hence the hypothesis \texttt{(H2epi)} is fulfilled. So the triple $(Com, As, Lie)$ is a good triple of operads and the structure Theorem holds for this triple. Translated in terms of unital-counital cocommutative bialgebras it gives:

\begin{thm}[CMM+PBW] \label{thm:ComAs}\index{CMM}\index{PBW}  Let $\KK$ be a field of characteristic zero.
For any cocommutative bialgebra $\HH$ the following are equivalent:

\noindent a)  $\HH$ is connected,

\noindent b) there is an isomorphism of bialgebras $\HH \cong U(\Prim \HH)$,

\noindent c) there is an isomorphism of connected coalgebras $\HH \cong S^c(\Prim \HH)$.
\end{thm}

Here the functor $U$ is the classical universal enveloping algebra functor from the category of Lie algebras to the category of unital associative algebras (or more accurately classical cocommutative bialgebras).

Of course, this is a classical result. In fact $(a) \Rightarrow (b)$ is the Cartier-Milnor-Moore Theorem  which first appeared in Pierre Cartier's seminar lectures \cite{Cartier}\footnote{In this paper a bialgebra is called a ``hyperalgebra".} and was later popularized by Milnor and Moore in \cite{MM}.

Then  $(b) \Rightarrow (c)$ is, essentially,  the Poincar\'e-Birkhoff-Witt Theorem. In fact it is slightly stronger since it not only gives a basis of $U(\gg)$ from a basis of the Lie algebra $\gg$ but it also provides an isomorphism of \emph{coalgebras} $U(\gg) \to S^c(\gg)$. This is Quillen's version of the PBW theorem, cf. \cite{Quillen} Appendix B. In this appendix Dan Quillen gives a concise proof of the PBW Theorem and of the CMM Theorem. The idempotent that he uses in his proof is the \emph{Dynkin idempotent}\index{Dynkin idempotent} 
$$x_{1} \ldots x_{n}\mapsto \frac{1}{n} [\ldots [[x_{1},x_{2}],x_{3}],\ldots , x_{n}].$$

We will show that the idempotent given by our proof (cf.~\ref{idempotent}) is the Eulerian idempotent (cf.\cite{Reutenauer, JLLeuler}). 

Theorem \ref{classicalPBW} applied to the triple $(Com,As,Lie)$ gives the most common version of the PBW Theorem:
$$ gr\ U(\gg) \cong S(\gg).$$
Observe that the implication  $(a) \Rightarrow (c)$ had been proved earlier by Jean Leray (cf.~\cite {Leray}), who had shown that the associativity hypothesis of the product was not necessary for this implication (see \ref{Leray} for an explanation in terms of triples of operads). 

For historical notes on the PBW theorem 
 one may consult the paper by P.-P. Grivel \cite{Grivel}.

\subsection{Eulerian idempotent}\label{Euler} \cite{JLLeuler}  Let $\HH$ be a connected cocommutative bialgebra (nonunital framework). The convolution of two linear maps $f,g: \HH \to \HH$ is defined as 
$$f\star g := \mu \circ (f\t g) \circ \dd\ .$$
 It is known that $\dd$ can be made explicit in terms of shuffles. By definition the (first) \emph{Eulerian idempotent}\index{Eulerian idempotent} $e^{(1)}: \HH \to \HH$ is defined as
$$e^{(1)}:=\log^{\star} (uc+J)= J - \frac{J^{\star 2}}{2} + \frac{J^{\star 3}}{3} - \ldots $$
where $J := \Id_{\HH}-uc$. For $\HH= \overline T(V)$, the nonunital tensor algebra, $e^{(1)}$ sends $V^{\t n}$ to itself and we denote by $e_{n}^{(1)}:V^{\t n}\to V^{\t n}$ the restriction to $V^{\t n}$. Explicitly, it is completely determined by an element $e_{n}^{(1)}= \sum_{\ss} a_{\ss}\ss \in \QQ[S_{n}]$ since, by the Schur Lemma, 
$$e_{n}^{(1)}\row x1n =  \sum_{\ss} a_{\ss}(x_{\ss(1)}\cdots x_{\ss(n)}),$$
for some coefficients $a(\ss)$ (here we let $\ss$ act on the right).

The higher Eulerian idempotents are defined as
$$e^{(i)} := \frac {(e^{(1)})^{\star i}}{i!}\ .$$
From the relationship between the exponential series and the logarithm series, it comes:
$$\Id_{n} = e_{n} ^{(1)}+ \cdots +e_{n} ^{(n)} \ . $$
\begin{prop} For any connected cocommutative bialgebra $\HH$ the \emph{universal idempotent} $e$ is equal to the Eulerian idempotent:
$$e := \Pi_{n\ge 2} (\Id -\oo^{[n]})= \sum _{n\ge 1}(-1)^n \frac{J^{\star n}}{n} =: e^{(1)}\ . $$
\end{prop}

\begin{proo} It suffices to prove this equality for $\HH = \overline T(V)$. From the definition of $\oo^{[n]}$ we get its expression in terms of shuffles. We get 
$$\oo^{[n]} = e^{(n)}+ e^{(n+1)}+ \cdots \ .$$
Hence we deduce
$$\Id -\oo^{[n]} = e^{(1)}+ \cdots + e^{(n-1)}\ .$$
Since the idempotents $e^{(i)}$ are orthogonal to each other (cf.~\cite{JLLeuler}) we get
$$e = \Pi_{n\ge 2} (\Id -\oo^{[n]}) = e^{(1)}(e^{(1)}+e^{(2)})(e^{(1)}+e^{(2)}+e^{(3)})( \cdots )\cdots = e^{(1)}.$$
\end{proo}

\subsection{Explicit formula for the PBW isomorphism}\label{explicit} Since the Eulerian idempotent can be computed explicitly in the symmetric group algebra, one can give explicit formulas for the isomorphism
$$T(V) \cong S^c(Lie(V))\ .$$
In low dimension we get:
$$x=x,$$
$$xy = \frac{1}{2}[x,y] + \frac{1}{2}(xy+yx),$$
\begin{eqnarray*}
xyz&=& \frac{1}{6}([[x,y],z] + [x,[y,z]])\\
&& + \frac{1}{4}\big(x[y,z] + [y,z]x+ y[x,z] + [x,z]y+ z[x,y] + [x,y]z\big)\\
&& + \frac{1}{6}\sum_{\ss\in S_{3}}\ss (xyz)\ . 
 \end{eqnarray*}
\subsection{Remark} In the case of classical bialgebras, not necessarily cocommutative, i.e.~$As^c\--As$-bialgebra with $\between= \between_{Hopf}$, the hypotheses {\texttt{(H0)}} and \Hone are also fulfilled. However the condition \Htwoepi is not fulfilled, since the map $\varphi : \To(V) \to \To^c(V)$ factors through $\overline{S}^c(V)$. This is due to the cocommutativity of the coproduct on the free associative algebra.

\subsection{The triple $(Com,Com,Vect)$}\label{comcom} As mentioned in \ref{quotient} if we mod out by relators in $Lie$, then we get a new triple of operads. For instance if we mod out by $\overline {Lie}$ (cf.~\ref{opideal} for the notation), then we get a good triple of operad:
$$(Com,Com,Vect).$$

In the unital framework the free commutative algebra $Com(V)$ is the symmetric algebra $S(V)$, which is the polynomial algebra $\KK[x_{1}, \ldots , x_{r}]$ when $V= \KK x_{1}\oplus \ldots \oplus\KK x_{r}$. Similarly the cofree coalgebra $Com^c(V)$ can be identified with $\KK[x_{1}, \ldots , x_{r}]$ and the coproduct is given by $\DD(x_{i}) = x_{i}\t 1 + 1 \t x_{i}$. Under these identifications the map 
$\varphi(V) : Com(V) \to Com^c(V)$ is not the identity, but is given by $x_{1}^{i_{1}}\cdots x_{r}^{i_{r}}\mapsto \frac{x_{1}^{i_{1}}}{i_{1}!}\cdots \frac{x_{r}^{i_{r}}}{i_{r}!}$. This phenomenon can be phrased differently as follows. On the vector space of polynomials in one variable one can put two different commutative algebra structures:
$$(I)\qquad x^px^q := x^{p+q},$$
$$(II)\qquad x^px^q := {{p+q}\choose{p}}x^{p+q},$$
where ${p+q}\choose{p}$ is the binomial coefficient $\frac{(p+q)!}{p! q!}$. Of course, over a characteristic zero field, they are isomorphic ($x^n\mapsto \frac{x^n}{n!}$). By dualization we obtain two coalgebra structures $(I^c)$ and $(II^c)$. In order to make  $\KK[x]$ into a Hopf algebra we have to combine either $(I)$ and $(II^c)$ or $(II)$ and $(I^c)$. 

The rigidity theorem for the cocommutative commutative connected bialgebras is the classical Hopf-Borel theorem recalled in the introduction of \ref{s:rigidity}. Let us recall that the classical version (the one which is used in algebraic topology) is phrased in the graded framework (cf.~\ref{graded}). Here we gave the claim in the nongraded framework.

\subsection{The triple $(Com, Parastat, NLie)$} A \emph{parastatistics algebra}\index{parastatistics algebra} is an associative algebra for which the Lie bracket verifies the relation
$$[[x,y],z]=0.$$
In a classical bialgebra the elements $[[x,y],z]$ are primitive, hence we are in the situation of Proposition \ref{quotient}.  The primitive type associated to the bialgebra type $As^c\pt Parastat$ is simply nilpotent algebras whose product is antisymmetric. The structure theorem was proved in \cite{LPo} (it follows easily from the classical one). This triple is interesting on several grounds. First, the parastatistics algebras (and their sign-graded version) appear naturally in theoretical physics. Second, the parastatistics operad is interesting from a representation theory point of view because $Parastat(n)$ is the sum of  \emph{one} copy of each irreducible type of $S_{n}$-representations. Third, the parastatistics algebras are exactly the algebras which are at the same time associative and Poisson (for the Lie bracket and the symmetrized operation).

\subsection{The triple $(Com, Mag, Sabinin)$}\label{Sabinin}

Let $Com^c\--Mag$ be the magmatic cocommutative bialgebra type. The compatibility relation is the Hopf compatibility relation. Let us recall that a \emph{magmatic algebra}\index{magmatic algebra} is a vector space equipped with a binary operation, without further hypothesis. This is the nonunital case. In the unital case we assume further that there is an element $1$ which is a unit on both sides. It is easy to show that the free magmatic algebra can be equipped with  a cocommutative cooperation as follows. 
Working in the unital framework we put on the tensor product of two unital magmatic algebras a  unital magmatic structure by
$$(a\t b)\cdot (a'\t b') = (a\cdot a')\t (b\cdot b')\ \textrm{and}\ 1_{A\t B}= 1_{A}\t 1_{B}.$$
So the free unital magmatic algebra $Mag_{+}(V)$ tensored with itself is still a unital magmatic algebra. There is a unique morphism 
$Mag_{+}(V)\to Mag_{+}(V)\t Mag_{+}(V)$
which extends the map 
$$V\to Mag_{+}(V)\t Mag_{+}(V), v\mapsto v\t 1 + 1 \t v.$$
 This cooperation is immediately seen to be coassociative and cocommutative. Restricting the whole structure to the augmentation ideal gives a $Com^c\--Mag$ structure on the free magmatic algebra $Mag(V)$. Hence hypothesis \texttt{(H1)} holds. Explicitly the cooperation $\dd$ is given by
 $$\dd(t;v_{1},\ldots , v_{n}) = \sum_{i=1}^{n-1}\sum_{\ss}(t_{(1)}^{\ss}; v_{\ss(1)}, \ldots, v_{\ss(i)})\t 
 (t_{(2)}^{\ss}; v_{\ss(i+1)}, \ldots, v_{\ss(n)})$$
 where  $\ss$ is an $(i,n-i)$-shuffle and the trees $t_{(1)}^{\ss}$ and $t_{(2)}^{\ss}$ are subtrees of  $t$ corresponding to the shuffle decomposition (cf.~\cite{Holtkamp05}).
 
 A coalgebra splitting $s$ to $\varphi(V)$ is obtained by
 $$s(v_{1}\ldots v_{n})= \sum_{\ss\in S_{n}}\frac{1}{n!}(\textrm{comb}^l_{n}; v_{\ss(1)}\ldots v_{\ss(n)})$$
 where $\textrm{comb}^l_{n}$ is the left comb (cf.~\ref{pbtree}). So
 Hypothesis \texttt{(H2epi)} is fulfilled. Hence the structure theorem holds for 
$Com^c\--Mag$-bialgebras. It was first proved by R.~Holtkamp in \cite{Holtkamp05}. Earlier studies on this case can be found in the pioneering work of M.~Lazard \cite{Lazard55} in terms of ``analyseurs" and also in \cite{GH}.

\subsection{Sabinin algebras} The problem is to determine explicitly the primitive operad $\Prim_{Com}Mag$. Results of Shestakov and Urmibaev \cite{SU} and of P\'erez-Izquierdo \cite{Perez06} show that it is the \emph{Sabinin} operad. A \emph{Sabinin algebra}\index{Sabinin algebra} can be defined as follows (there are other presentations). The generating operations are: 
$$\langle x_1, \ldots , x_m; y,z\rangle,\quad  m\geq 0,$$
$$\Phi(x_1,\ldots , x_m ; y_1,\ldots , y_n), \quad m\geq 1, n\geq 2,$$
 with symmetry relations 
$$\langle x_1, \ldots , x_m; y,z\rangle=-\langle x_1, \ldots , x_m; z,y\rangle,$$
$$\Phi(x_1,\ldots , x_m ; y_1,\ldots , y_n) = \Phi(\omega(x_1,\ldots , x_m) ; \theta(y_1,\ldots , y_n)), \quad \omega\in S_m, \theta\in S_n,$$
 and the relations are
$$\langle x_1, \ldots , x_r,u,v,x_{r+1},\ldots , x_m; y,z\rangle-\langle x_1, \ldots , x_r,v,u,x_{r+1},\ldots , x_m; y,z\rangle$$
$$+ \sum_{k=0}^r\sum_{\ss}\langle x_{\ss(1)}, \ldots , x_{\ss(k)}; \langle x_{\ss(k+1)},\ldots , x_{\ss(r)};u,v\rangle, x_{r+1}, \ldots , x_m; y,z\rangle$$
where $\ss$ is a $(k,r-k)$-shuffle,
$$K_{x,y,z}\big(\langle x_1, \ldots , x_r,x; y,z\rangle+\hfill  $$
$$\hfill \sum_{k=0}^r\sum_{\ss}\langle x_{\ss(1)}, \ldots , x_{\ss(k)}; \langle x_{\ss(k+1)},\ldots , x_{\sigma(r)};y,z\rangle, x\rangle\big)=0$$
where $K_{x,y,z}$ is the sum over all cyclic permutations. 

Observe that there is no relation between the generators $\langle -;-\rangle$ and the generators $\Phi$. The functor $F:Mag\alg \to Sabinin\alg$ is constructed explicitly in \cite{SU} (also recalled in \cite{Perez06}). For instance $\langle y,z\rangle = -y\cdot z + z \cdot y$ and $(x\cdot y)\cdot z - x\cdot (y\cdot z) = -\frac{1}{2}\langle x;y, z \rangle +  \Phi (x;y,z)$. It is easy to check that in $Mag(V)$ the two operations ``bracket" and ``associator" are not independent but related by the \emph{nonassociative Jacobi identity}\index{magmatic Jacobi identity}(cf.~\cite{Holtkamp05}):
$$[[x,y],z]+[[y,z],x]+[[z,x],y] = \sum_{\textrm{sgn}(\ss)\ss\in S_3} \ss\  as(x,y,z)\ .$$
In the preceding presentation it corresponds to the cyclic relation with $r=0$.

So, there is a good triple of operads 
$$(Com, Mag, Sabinin)\ .$$
  As a consequence, the generating series is $ f^{Sab}(t)=\log(1+  (1/2)(1-\sqrt{1-4t}))$ and the dimension of the space of $n$-ary operations is
$$\dim Sab(n) =1,1,8,78, 1104, \ldots\  .$$

\subsection{Remarks}\label{Leray} Since the $Com^c\--Mag$-bialgebra type satisfies {\texttt{(H0)}}, \Hone and \Htwoepi, any connected bialgebra is cofree. This result has been proved earlier by Jean Leray in \cite {Leray}. See also \cite {Oudom99, Fresse98} for a different generalization.

\subsection{Quotients of $Com^c\--Mag$} The classical type $(Com, As, Lie)$ is a quotient of the triple $(Com, Mag, Sabinin)$ (quotient by the associator, which is a primitive element and apply Proposition \ref{quotient}). In fact we have the following commutative diagram of operads:
$${\xymatrix{
Sabinin\ar@{->>}[r]\ar[d] & Lie \ar@{->>}[r]\ar[d] &  NLie \ar@{->>}[r]\ar[d]& Vect\ar[d]  \\
Mag\ar@{->>}[r] & As \ar@{->>}[r] & Parastat \ar@{->>}[r] & Com \\
}}$$

It may be worth to study other quotients of $Mag$ by an ideal $J$ (for instance $PreLie$ since $\langle x;y,z\rangle$ is the pre-Lie relator),  and find a small presentation of the quotient $Sabinin / (J)$. Some results in this direction have been done for $Malcev$-algebras\index{Malcev algebra} in \cite{PIS}. It fits into this framework, since the Malcev relators:
$$x\cdot y-y\cdot x$$
$$((x\cdot y)\cdot z) \cdot t + (x \cdot (y \cdot z)) \cdot t -  (x \cdot y) \cdot (z \cdot t) +  x \cdot ((y \cdot z) \cdot t) +  x \cdot (y \cdot( z \cdot t))$$
are primitive in $Mag(V)$. Other examples with Moufang algebras, Bol algebras and Lie triple systems should come into play.

\subsection{Poisson bialgebras} Let us mention that there is a notion of Hopf-Poisson algebras (see for instance \cite{FM}), that is $As^c\--Pois$-bialgebra and also cocommutative Hopf-Poisson algebras, that is $Com^c\--Pois$-bialgebra.
The compatibility relation for the pair $(\dd, \cdot)$ (where $a\cdot b$ is the commutative operation) is Hopf and the compatibility relation for the pair $(\dd, [\,  ])$ (where $[\,  ]$ is the Lie bracket)  is given by

$$\cpbDeuxDeux{[\, ]}{*{}} = 
\cpbGG{*{}}{[\, ]}{*{}}{\cdot}+\cpbGG{*{}}{\cdot}{*{}}{[\,]}\ .$$

As in the classical case there is a good triple $(Com, Pois, Lie)$ since it is well-known that the free Poisson algebra $Pois(V)$ is precisely $Com(Lie(V))$.

\subsection{A conjectural triple $(Com, ?? , preLie)$}\footnote{Note added in proof. This problem has been settled in \cite{LRcha}.} In \cite{CL} F.\ Chapoton and M.\ Livernet showed that the symmetric algebra over the free pre-Lie algebra in one generator can be identified with the dual of the Connes-Kreimer Hopf algebra \cite{CK}. The study of this case, as done in \cite{GO} for instance, suggests the existence of a triple of the form $(Com, \AA , preLie)$ where the operad $\AA$ is the unknown. It is expected that $\AA(\KK)\cong \HH_{CK}^*$ and that this triple admits $(Com, As, Lie)$ as a  quotient. It is also interesting to remark that $\dim \AA(n)$ should be equal to $(n+1)^{n-1}$ which is also $\dim Park(n)$ (parking functions).

\section{Unital infinitesimal bialgebras}\label{uib} In this section we study some triples $(As, \AA, \PP)$ which are over $(As, As, Vect)$ with compatibility relation the nonunital infinitesimal relation.

\subsection{The (non)unital infinitesimal compatibility relation} On the tensor algebra $T(V)$ the product is the concatenation product $\mu$. Let us equip it with the deconcatenation coproduct given by
$$\DD(v_{1}\cdots v_{n}):=\sum_{i=0}^{n} v_{1}\cdots v_{i}\t v_{i+1}\cdots v_{n} .$$
The pair $(\DD,\mu)$ does not satisfy the Hopf compatibility relation, but does satisfy the \emph{unital infinitesimal}\index{unital infinitesimal} relation:
$$\DD(xy) = -x\t y + x_{(1)}\t x_{(2)}y + xy_{(1)}\t y_{(2)} ,$$
where $\DD(x)= x_{(1)}\t x_{(2)}$.

Since we want to work in the nonunital framework, we need to introduce the reduced coproduct $\dd$ defined by the equality $\DD(x)= x\t 1 + 1\t x + \dd(x)$.

The compatibility relation for the pair $(\dd,\mu)$ is the \emph{nonunital infinitesimal}\index{nonunital infinitesimal} (n.u.i.) compatibility relation $\between_{nui}$:
$$\dd(xy) = x\t y + \dd(x)(1\t y ) + (x\t 1) \dd(y),$$
Diagrammatically it reads (cf.  \ref{genbialg} example 2):
$$\cpbdeuxdeux = \cpbA +\cpbC +\cpbE\quad .$$

\subsection{The triple $(As, As, Vect)$}\label{AsAsVect} By definition a \emph{nonunital infinitesimal bialgebra} ($As^c\--As$-bialgebra) is a vector space equipped with an associative operation and a coassociative cooperation satisfying the n.u.i.~compatibility relation. 
Hypothesis {\texttt{(H0)}} is obviously fulfilled. From the above discussion it follows that  hypothesis \Hone is also fulfilled. 

The next Proposition shows that hypothesis \Htwoiso is fulfilled.
So we get the rigidity theorem for the triple $(As,As,Vect)$. It was first announced in \cite{LRnote} and proved in \cite{LRstr} where details can be found. Let us just recall that the universal idempotent in this case is given by the geometric series:
$$e = \sum_{n\geq 1}(-1)^{n-1} \id^{\star n}$$
where $\star$ is the convolution product $f\star g := \mu \circ (f\t g) \circ \dd$.

\begin{prop} The map $\varphi(V):\To(V) \to \To(V)^c$ is induced by the identification of the generator of $As_{n}$ with its dual. 
\end{prop}
\begin{proo} The prop $As^c\pt As$ is nonsymmetric and $As_{n}$ is one-dimensional. Let $\mu_{n}$ denote the generator of $As_{n}$, so $\mu_{n}\row x1n = x_{1}\ldots x_{n}$. In order to compute its image by $\varphi$ it suffices to compute $\dd_{n}\circ \mu_{n}\row x1n$ where $\dd_{n}$ is the dual of $\mu_{n}$, that is the generator of $As_{n}^c$ (cf. \ref{verification}). From the compatibility relation we get
$$ \dd_{n}\circ \mu_{n}\row x1n = x_{1}\t \cdots \t x_{n}\in V^{\t n}\subset As(V)^{\t n}$$
where $V=\KK x_{1}\oplus \cdots \oplus \KK x_{n}$. Hence $\varphi(V)(\mu_{n})=\dd_{n}$ and $\varphi$ is an isomorphism.
\end{proo}

\subsection{The triple $(As, Mag, MagFine)$}\label{MagFine} This triple has been studied and proved to be good  in \cite{HLR}. The compatibility relation is $\between_{ui}$. The coproduct $\dd$ on $Mag(V)= \oplus_{n\geq 1}\KK[PBT_{n}]\t V^{\t n}$, where $PBT_{n+1}$ is the set of planar binary trees with $n$ internal vertices (cf.~\ref{pbtree}), can be constructed as follows. Let $t$ be p.b.~trees whose leaves are numbered from left to right beginning at $0$. We cut the tree along  the path going from the $i$th vertex (standing between the leaves $i-1$ and $i$) to the root. It gives two trees denoted $t^{i}_{(1)}$ and $t^{i}_{(2)}$.We have
$$\dd(t; v_{0}\cdots v_{n}) = \sum_{i=1}^{n}(t^i_{(1)};v_{0}\cdots v_{i-1})\t (t^i_{(2)};v_{i}\cdots v_{n})\ .$$
Hypothesis \Hone can be proved either by using the explicit form of $\dd$ or by an inductive argument as explained in \ref{verification}.
The map $\varphi(V): Mag(V)\to As^c(V)$ sends $t$ to the generator $1_{n}$ of $As_{n}^c$. We choose the splitting $s_{n}: As_{n}\to Mag_{n}$ given by $s(1_{n})= comb^l_{n}$ (left comb). So Hypothesis \Htwoepi is fulfilled and we have a good triple of operads:
$$(As, Mag, \Prim_{As}Mag)\ .$$
It is proved in \cite {HLR} that the primitive operad $\Prim_{As}Mag$ is generated by $n-2$ operations in arity $n$ and that they satisfy no relations. So this operad is a magmatic operad (free operad). Since the dimension of $(\Prim_{As}Mag)_{n}$ is the Fine number, it is called the \emph{magmatic Fine operad}\index{magmatic Fine operad}, denoted $MagFine$. So there is a good triple of operads
$$(As, Mag, MagFine)\ .$$
\subsubsection{Relationship with previous work} As a byproduct of the structure theorem for $As^c\--Mag$-bialgebras we have that a connected coassociative algebra equipped with a magmatic operation satisfying the n.u.i.~relation is cofree. Dually we have the following: an associative algebra equipped with a comagmatic operation, which is connected and satisfies the n.u.i.~relation is free. A very similar result has been shown by I.~Berstein in \cite{Berstein65}, who proved that a cogroup (in fact comonoid) in the category of associative algebras in free. J.-M.~Oudom remarked in \cite{Oudom99} that coassociativity of the cooperation is not even necessary to prove the  freeness. The difference with our case is in the compatibility relation. See also \cite{BH} for similar results in this direction.

\subsubsection{Quotient triples} Of course the quotient triple of $(As, Mag, MagFine)$ is $(As,As, Vect)$.  It would be interesting to find a small presentation 
of the intermediate quotient by the pre-Lie relator\index{pre-Lie relator}
$$\langle x;y,z\rangle:=(x\cdot y)\cdot z - x\cdot (y\cdot z) - (x\cdot z)\cdot y + x\cdot (z\cdot y)\ ,$$
which gives the good triple
$$\big(As, PreLie, \Prim_{As}PreLie)\ .$$

\subsection{The triple $(As, 2as, MB)$}\label{2as} (cf.~\cite{LRstr}). By definition a \emph{2-associative algebra}\index{2-associative algebra} or $2as$-algebra for short, is a vector space $A$ equipped with two associative operations denoted $a\cdot b$ and $a*b$. In the unital case we assume that $1$ is a unit for both operations. By definition a \emph{$As^c\--2as$-bialgebra} is a $2as$-algebra equipped with a coassociative cooperation $\dd$, whose compatibility relations are as follows:

a) $\cdot$ and $\dd$ satisfy the n.u.i.~compatibility relation (cf.~\ref{AsAsVect}),

b) $*$ and $\dd$ satisfy the Hopf compatibility relation (cf.~\ref{ex:comas}).

The free $2as$-algebra can be explicitly described in terms of planar trees and one can show that it has a natural $As^c\--2as$-bialgebra structure. Hypotheses {\texttt{(H0)}}, \Hone and  \Htwoepi are also easy to check (cf.~loc.cit.) and so there is a good triple of operads 
$(As, 2as, \Prim_{As}2as)$.
 It has been first proved in \cite{LRstr} where the primitive operad $\Prim_{As}2as$ has been shown to be the operad of multibrace algebras (abbreviated into $MB$-algebras, denoted $\Bi$-algebras in \cite{LRstr}). This is a very important operad since there is an equivalence between the category of $MB$-algebras and the category of cofree Hopf algebras $\HH$ equipped with  an isomorphism $\HH\cong T^c(\Prim \HH)$. 

Let us give some details on this equivalence. By definition a \emph{$MB $-algebra}\index{$MB $-algebra} (cf.~\cite{GV, LRstr}) is a vector space $A$ equipped with $(p+q)$-ary operations $M_{pq}$ satisfying some relations. Let $T^c(A)$ be the cofree counital coalgebra on $A$. Let 
$$*:T^c(A)\t T^c(A) \to T^c(A) $$
be the unique coalgebra morphism which extends the operations $M_{pq}:A^{\t p}\t A^{\t q} \to A$. Then the relations satisfied by the operations $M_{pq}$ imply that $(T^c(A), *, deconcatenation)$ is a cofree Hopf algebra. In the other direction, any cofree Hopf algebra determines a $MB $-structure on its primitive part. The details are to be found in \cite{LRstr}.

Hence we deduce that
$$(As, 2as, MB )$$
 is a good triple of operads (first proved in loc.cit.). One of the outcomes of this result was to give an explicit description of the free $MB $-algebra. Indeed the operad $2as$ can be explicitly described in terms of planar rooted trees. Thanks to this description and the structure theorem, one can describe the operad $MB $ in terms of trees (cf.~loc.cit.). Observe that moding out by the primitives gives the triple $(As, As, Vect)$. 

Since the functor $As\alg \to 2as\alg$ admits an obvious splitting (forgetful map), we can use it to construct the splitting of $\varphi$. Hence the idempotent $e$ is the same as in the case of the triple $(As,As,Vect)$. It was shown in \cite{LRstr} that the universal idempotent is given by the geometric series:
$$e = \Id - \Id \star \Id + \cdots +(-1)^{n-1} \Id^{\star n} + \cdots \ .$$
Here $\star$ stands for the convolution product.

\subsection{The triple $(As, 2as, Mag^{\infty})$}\label{as2asmag} Consider $As^c\--2as$-bialgebras with compatibility relations being both the n.u.i.~compatibility relation. So this type of bialgebras is different from the one described in \ref{2as}. It is immediate to check {\texttt{(H0)}}, \Hone and \Htwoepi are fulfilled. So there is a good triple of operads $(As, 2as, \Prim_{As}2as)$. With some more work (cf.~\cite{Ronco06}) one can show that the primitive operad is the operad $Mag^{\infty}$ which has one generating operation in each arity $n\geq 2$ and no relation. So there is a good triple of operads
$$(As, 2as, Mag^{\infty}).$$

\section{Dendriform, dipterous and Zinbiel bialgebras}\label{s:dendriform} In this section we study some triples $(As, \AA, \PP)$ which are over $(As, Zinb, Vect)$ with compatibility relation the semi-Hopf relation. Here $Zinb$ is the operad of Zinbiel algebras. The triple $(As, Zinb, Vect)$ is used by E.\ Burgunder in her analysis of a noncommutative version of the Kontsevich graph complex, cf. \cite{Burgunder08}. The triple 
$(As, Dipt, MB )$ is particularly interesting since it gives a structure theorem for cofree Hopf algebras, cf. \cite{LRstr}. The triple $(As, Dend, Brace)$ was studied by M.\ Ronco in \cite{Ronco00, Ronco02}. It gives rise to a noncommutative version of the Connes-Kreimer Hopf algebra, cf.\ \cite{LRcha}. The triple $(As, TriDend, \Prim_{As}TriDend)$ admits a quotient $(As, CTD, Com)$ which is strongly related to the quasi-shuffle algebras, cf. \ \cite{JLLctd}.

\subsection{Zinbiel algebra and semi-Hopf compatibility relation}\label{Zinbiel}  By definition a \emph{Zinbiel algebra}\index{Zinbiel algebra} is a vector space $A$ equipped with an operation denoted $a\g b$ satisfying the Zinbiel relation
$$(a\g b) \g c = a\g (b\g c + c \g b).$$
We note immediately that the operation $a*b:= a\g b + b \g a$ is associative (and commutative of course). The terminology comes from the fact that the Koszul dual operad is the operad of Leibniz algebras (cf.~\cite{JLLdig}). 

By definition a \emph{$As^c\--Zinb$-bialgebra}, or \emph{Zinbiel bialgebra},  is a Zinbiel algebra equipped with a coassociative cooperation $\dd$, whose compatibility relation is nonunital \emph{semi-Hopf}\index{semi-Hopf}, denoted $\between_{semiHopf}^l$ :

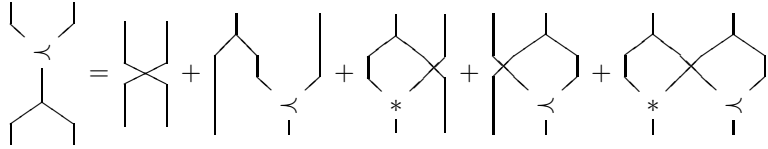
$$\cpbDeuxDeux{\g}{*{}}= \cpbB + \cpbCC{*{}}{\g} + \cpbDD{*{}}{{*}}+\cpbFF{*{}}{\g}+ \cpbGG{*{}}
{{*}}{*{}}{\g}$$

Observe that, as a consequence, the compatibility relation for the pair $(\dd, * )$ is Hopf (nonunital setting). It was obtained as a consequence of the semi-Hopf relation in the unital framework, given by
$$\DD(x\g y) = \DD(x)\g \DD(y),$$
where the tensor product of the bialgebra with itself has been equipped with following Zinbiel structure:
$$(a\t b) \g (a'\t b') = a*a' \t b\g b',$$
whenever it is defined, and
$$(a\t 1) \g (a'\t 1) = a\g a' \t 1,$$
otherwise (cf.~M.~Ronco \cite{Ronco00} and \ref{split}). The behavior of $\g$ with respect to the unit is given by
$$1\g x = 0, \qquad x\g 1 = x \ . $$
The free Zinbiel algebra over $V$ is the reduced tensor module $\To(V)$ (cf.\ \cite{JLLcup}) and the relationship between the tensors and the Zinbiel algebra structure is given by
$$v_{1}\cdots v_{n} = v_{1}\g (v_{2}\g (\cdots(v_{n-1} \g v_{n}) \cdots )).$$
Explicitly, the Zinbiel product is given by the \emph{half-shuffle}\index{half-shuffle}:
$$v_{1}\ldots v_{p}\g v_{p+1}\cdots v_{n}= v_{1}\sum_{\ss\in SH(p-1,n-p)}\ss (v_{2}\cdots v_{n})$$
where $SH(p-1,n-p)$ is the set of $(p-1,n-p)$-shuffles. As a consequence $(\To(V), *)$ is the (nonunital) shuffle algebra.

It can be shown that the free Zinbiel algebra satisfies both hypotheses \Hone and \Htwoiso. It is a consequence of the work of M.~Ronco \cite{Ronco00, Ronco02}, see E.~Burgunder \cite{Burgunder08} for a self-contained proof. Hence 
$$(As, Zinb, Vect)$$
 is a good triple of operads.
This example is interesting because it shows that, for a certain algebraic structure $\AA$, ($\AA = As$ here), there can be several different coalgebraic structures $\CC$ for which $(\CC, \AA, Vect)$ is a good triple. Here $\CC= As$ or $Zinb$.

We can revert the roles of $As$ and $Zinb$ in this example (cf.~\ref{coalgebraic}) and so there is a notion of $Zinb^c\--As$-bialgebra.  As a consequence 
$$(Zinb, As, Vect)$$
 is also a good triple (cf.~\ref{coalgebraic}). This result plays a role in the analysis of the Leibniz homology of the Lie algebra of derivations of an operad (cf.~\cite{Burgunder08}).

\subsection{Dipterous bialgebras}\label{ss:dipterous}  By definition a \emph{dipterous algebra}\index{dipterous algebra} is a vector space $A$ equipped with two binary operations denoted $a* b$ and $a\g b$ verifying the relations:
 \begin{eqnarray*}
 (x\g y)\g z  &=& x\g (y* z)\ , \\
(x* y)* z  &=& x* (y* z) \ .
\end{eqnarray*}
 By definition a \emph{$As^c\--Dipt$-bialgebra}, or \emph{dipterous bialgebra},  is a dipterous algebra equipped with a coassociative cooperation $\dd$, whose compatibility relation is Hopf with $*$ and semi-Hopf with $\g$. In fact one can put a unit on a dipterous algebra by requiring that $1$ is a unit for $*$ and that 
 $$1\g a = 0, \quad a\g 1=a,$$
 ($1\g 1$ is not defined). This is a particular case of a multiplicative  operad, see \ref{Hopfoperad} and \cite{JLLsci}.
 
The free dipterous algebra can be described in terms of planar trees. It can be shown that the free dipterous algebra satisfies both hypotheses \Hone and \Htwoiso (cf.~\cite{LRnote}). The primitive operad was proven to be the operad $MB $ introduced in loc.cit.~(cf.~\ref{2as}). Hence 
$$(As, Dipt, MB )$$
 is a good triple.  

\subsection{Dendriform bialgebras}\label{ss:dendriform} Because of the importance of the dendriform and brace notions, we give more details on this case. The results are due to M.~Ronco \cite{Ronco00, Ronco02}.

\subsubsection{Dendriform algebras} A \emph{dendriform algebra}\index{dendriform algebra}  $A$  is determined by two binary operations $ A\t A \to A$ called \emph {left} $(a,b)\mapsto a\g b$ and \emph{right} $(a,b)\mapsto a\d b$, satisfying the following three relations
 \begin{eqnarray*}
 (x\g y)\g z  &=& x\g (y* z)\ , \\
(x\d y)\g z  &=& x\d (y\g z) \ ,\\
(x* y)\d z &=& x\d (y\d z) \ ,
\end{eqnarray*}
where $x*y:= x\g y + x \d y$. From these axioms it follows that the operation $*$ is associative, hence a dendriform algebra is an example of dipterous algebra. The operad $Dend$ is obviously nonsymmetric. It has been shown in \cite{JLLdig} that  $Dend_{n}= \KK[PBT_{n+1}]$ where $PBT_{n+1}$ is the set of planar binary rooted trees with $n+1$ leaves (cf.\ \ref{pbtree} for more details on trees).

A \emph{unital dendriform algebra} is a vector space $A_{+}:=\KK\ 1 \oplus A$ where $A$ is a dendriform algebra and where the left and right products have been extended by:
$$1\g a = 0,\quad  a \g 1 = a ,\quad  1\d a = a,\quad  a \d 1 = 0 .$$
Observe that $1$ is then a unit for the associative operation $*$.  If $A_{+}$ and  $B_{+}$ are two dendriform algebras, then we can put a unital dendriform structure on their tensor product $A_{+}\t B_{+}$ by
$$(a\t b)\g (a'\t b') := (a*a')\t (b\g b'), \textrm { whenever } a\neq 1 \textrm{ and } a'\neq 1,$$
and a similar formula with $\d$ in place of $\g$. 

\subsubsection{Dendriform bialgebras} By definition a \emph{$As^c\--Dend$-bialgebra}, or \emph{dendriform bialgebra},  is a dendriform algebra $A$ equipped with a coassociative cooperation $\dd$, whose compatibility relations are as follows.

For the pair $(\dd,\d)$  it is given by $\between_{semiHopf}^r$:

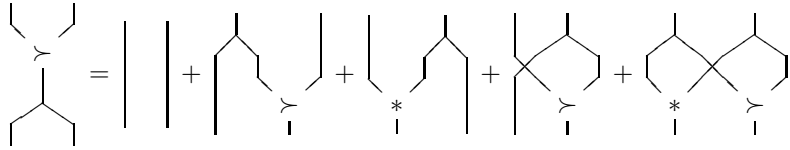
$$\cpbDeuxDeux{\d}{*{}}= \cpbA + \cpbCC{*{}}{\d} + \cpbEE{*{}}{{*}}+\cpbFF{*{}}{\d}+ \cpbGG{*{}}
{{*}}{*{}}{\d}$$

and  for the pair $(\dd,\g)$ it is given by $\between_{semiHopf}^l$ :

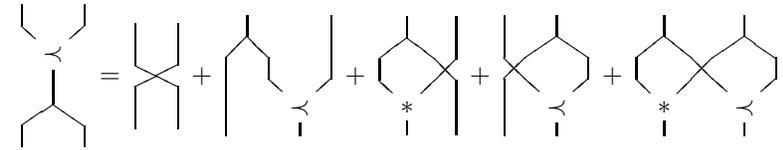
$$\cpbDeuxDeux{\g}{*{}}= \cpbB + \cpbCC{*{}}{\g} + \cpbDD{*{}}{{*}}+\cpbFF{*{}}{\g}+ \cpbGG{*{}}
{{*}}{*{}}{\g}$$

Equivalently, there is a coassociative diagonal $\DD: A_{+}\to A_{+}\t A_{+}$, such that $\DD(1)=1\t 1$ and 
$$\DD(a\g b) = \DD(a)\g \DD(b), \quad \DD(a\d b) = \DD(a)\d \DD(b).$$
The relationship between $\DD$ and $\dd$ is given by $\DD(a) = a\t 1 + 1 \t a + \dd(a)$. It follows that hypothesis {\texttt{(H0)}} is true.

Though $As$ and $Dend$ are nonsymmetric, the prop of dendriform bialgebras is not nonsymmetric because of the form of the compatibility relation.
Observe that the sum of these two relations gives the Hopf compatibility relation for the pair $(\dd, * )$. 

Let $Dend(V)_{+}$ be the free unital dendriform algebra. It has a functorial bialgebra structure constructed as follows. There is a unique map $\DD: Dend(V)_{+}\to  Dend(V)_{+}\t  Dend(V)_{+}$ which extends the map $V 
\to  Dend(V)_{+}\t  Dend(V)_{+}, v\mapsto v\t 1 + 1 \t v$ because $Dend(V)$ is free. The map $\DD$ is easily shown to be cocommutative. Hence the hypothesis \Hone is true.

\subsubsection{Structure theorem for dendriform bialgebras} The map of \sms $Dend\to As^c$ is obviously surjective. Let us choose the splitting which maps the cooperation $\mu_{n}$ to the operation $(\g)^n$ defined inductively as
$$(\g)^1(x_{0}, x_{1}) = x_{0}\g x_{1}, \quad (\g)^n(x_{0},\ldots ,  x_{n}) = x_{0}\g  (\g)^{n-1}(x_{1},\ldots ,  x_{n}) .$$

Under this choice the versal idempotent is:
$$e = \sum_{n\geq 1} (-1)^{n+1} (\g)^n\circ \dd^{n},$$
where $\dd^{n}$ is the iteration of $\dd$. This is exactly the idempotent used by M.~Ronco in her proof of the structure theorem for dendriform bialgebras.

Since the hypotheses {\texttt{(H0)}}, \Hone and \Htwoepi are true, the structure theorem holds for dendriform bialgebras. The good triple $(As, Dend, \Prim_{As}Dend)$ will be completely understood once we make the primitive operad explicit.

 \subsubsection{Brace algebras}\label{brace} The \emph{Brace} operad admits one $n$-ary operation $\{-;-,\cdots, -\}$, for all $n\geq 2$, as generators and the relations are:
 $$Br_{n,m}:  \{\{x;y_1,\ldots , y_n\};z_1,\ldots ,z_m \} = \sum\{x; \ldots ,\{y_1; \ldots\}, \ldots , \{y_n; \ldots , \}, \ldots \}.$$
 On the right-hand side the dots are filled with the variables $z_i$'s (in order) with the convention $\{y_k; \emptyset \}= y_k$. The sum is over all the possibilities of inserting the variables $z_{i}$'s. Equivalently a brace algebra is a $MB $-algebra for which $M_{pq}= 0$ for $p\geq 2$ and $M_{1q}=\{-;-,\ldots , -\}$.  The first nontrivial relation, which relates the 2-ary operation and the 3-ary operation  reads
$$Br_{1,1}:\qquad \{\{x;y\} ;z\} - \{x; \{y; z\}\}  = \{x;y,z\} + \{x; z, y\}\ .$$
As a consequence we deduce that the associator of the 2-ary operation is right-symmetric:
$$\{\{x;y\} ;z\} - \{x; \{y; z\}\}  = \{\{x;z\} ;y\} - \{x; \{z; y\}\}  \ .$$
So the binary brace operation is in fact a pre-Lie operation.
 
 There is a functor $Dend\alg \to Brace\alg$ given by the following formulas:
 $$\{x;y_{1},\ldots ,y_{n}\}:= \sum_{i=0}^{n} (-1)^{i} \omega_{\g}(y_{1},\ldots ,y_{i})\d x \g  \omega_{\d}(y_{i+1},\ldots ,y_{n}) ,$$
 where $ \omega_{\g}(y_{1})= y_{1},  \omega_{\g}(y_{1},\ldots ,y_{i})= y_{1}\g \omega_{\g}(y_{2},\ldots ,y_{i})$ and $ \omega_{\d}(y_{1})= y_{1},  \omega_{\d}(y_{1},\ldots ,y_{i})= \omega_{\g}(y_{1},\ldots ,y_{i-1})\d y_{i}$.
In low dimensions it reads
 $$\{x;y\} := x\g y - y \d x, $$
 and
 $$\{x;y,z\}= x\g (y\d z) - y\d x \g z + (y\g z)\d x\ . $$
  
  One can  verify that a brace product of primitive elements in a dendriform bialgebra is still primitive. Hence there is a morphism  $Brace(V) \to  \Prim_{As}Dend(V)$, which can be extended as a map 
  $$\Theta (V) : T^c(Brace(V)) \to T^c(  \Prim_{As}Dend(V)) \cong Dend(V).$$
  
  On the other hand,  $T^c(Brace(V))$ can be equipped with an associative product and a left product (compare with the dipterous case). But, because  $M_{pq}= 0$ for $p\geq 2$ the extra relation $(x\d y)\g z  = x\d (y\g z) $ holds. Therefore $T^c(Brace(V))$ is a dendriform algebra. So there is a functorial map 
  $$\Xi(V) : Dend(V) \to T^c(Brace(V)).$$
  One can show that $\Theta$ and $\Xi$ are inverse to each other. 
As a consequence of this discussion   $(As, Dend, Brace)$ is a good triple of operads. 

Comments. The Hopf structure of the free dendriform algebra was first constructed in \cite{LRHopf}. See \cite{Aguiar06} for an alternative basis with nice behavior with respect to the coproduct. The primitive operad $\Prim_{As}Dend$ was first  identified to be the brace operad by Mar\'\i a Ronco in  \cite{Ronco00}. The structure theorem  was proved in \cite{Ronco02} and in \cite{Chapoton02}. It was the first example outside the classical framework and the one which motivated this theory.

Since $Dend$ is a quotient of $Dipt$ by the relation $(x\d y)\g z  = x\d (y\g z) $ and since the operation $(x\d y)\g z  - x\d (y\g z) $ in $Dipt$ is primitive, it follows from Proposition \ref{quotient} that $(As, Dend, MB /\sim)$ is a good triple. The quotient  $MB /\sim$  turns out to be the operad $Brace$.

If we mod out by the primitive operation $x\g y - y \d x$, then we get the good triple $(As, Zinb, Vect)$.

\subsection{Tridendriform algebra}\label{tridendriform} The notion of dendriform algebra admits several generalizations. One of them is  the notion of \emph{tridendriform algebra}\index{tridendriform algebra} (originally called dendriform trialgebra in \cite{LRtri}). It has three generating operations denoted $\g$ (left), $\d$ (right), and $\cdot$ (dot or middle). They  satisfy the following 11 relations (one for each cell of the pentagon, see loc.cit.):

\begin{eqnarray*}
(x \g y) \g z &=& x \g (y * z)\ , \\
(x \d y) \g z &=& x \d (y \g z)\ , \\
(x * y) \d z &=& x \d (y \d z)\ , \\
(x \d y) \cdot z &=& x \d (y \cdot z)\ , \\
 (x \g y) \cdot z &=& x \cdot (y \d z)\ , \\
  (x \cdot y) \g z &=& x \cdot (y \g z)\ , \\
(x \cdot y) \cdot z &=& x \cdot (y \cdot z)\ , \\
\end{eqnarray*}
where $x*y := x\g y + x\d y  + x\cdot y$. 
  
The operad $Tridend$ is obviously binary, quadratic and nonsymmetric. The free tridendriform algebra on one generator has been shown to be linearly generated by the set of all planar rooted trees in \cite{LRtri}.

Using the existence of a partial unit one can put a structure of $As^c\--Tridend$-bialgebra structure on $Tridend(V)$ as in \cite{JLLsci}. The coefficients $\aa$ and $\bb$ (cf.~\ref{split}) are given by:
$$x\g 1=x=1\d x,\ \textrm{and}\ 1\g x= x\d 1 = 1\cdot x = x\cdot 1 = 0\ .$$
These choices are coherent with the operad structure of $Dend$ and therefore, by \cite{JLLsci} (see also \ref{split}), there is a well-defined notion $As^c\--Dend$-bialgebra for which the hypotheses {\texttt{(H0)}} and \Hone are fulfilled. 

Hypothesis \Htwoepi is easy to check, and therefore the triple $$(As, Tridend, \Prim_{As}Tridend)$$
 is good. The operad $\Prim_{As}Tridend$ can be described explicitly as a mixture of the brace structure and the associative structure, cf.~\cite{PR}.
 
One of the interesting points about the good triple $(As, Tridend, Brace+As)$ is its quotient $(As, CTD, Com)$, where $CTD$ stands for the Commutative TriDendriform algebra operad. The commutativity property is
$$x\g y = y\d x,\ \textrm{and}\ x\cdot y = y \cdot x \ .$$
Hence a \emph{$CTD$-algebra} can be described by two generating operations $x\g y$ and $x\cdot y$ (the second one being symmetric), satisfying the relations:
\begin{eqnarray*}
(x \g y) \g z &=& x \g (y * z)\ , \\
  (x \cdot y) \g z &=& x \cdot (y \g z)= (x\g z)\cdot y\ , \\
(x \cdot y) \cdot z &=& x \cdot (y \cdot z)\ . \\
\end{eqnarray*}
The good triple 
$$(As,CTD,Com)$$
 has been studied in \cite{JLLctd}\footnote{The relation $x \cdot (y \g z)= (x\g z)\cdot y$ is missing in the axioms of a CTD-algebra given in \cite{JLLctd}.} and shown to be related with the quasi-shuffle algebras in the following sense. The left adjoint $U:Com \to CTD\alg$ gives a classical bialgebra $U(R)$ which is the quasi-shuffle algebra over the commutative algebra $R$.

\section {$Lie^c\--Lie$-bialgebras}\label{LieLie}

In this section we work over a characteristic zero field. 
We introduce the notion of $Lie^c\--Lie$-bialgebra, different from the classical notion of Lie bialgebra, and we prove a rigidity theorem for $Lie^c\--Lie$-bialgebras.

\subsection{Definition}\label{compatibilityLieLie} A \emph{$Lie^c\--Lie$-bialgebra}\index{Lie@$Lie^c\--Lie$-bialgebra}  is a vector space $A$ which is a Lie algebra for the bracket $[x,y]$, a Lie coalgebra for $\dd_{[\, ]}$ and whose compatibility relation is $\between_{Lily}$\index{Lily}:

$$\cpbdeuxdeux = 2\Big(\cpbA - \cpbB\Big) + \frac{1}{2}\Big(\cpbC + \cpbD + \cpbE + \cpbF\Big)$$

Here  $\vcenter{\xymatrix@R=1pt@C=1pt{
\ar@{-}[ddrr]&&&&\ar@{-}[ddll]\\
&&&&\\
&&*{}\ar@{-}[dd]&&\\
&&&&\\
&&&&
}}$
  stands for the bracket $[-, -]$ and 
 $\vcenter{\xymatrix@R=1pt@C=1pt {
&&&&\\
&&&&\\
&&*{}\ar@{-}[uu]&&\\
&&&&\\
\ar@{-}[uurr]&&&&\ar@{-}[uull]
}}$ 
stands for the cobracket $\dd_{[\, ]}$.

Observe that the notion of $Lie^c\--Lie$ bialgebra is completely different from the notion of Lie bialgebras, since, in this latter case, the compatibility relation is the cocycle condition $\between_{biLie}$\index{biLie compatibility relation}:

$$\cpbdeuxdeux = \cpbC + \cpbD + \cpbE + \cpbF$$
 
 In particular for $\between_{Lily}$ there is a $\Phi_{1}$-term, so there is some chance for the existence of a rigidity theorem.

In order to show that the free Lie algebra $Lie(V)$ is equipped with   a structure of $Lie^c\--Lie$-bialgebra, we are going to use the tensor algebra $\To(V)$ for $V= \KK x_1 \oplus \cdots \oplus  \KK x_n$. Hence $\To(V)$ is the space of noncommutative polynomials without constant term in the variables $x_{i}$'s. The coproduct $\dd$ on $\To(V)$ is the deconcatenation coproduct (cf.~\ref{AsAsVect}). Recall that $Lie(V)$ is made of the Lie polynomials, that is the polynomials generated by the $x_{i}$'s under the bracket operation. The degree of a homogeneous polynomial $X$ is denoted $\vert X\vert$. We use the involution $X \mapsto  \overline X$ on $\To(V)$ which is the identity on the $x_{i}$'s and satisfies $\overline {(XY)} = \overline Y \  \overline X$.

\begin{lemma}\label{bracketsymmetry} If $X\in Lie(V)$, then $\overline X = -(-1)^{\vert X\vert} X$  and 
$$2\dd (X) = X_{1} \t X_{2} - (-1)^{p+q}\overline X_{2} \t \overline X_{1},$$
where $\dd (X) =: X_{1} \t X_{2}$ (Sweedler notation)  and $p= \vert X_{1}\vert, q= \vert X_{2}\vert$.
\end{lemma}
\begin{proo} The proof is by induction on the degree $n$ of $X$. We assume that these formulas are true for $X$, and then we prove them for $[X,z]$ where $z$ is of degree $1$.

For the first formula we get 
$$\overline {[X,z]} = \overline {(Xz-zX)} = z\overline X - \overline X z =  -(-1)^{n}([z,X])= 
 -(-1)^{n+1} [X,z] ,$$
as expected.

For the second formula, the n.u.i.~compatibility relation and the induction hypothesis give:
\begin{eqnarray*}
\dd([X,z]) &=& \dd(Xz-zX)\\
& =& X\t z + X_{1}\t X_{2}z  - z\t X - zX_{1}\t X_{2}\\
&=& X\t z - z\t X+\frac{1}{2} \Big(  X_{1}\t X_{2}z - zX_{1}\t X_{2}\\
&& - (-1)^{p+q} \overline X_{2}\t \overline X_{1}z +  (-1)^{p+q} z\overline X_{2}\t \overline X_{1}\Big).
\end{eqnarray*}

On the other hand we have

\medskip

$\begin{array}{l}
[X,z]_{1}\t [X,z]_{2} - (-1)^{p+q}\overline{[X,z]_{1}}\t \overline{[X,z]_{2}}\\
= \frac{1}{2}\big( X\t z - z\t X + X_{1}\t X_{2}z - zX_{1}\t X_{2} -(-1)^{n+1}z\t \overline X \\
\qquad + (-1)^{n+1}\overline X \t z -  (-1)^{p+q+1}z\overline X_{2}\t \overline X_{1}  + (-1)^{p+1+q}\overline X_{2}\t \overline X_{1}z\big).
\end{array}$

\medskip

The two expresssions are equal, because, since $X$ is a Lie polynomial, we have $X=-(-1)^n\overline X$.
\end{proo}

\begin{prop}\label{Lieinternal} Let $\dd_{[,]}:= \dd - \tau \dd$. The image of $Lie(V)$ by $\dd_{[,]}$ is in $Lie(V)\t Lie(V)$.
\end{prop}
\begin{proo}  The proof is by induction on the degree $n$ of $X\in Lie(V)$. It is immediate for $n=1$. Suppose that $X\in Lie(V), \dd (X) = X_{1} \t X_{2}$ and $ X_{1} , X_{2}\in Lie(V)$. We observe that, by Lemma \ref{bracketsymmetry} we have 
$$ \dd (X) =: X_{1}\t X_{2}= \frac{1}{2}\big(X_{1} \t X_{2} -   X_{2} \t X_{1}\big).$$
We are going to show that, for any element $z$ of degree $1$, we have $\dd_{[,]}([X,z]) \in Lie(V)\t Lie(V)$. We compute:
\begin{eqnarray*}
\dd_{[,]}([X,z])&=& (\dd - \tau \dd)(Xz - zX)\\
&=& X\t z - z\t X+ X_{1}\t X_{2}z  - zX_{1}\t X_{2} \\
&&  - z\t X + X\t z - X_{2}z\t X_{1}  + X_{2}\t zX_{1}\\
&=& 2\big( X\t z - z\t X\big) \\
&&+ \frac{1}{2}\big(  X_{1}\t X_{2}z  - zX_{1}\t X_{2}  - X_{2}z\t X_{1}  + X_{2}\t zX_{1}\\
&& -X_{2}\t X_{1}z  + zX_{2}\t X_{1}  + X_{1}z\t X_{2}  - X_{1}\t zX_{2}\big)\\
&=&  2\big( X\t z - z\t X\big)\\
&& + \frac{1}{2}\big(  X_{1}\t [X_{2},z]  - [z,X_{1}]\t X_{2}  - [X_{2},z]\t X_{1}  + X_{2}\t [z,X_{1}]\big).
\end{eqnarray*}
So we have proved that $\dd_{[,]}([X,z]) \in Lie(V)\t Lie(V)$.
\end{proo}

\begin{prop}\label{LieV}  On $Lie(V)$ the bracket operation $[x,y]$ and the bracket cooperation $\dd_{[,]}$ satisfy the compatibility relation  $\between_{Lily}$  (cf.~\ref{compatibilityLieLie}).
\end{prop}
\begin{proo}   
Let $X,Y\in Lie(V)$. We compute $\dd_{[,]}([X,Y])$ :
\begin{eqnarray*}
\dd_{[,]}([X,Y])&=&  (\dd - \tau \dd)(XY - YX)\\
&=&  \dd(XY) - \dd (YX) - \tau\dd (XY)  + \tau \dd(YX)\\
&=&+X\t Y +X_{1}Y\t X_{2}+ XY_{1}\t Y_{2} \\
&&-Y\t X -Y_{1}X\t Y_{2}-YX_{1}\t X_{2} \\
&&-Y\t X -X_{2}\t X_{1}Y- Y_{2}\t XY_{1} \\
&& +X\t Y +Y_{2}\t Y_{1}X+ X_{2}\t YX_{1} \\
&=& 2\big( X\t Y - Y \t X \big)\\
&&  +[X_{1},Y]\t X_{2}+ [X,Y_{1}]\t Y_{2} + X_{1}\t [X_{2},Y]+ Y_{1}\t [X,Y_{2}] \\
&=& 2\big( X\t Y - Y \t X \big) +\frac{1}{2}\Big( [X_{[1]},Y]\t X_{[2]}\\
&& + [X,Y_{[1]}]\t Y_{[2]} + X_{[1]}\t [X_{[2]},Y]+ Y_{[1]}\t [X,Y_{[2]}]\Big) . \\
\end{eqnarray*}
Observe that, in this computation, we have used the fact that, for any element $Z\in Lie(V)$, the element $\dd(Z)=Z_{1}\t Z_{2}$ is antisymmetric, that is $Z_{1}\t Z_{2}=-Z_{2}\t Z_{1}$ (cf.~\ref{bracketsymmetry}). As a consequence we have $\frac{1}{2}\dd_{[,]}(Z)=\dd(Z)$.
\end{proo}

\begin{thm} In characteristic zero, the prop $Lie^c\--Lie$ satisfies the hypotheses {\texttt{(H0)}} \Hone \Htwoiso, therefore the triple $(Lie, Lie, Vect)$ (with $\between = \between_{Lily}$) is a good triple. Hence any connected $Lie^c\--Lie$-bialgebra is both free and cofree.
\end{thm}
\begin{proo} It is clear that the compatibility relation $ \between_{Lily}$ is distributive. By Propositions \ref{Lieinternal} and \ref{LieV} the hypothesis \Hone is fulfilled. Let us prove \Htwoiso.

We have seen in the last proof that $\dd_{[,]}(Z)=2\dd(Z)$ when $Z\in Lie(V)$. The cooperation $\dd$ induces the isomorphism map $\varphi_{As}(V): As(V) \to As^c(V)$ which identifies the generator of $As_{n}$ with its dual. Hence, restricted to $Lie(V)$ it is injective. So the map $\varphi_{Lie}(V):Lie(V)\to Lie^c(V)$, induced by $\dd_{[,]}$  is injective. Since $Lie(n)$ and $Lie^c(n)$ have the same dimension, it is an isomorphism.
\end{proo}

\subsection{The conjectural triple $(Lie, PostLie, \Prim_{As}PostLie)$} By definition, cf.~\cite{Vallette06},  a \emph{PostLie algebra}\index{PostLie algebra} is a vector space $A$ equipped with two operations 
$x\circ y$ and $[x, y]$ which satisfy the relations
$$\displaylines{[x, y]=-[y, x],\cr
[x,[y,z]] + [y,[z,x]]+ [z,[x,y]] =0,\cr
(x\circ y)\circ z - x\circ (y\circ z) - (x\circ z) \circ y + x\circ (z\circ y) = x \circ [y,z],\cr
[x,y]\circ z = [x\circ z, y] + [x,y\circ z]\ .}
$$
It is shown in loc.cit.~that the operad $PostLie$ is the Koszul dual of the operad $ComTrias$. A PostLie algebra is a Lie algebra for the bracket $[x,y]$. But it is also a Lie algebra for the operation $\{x,y\}:= x\circ y - y\circ x + [x,y]$ (cf.~loc.cit). We conjecture that there exists a notion of $Lie^c\--PostLie$ bialgebra such that the free PostLie algebra is such a bialgebra. Hopefully  there is a good triple of operads $(Lie, PostLie, \Prim_{As}PostLie)$. The isomorphism of   $PostLie=Lie\circ PBT$, where $PBT=\oplus_{n}\KK[PBT_{n}]$ proved in \cite{Vallette05} is an evidence in favor of this conjecture. It is not clear what is the algebraic structure one should put on $PBT$ to make it work ($Mag$ is one option out of many).

\section {$NAP^c\--\AA$-bialgebras}\label{NAPA}

Triples of the form $(NAP^c, \AA, \Prim_{NAP}\AA)$ come from the work of Muriel Livernet \cite{Livernet06}.

\subsection{Pre-Lie algebras}\label{preLie} By definition a \emph{pre-Lie algebra}\index{pre-Lie algebra}  is a vector space $A$ equipped with a binary operation  $a\cdot b$ which satisfies the following relation
$$(x\cdot y)\cdot z - x\cdot (y\cdot z) = (x\cdot z)\cdot y - x\cdot (z\cdot y)$$
(right-symmetry of the associator). The free pre-Lie algebra has been described in terms of abstract trees in \cite{CL}.

\subsection{$NAP$-algebra}\label{NAP}\index{NAP} By definition a \emph{non-associative permutative algebra}\index{non-associative permutative algebra}, or $NAP$ algebra for short, is a vector space $A$ equipped with a binary operation denoted $a b$ which satisfies the following relation
$$(xy)z = (xz)y\ .$$
In fact we are going to use the notion of $NAP$-coalgebra, whose relation is pictorially as follows:
$$
{\vcenter{\xymatrix@R=2pt@C=2pt{
&&*{}\ar@{-}[d]&&                        &&&*{}\ar@{-}[d]&&   \\
&&*{}\ar@{-}[dl]\ar@{-}[dr]&&         &&&*{}\ar@{-}[dl]\ar@{-}[dr]&&\\
&*{}\ar@{-}[d]&&*{}\ar@{-}[ddd]&    =     &&*{}\ar@{-}[d]&&*{}\ar@{-}[d]&\\
&*{}\ar@{-}[dl]\ar@{-}[dr]&&&         &&*{}\ar@{-}[dl]\ar@{-}[ddrr]&&*{}\ar@{-}[ddll]&\\
*{}\ar@{-}[d]&&*{}\ar@{-}[d]&&                        &*{}\ar@{-}[d]&&&&\\
*{}&&*{}&*{}&  &*{}&*{} &*{}&*{} &*{} 
}}}
$$

\subsection{$NAP^c\--PreLie$-bialgebra}\label{NAPpreLie} By definition a \emph{$NAP^c\--PreLie$-bialgebra} is a pre-Lie algebra equipped with a NAP cooperation $\dd$, whose compatibility relation $\between_{Liv}$ is as follows:

$$\cpbdeuxdeux = \cpbA +\cpbC +\cpbD  .$$

It has been shown by M. Livernet in \cite{Livernet06}, Proposition 3.2, that the free pre-Lie algebra is naturally a $NAP^c\--PreLie$-bialgebra. She also proved the rigidity theorem for $NAP^c\--PreLie$-bialgebras by providing an explicit idempotent. This result follows also from our general result, since the coalgebra map $PreLie(V) \to NAP^c(V)$ is an isomorphism (cf.~loc.cit.). The explicit description of the universal idempotent in terms of generating operations and cooperations (as described in \ref{explidemp}) is to be found in loc.cit.

\subsection{$NAP^c\--Mag$-bialgebras} The compatibility relation for $NAP^c\--Mag$-bialgebras is $\between_{Liv}$. Hypothesis {\texttt{(H0)}} is clearly fulfilled. Here is a proof of Hypothesis \Hone.

\begin{prop}\label{NAPMag} On the free magmatic algebra $Mag(V)$ there is a well-defined cooperation $\dd$ which satisfies the $NAP^c$ relation, that is $(\dd\t \Id) \dd = (\Id \t \tau)(\dd \t \Id) \dd$, and the Livernet compatibility relation $\between_{Liv}$.
\end{prop}
\begin{proo} We use the inductive method described in \ref{Hopfoperad}. We let 
$$\dd:Mag(V) \to Mag(V)\t Mag(V)$$
be the unique linear map which sends $V$ to $0$ and which satisfies the compatibility relation $\between_{Liv}$. Here the tensor product $Mag(V)\t Mag(V)$ is equipped with its standard magmatic operation. In low dimension we get
\begin{eqnarray*}
\dd(x\cdot y) &=&x\t y ,\\
\dd((x\cdot y)\cdot z) &=&x\cdot y \t z + x \t y\cdot z + x\cdot z \t y , \\
\dd(x\cdot (y\cdot z)) &=& x \t y\cdot z . \\
\end{eqnarray*}
Remark that the pre-Lie relator is primitive. We show that $\dd$ satisfies the $NAP^c$-relation  (see the diagram above) inductively, by using the natural filtration of $Mag(V)=\oplus_{n\geq 1}\KK[PBT_{n}]\t V^{\t n}$. 
Applying $\between_{Liv}$ twice we get 

$$\cpbdeuxTrois = \cpbMa + \cpbMb + \cpbMc$$
and then
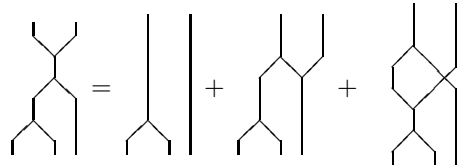
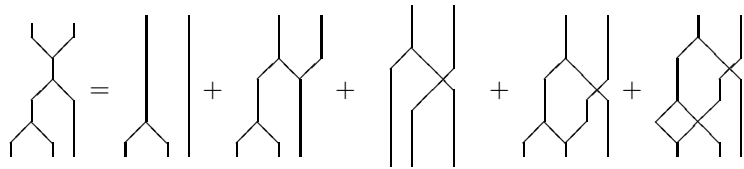
$$\cpbdeuxTrois = \cpbMa + \cpbMb + \cpbMd +\cpbMe +\cpbMf$$

Let us denote by $\oo_{1}+\oo_{2}+\oo_{3}+\oo_{4}+\oo_{5}$ the five terms on the right-hand side of the last line. We check that $(\Id \t \tau)\oo_{1}= \oo_{3}$, and that,  under the $NAP^c$ relation, we have $(\Id \t \tau)\oo_{2}= \oo_{4}$ and $
(\Id \t \tau)\oo_{5}= \oo_{5}$. Hence we have proved that the $NAP^c$-relation holds. By Theorem \ref{thm:primitive} there is a triple of operads $(NAP, Mag, \Prim_{NAP}Mag)$.
\end{proo}

\subsection{Conjecture on $NAP^c\--Mag$-bialgebras} We mentioned in the proof of Proposition \ref{NAPMag} that the pre-Lie relator is primitive in the bialgebra $Mag(V)$. Moding out by the ideal generated by this pre-Lie relator gives the $NAP^c\--PreLie$-bialgebra $PreLie(V)$. It follows that the map $\varphi(V)$ described in \ref{phi} is the composite
$$\varphi(V): Mag(V) \epi PreLie(V)\cong NAP^c(V)\ .$$

\noindent {\bf Conjecture}. The coalgebra map $\varphi(V): Mag(V) \epi  NAP^c(V)$ admits a coalgebra splitting.

It would follow that there is a good triple of operads 
$$(NAP, Mag, \Prim_{NAP}Mag)$$
with quotient triple $(NAP, PreLie, Vect)$. The operad $ \Prim_{NAP}Mag$ has no generating operation in arity 2, but has a generating operation in arity 3, namely the \emph{pre-Lie relator}\index{pre-Lie relator}
$$\langle xyz\rangle :=(x\cdot y)\cdot z - x\cdot(y\cdot z) - (x\cdot z)\cdot y + x\cdot(z\cdot y)\ .$$

\section{Some examples of the form $(\AA, \AA, Vect)$} We have already discussed the cases $(Com, Com, Vect)$ and $(As, As, Vect)$. We give here some more examples, some of them being already in the literature. The main point is to unravel the compatibility relation. The triple $(Dup,Dup,Vect)$ will be treated in Chapter 5.

\subsection{The triple $(Mag, Mag, Vect)$}\label{bimag} The free magmatic algebra $Mag(V)= \oplus_{n\geq 1}\KK[PBT_{n}]\t V^{\t n}$  inherits a comagmatic coalgebra structure under the identification of the basis $PBT_{n}$ of $Mag_{n}$ with its dual. An immediate inspection shows that the magmatic operation and the comagmatic cooperation are related by the magmatic compatibility relation $\between_{mag}$ :
$$\cpbdeuxdeux =   \cpbA $$
Hence Hypotheses {\texttt{(H0)}} and \Hone are fulfilled for $Mag^c\--Mag$-bialgebras. Since the map
$\varphi(V): Mag(V)\to Mag^c(V)$ is easily seen to be the identification of the basis of $Mag_{n}$ with its dual, it is an isomorphism and Hypothesis \Htwoiso is fulfilled. By Theorem \ref{thm:rigidity} 
$$(Mag, Mag, Vect)$$
is a good triple of operads. 

\noindent{\bf Exercise.} Describe the idempotent $e$ explicitly in terms of the generating operation and the generating cooperation. The answer is to be found in \cite{Burgunder08}.

\subsection{The triple $(2as, 2as, Vect)$}\label{2as2as} The operad $2as$ admits a basis made of planar trees. In fact, for $n\geq 2$, the space $2as_{n}$ is spanned by two copies of the set of planar trees with $n$ leaves. So it is immediate to describe the $2as$-coalgebra structure on the same space. We put the 
compatibility relations given by the following tableau:
\begin{displaymath}
\begin{array}{c||c|c}
 & *&\cdot \\
\hline
\hline
\dd_{\cdot} &Hopf & n.u.i. \\
\hline
\dd_{*} & n.u.i. & Hopf 
\end{array}
\end{displaymath}

It was shown in \cite{LRstr} (see also \ref{2as}) that $2as(V)$ equipped with $\dd_{\cdot}$ satisfies the first row of the compatibility relations. Inverting the role of $\cdot$ and $*$ it is clear that there is an associative cooperation $\dd_{*}$ which satisfies the second row of compatibility relations. Therefore the free $2as$-algebra is a $2as^c\--2as$-bialgebra. It is interesting to note that, in this case, the isomorphism between the free $2as$-algebra and the free $2as$-coalgebra is not given by identifying the basis with its dual.


\subsection{The triple $(\AA, \AA, Vect)$\label{AAVect} for a multiplicative operad $\AA$} Let $\AA$ be a binary operad with split associativity (cf. \ref{split} and \cite{JLLsci}). We denote by $x*y$ the associative operation and by $\aa(\circ), \bb(\circ)$ the coefficients such that 
$$x\circ 1 = \aa(\circ) x\quad \textrm{and}\quad 1\circ x=  \bb(\circ) x\ .$$
We always assume $\aa(*)=1_{\KK}= \bb(*)$. Let us suppose moreover that all the relations in $\AA$ are \emph{ generalized associativity}, that is
$$(x\circ_{1} y)\circ_{2}z = x\circ_{3} (y\circ_{4}z)\ ,$$
where $\circ_{i}\in \AA(2)$. 
For each generating operation $\circ$ we denote by $\DD_{\circ}$ the associated cooperation in the unital framework. By definition a $\AA^c\--\AA$-bialgebra, called bi$\AA$-bialgebra, has the following compatibility relations:
$$\DD_{\circ}(x\bullet y)= \DD_{*}(x)\bullet \DD_{\circ}(y)$$
where, by definition,
$$(x\t y)\bullet (x'\t y') = (x* x')\t (y\bullet y')$$
(Ronco's trick, see \ref{split} for the convention when $y=1=y'$). In order to formulate this compatibility relation in the nonunital framework we need to introduce the reduced cooperations $\dd_{\circ}$, defined by the equality
$$\DD_{\circ}= \aa(\circ) x\t 1 + \bb(\circ) 1\t x + \dd_{\circ}\ .$$
The compatibility relations become:
$$\cpbDeuxDeux{\bullet}{\circ}= \bb(\circ)\bb(\bullet) \cpbA +\   \aa(\circ)\aa(\bullet) \cpbB +
\bb(\circ) \cpbCC{{*}}{\bullet}+ \aa(\circ)\aa(\bullet) \cpbDD{{*}}{{*}}$$
$$ +\  \bb(\bullet)  \cpbEE{\circ}{{*}} +\ \cpbFF{\circ}{\bullet} +\  \cpbGG{{*}}{{*}}{\circ}{\bullet} $$
 
 \begin{prop} If $\AA$ is a multiplicative operad with generalized associativity relations, then the free $\AA$-algebra has a natural structure of bi$\AA$-bialgebra.
\end{prop}
\begin{proo} Let us work in the unital framework (cf.~\ref{split}). First, we construct  
 a cooperation $\DD_{\circ}$ on $\AA(V)_{+}$ for each generating binary operation $\circ$ by induction as follows. First, $\DD_{\circ}(v)= \aa(\circ) v\t 1 + \bb(\circ) 1\t v$. Second, we use the compatibility relations to define $\DD_{\circ}$ on $\AA(V)_{2}$, then on $\AA(V)_{3}$ and so forth. So the maps $\DD_{\circ}$ are uniquely defined and satisfy the compatibility relations. Let us use the inductive argument to prove the generalized associativity relations.
 
 One one hand we have:
\begin{eqnarray*}\DD_{\circ}((x\circ_1 y)\circ_2 z) &=& \DD_{*}(x\circ_1 y)\circ_2 \DD_{\circ}(z)\\
&=& (\DD_{*}(x)\circ_1 \DD_{*}(y))\circ_2 \DD_{\circ}(z)
\end{eqnarray*}
On the other hand we have:
\begin{eqnarray*}\DD_{\circ}(x\circ_3 ( y\circ_4 z)) &=& \DD_{*}(x)\circ_3  \DD_{\circ}(y\circ_4  z)\\
&=& \DD_{*}(x)\circ_3  (\DD_{*}(y)\circ_4\DD_{\circ}( z))
\end{eqnarray*}
Assuming that the generalized associativity relations hold in some dimension (including the associativity of $*$), we prove from this computation that they hold one step further.

Again by induction we can show, by a straightforward verification,  that these cooperations do satisfy the $\AA^c$ relations. So $\AA(V)_{+}$ is a unital-counital $\AA^c\--\AA$-bialgebra and, by restriction, $\AA(V)$ is a $\AA^c\--\AA$-bialgebra.
\end{proo}

\subsection{Examples}  The operads 
$$As, Dend, Dipt, 2as, Dup, Dup^!, Tridend, Dias, Trias,  Quad, Ennea$$
 (cf. \cite{JLLdig, LRnote, LRstr, LRtri, ALo, Leroux04, EG}) are examples of multiplicative operads  with generalized associative relations.  
 In some of these examples the map $\varphi$ is an isomorphism. However it is not always true: $As$ is a counter-example. In the following section we study in more details the case $\AA=Dend$. Observe that in the case of $2as$ (resp. $Dup$) we get a type of bialgebras which is different from the type studied in \ref{2as2as}  (resp. \ref{DupDup}) since the compatibility relations are different.

\subsection{The triple $(Dend, Dend, Vect)$} Let us make explicit the particular case $\AA=Dend$ which has been treated in details in \cite{Foissy}. Recall that the coefficients $\aa$ and $\bb$ are given by the relations
$$1\g x = 0 = x\d 1\quad \textrm{and} \quad x\g 1 = x = 1\d x .$$
The compatibility relations for the reduced cooperations read as follows.
\noindent For the pair $(\dd_{\d},\d )$:

$$\cpbDeuxDeux{\d}{\d} = \cpbA + 0 + \cpbCC{{*}}{\d} + 0 + \cpbEE{\d}{{*}} + \cpbFF{\d}{\d} + 
\cpbGG{{*}}{{*}}{\d}{\d}$$

\noindent for the pair $(\dd_{\g},\d)$:

$$\cpbDeuxDeux{\d}{\g} = 0 + 0 + 0 + 0 + \cpbEE{\g}{{*}} + \cpbFF{\g}{\d} + \cpbGG{{*}}{{*}}{\g}{\d}$$

\noindent for the pair $(\dd_{\d},\g)$:

$$\cpbDeuxDeux{\g}{\d} =0 + 0  + \cpbCC{{*}}{\g} + 0 + 0 + \cpbFF{\d}{\g} + \cpbGG{{*}}{{*}}{\d}{\g}$$

\noindent for the pair $(\dd_{\g},\g)$:

$$\cpbDeuxDeux{\g}{\g} = 0 +  \cpbB  + 0 +  \cpbDD{{*}}{{*}}  + 0 + \cpbFF{\g}{\g} + \cpbGG{{*}}{{*}}{\g}{\g}$$

It has been shown by L.~Foissy in \cite{Foissy} that the triple $(Dend, Dend, Vect)$ is good by using the explicit description of the free dendriform algebra \cite{JLLdig} (compare with \ref{split}). So there is a rigidity theorem in this case.

\subsection{The triple $(Nil, Nil, Vect)$}\label{NilNil}\index{Nil-algebra} By definition a \emph{$Nil$-algebra}\index{Nil-algebra} is a vector space $A$ equipped with a binary operation $a\cdot b$ such that any triple product is 0:
$$(x\cdot y)\cdot z = 0 = x\cdot (y\cdot z)\ .$$
Hence the operad $Nil$ is binary, quadratic and nonsymmetric. We have $Nil_{1}=\KK, Nil_{2}=\KK$ and $Nil_{n}=0$ when $n\geq 3$. 

By definition a $Nil^c\--Nil$-bialgebra is determined by the following compatibility relation $\between_{nil}$:
$$
\cpbdeuxdeux = \cpbA -\cpblosg - \cpblosd + \cpbloslos \ .
$$

On $Nil(V)= V \oplus V^{\t 2}$ the cooperation $\dd$ is given by $\dd(x)=0$ and $\dd(x\cdot y)= x\t y$. We obviously have $(\dd\t \id)\dd=0=(\id\t\dd)\dd$ as expected. 

From the explicit formula of $\dd$ it follows that $\varphi(V):Nil(V)\to Nil^c(V)$ is an isomorphism. Therefore $(Nil, Nil, Vect)$ is a good triple of operads.

\subsection{The triple $(Nil, Mag, Mag^{0,2,1})$}\label{NilMag} By definition a \emph{$Nil^c\pt Mag$-bialgebra} is determined by the compatibility relation $\between_{Nil}$. It is easy to check that the hypotheses \Hone and \Htwoepi are fulfilled. The primitive operad can be shown to be free (i.e.~magmatic) generated by two ternary operations and one quaternary operation. They are 
$$(x\cdot y)\cdot z\ ,\quad  x\cdot (y\cdot z)\ \textrm{ and } (x\cdot y)\cdot (z\cdot t)\ .$$
We denote by $Mag^{0,2,1}$ the associated operad. Hence there is a good triple
$$(Nil, Mag, Mag^{0,2,1})\ .$$
Its quotient triple is, of course, $(Nil, Nil, Vect)$.

\medskip

\noindent{\bf Question.} Is there a compatibility relation which makes $(Nil^3, Nil^3, Vect)$ into a good triple of operads ? Here $Nil^3$ is the operad of algebras equipped with a binary operation for which any quadruple product is 0.

\subsection{The triple $(Mag^{\infty}_{+}, Mag^{\infty}_{+}, Vect)$} By definition a $Mag^{\infty}_{+}$-algebra is a vector space equipped with an $n$-ary operation $\mu_{n}$ for any integer $n\geq 2$ and also equipped with a unit $1$. Moreover we suppose that
$$\mu_{n}(a_{1} \cdots a_{i-1} 1 a_{i+1} \cdots a_{n}) = \mu_{n-1}(a_{1} \cdots a_{i-1} a_{i+1} \cdots a_{n}) .$$
The compatibility relations are of the following form (indicated here in low dimensions)

$\cpbtroistrois = \cpbtroisID $\quad  and\quad  $\cpbtroisdeux =0= \cpbdeuxtrois$ .

This triple  has been treated in \cite{BurgunderHoltkamp} along the same lines as the general case.

\section{Pre-Lie algebras and a conjectural triple}

\subsection{From Pre-Lie algebras to Lie algebras} Let $A$ be a pre-Lie algebra, cf.~\ref{preLie}. It is immediate to check that the antisymmetrization $[x,y]:=xy-yx$ of this operation is a Lie bracket. Therefore there is defined a  forgetful functor $F:PreLie\alg \to Lie\alg$ which associates to a pre-Lie algebra $A$ the Lie algebra $(A, [-,-])$.

\subsection{The conjectural triple $(??, PreLie, Lie)$} In \cite{Markl07} Markl studied this functor. He mentioned the possible existence of connection with some triple of operads. Indeed it is very likely that there exists a notion of generalized bialgebras $\CC^c\--PreLie$ giving rise to a good triple of operads $(\CC, PreLie, Lie)$. Not only we have to find the operad $\CC$ but also the compatibility relations. The operad $\CC$ would have at least one binary generating operation verifying the symmetry $xy=yx$ and one ternary operation verifying the symmetry $(x,y,z)=(y,z,x)$ (and probably more generators  in higher arity). The compatibility relation between the binary coproduct and the pre-Lie product is probably of Hopf type. In low degrees the dimension of $\CC(n)$ are $(1, 1, 4=3+1, 23, 181)$.

It is also natural to ask the same question for the functor from Lie admissible algebras to Lie algebras.

\subsection{Symmetrizing the pre-Lie product} When symmetrizing the product of an associative algebra we get the notion of Jordan algebra. By an argument on the Poincar\' e series of the operads it can be seen that the forgetful functor from associative algebras to Jordan algebras cannot be part of a triple of operads. 

Let us denote by $\XXX$ the operad governing the type of algebras determined by the symmetrization of a pre-Lie product, that is the operation $x\cdot y := \{x,y\}+\{y,x\}$.  It might happen that the forgetful functor from pre-Lie algebras to $\XXX$-algebras is part of a triple of operads (or at least that there exists a factorization $preLie = F\circ \XXX$ for some \sm $F$).

\section{Interchange bialgebra}\label{IC}\index{IC} We introduce the notion of interchange algebra and interchange bialgebra and we prove that hypothesis \Hone holds. This example is extracted from the paper \cite{Lafont97} by Yves Lafont.

\subsection{Interchange algebra and bialgebra} By definition an \emph{interchange algebra}\index{interchange algebra} is a vector space $A$ equipped with two operations $\circ$ and $\bullet$ satisfying the \emph{interchange law}:
$$(a\circ b)\bullet (c\circ d ) = (a\bullet c) \circ (b\bullet d)\ .$$
This law is quite common in category theory and algebraic topology since it is part of the axioms for a double category  ($\circ=$ horizontal composition, $\bullet=$ vertical composition of bimorphisms). It can be used to prove the commutativity of the higher homotopy groups (Eckmann-Hilton argument).
Observe that this relation is not quadratic. It is closely related to the notions of rack, quandle and left-distributivity.

By definition a \emph{bi-interchange bialgebra} ($IC^c\pt IC$-bialgebra) $\HH$ is both an interchange algebra and an interchange coalgebra with compatibility relations as follows:

$$\cpbDeuxDeux{\circ}{\circ} = \cpbA \qquad ,\qquad  \cpbDeuxDeux{\bullet}{\bullet} = \cpbB , $$
$$\cpbDeuxDeux{\circ}{\bullet} = \cpbGG{\bullet}{\circ}{\bullet}{\circ} \qquad ,\qquad
\cpbDeuxDeux{\bullet}{\circ} = \cpbGG{\circ}{\bullet}{\circ}{\bullet} .$$

\begin{prop}\label{ICH1} The free interchange algebra $IC(V)$ has a natural $IC^c\pt IC$-bialgebra structure.
\end{prop}
\begin{proo} Since the operad $IC$ is set-theoretic, $IC(V)\t IC(V)$ is still an $IC$-algebra and one can define maps $\dd_{\circ}$ and $\dd_{\bullet}:IC(V)\to IC(V)\t IC(V)$ sending $V$ to $0$ and satisfying the compatibility relations, cf.~\ref{Hopfoperad}. By induction they can be shown to satisfy the interchange law (cf.~strategy number (2) in \ref{verification}). This proof can be explained in terms of rewriting systems.
\end{proo}

\section{The $\langle k \rangle$-ary case}\label{kary}

In the preceding examples the generating operations and cooperations were all binary (except sometimes for the primitive operad). In this section we give some examples with $k+1$-ary operations and cooperations for $k\geq 1$. There are many more, not yet explored.

\subsection{Associative $k+1$-ary algebras}\label{Askary} Let $k$ be an integer greater than or equal to 1. Let $\CC$ and $\AA$ be two operads generated by $k+1$-ary operations. Here are two important examples taken from \cite{Gnedbaye97}. 

A \emph{totally associative $k+1$-ary algebra} is a vector space $A$ equipped with a $k+1$-ary operation $(a_{0}\cdots a_{k})$ satisfying the relations
$$\big((a_{0}\cdots a_{k})a_{k+1}\cdots a_{2k}\big)= \big(a_{0}\cdots (a_{i}\cdots a_{i+k})\cdots a_{2k}\big)$$
for any $i=0, \ldots, k$. The operad is denoted $tAs^{\langle k \rangle}$.

A \emph{partially associative $k+1$-ary algebra} is a vector space $A$ equipped with a $k+1$-ary operation $(a_{0}\cdots a_{k})$ satisfying the relations
$$\sum_{i=0}^{k}(-1)^{ki}\big(a_{0}\cdots (a_{i}\cdots a_{i+k})\cdots a_{2k}\big)= 0$$
for any $i=0, \ldots, k$. These two operads were shown to be Koszul dual to each other by V.~Gnedbaye in \cite{Gnedbaye97}.

In this context, a $\CC^c\--\AA$-bialgebra (or generalized bialgebra) is a vector space $\HH$ equipped with a $\CC$-coalgebra structure, a $\AA$-algebra structure, and each pair $(\dd,\mu)$ of a generating operation and a generating cooperation is supposed to satisfy a distributive compatibility relation. Observe that in this case the $\Phi_{1}$-term is an element of the group algebra $\KK[S_{k+1}]$. 

Here is an example for $k=2$, denoted $\between_{tAs^{\langle 2 \rangle}}$:

$\cpbtroistrois = \cpbtroisID + \cpbtroisAB + \cpbtroisAC + \cpbtroisBA $

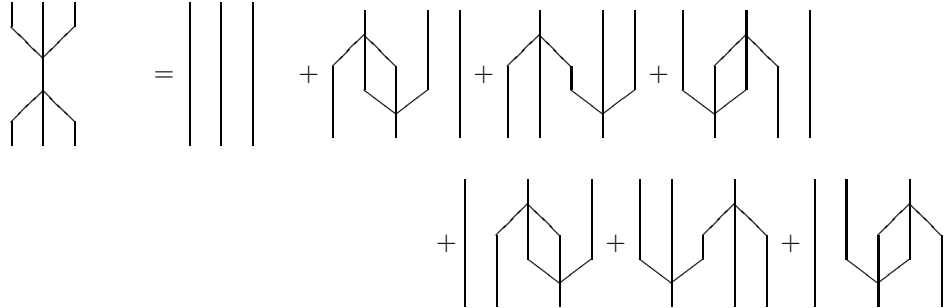
\hfill $+ \cpbtroisBC + \cpbtroisCA + \cpbtroisCB$

We let the reader figure out the similar relation for higher $k$'s. 

\subsection{The triple $(tAs^{\langle k\rangle}, tAs^{\langle k\rangle}, Vect)$}\label{k-ary} By definition a 
$tAs^{\langle k\rangle}$-bialgebra is a vector space $\HH$ equipped with a structure of $tAs^{\langle k\rangle}$-algebra, a structure of $tAs^{\langle k\rangle}$-coalgebra, related by the compatibility relation $\between_{tAs^{\langle k\rangle}}$ described above. For $k=1$ this is the nonunital infinitesimal compatibility relation.

The free totally associative $(k+1)$-ary algebra over $V$ is $tAs^{\langle k\rangle}(V)= \bigoplus_{n\geq 0}V^{\t 1+kn}$. We put a structure of $tAs^{\langle k\rangle}$-coalgebra on it by dualizing the natural basis. Then it is easy to prove by induction that the compatibility relation is precisely $\between_{tAs^{\langle k\rangle}}$. The map $\varphi(V)$ is an isomorphism, hence the triple 
$$(tAs^{\langle k\rangle}, tAs^{\langle k\rangle}, Vect)$$
 is good and the rigidity theorem holds.

\subsection{The triple $(tCom^{\langle k\rangle}, tCom^{\langle k\rangle}, Vect)$}\label{Comk-ary}
By definition a  $tCom^{\langle k\rangle}$-algebra is a totally associative $(k+1)$-ary algebra which is \emph{commutative} in the sense
$$(a_{0}\cdots a_{k}) = (a_{\ss(0)}\cdots a_{\ss(k)})$$
for any permutation $\ss \in S_{k+1}$.

\noindent {\bf Exercise}. Find the compatibility relation which gives a good triple of operads 
$(tCom^{\langle k\rangle}, tCom^{\langle k\rangle}, Vect)$.

It would be also interesting to work out the cases $pAs^{\langle k\rangle}$ and $pCom^{\langle k\rangle}$, and also the triple 
$$( tCom^{\langle k\rangle},  tAs^{\langle k\rangle},  tLie^{\langle k\rangle}),$$
 which is $(Com, As, Lie)$ for $k=1$.


\chapter{Duplicial bialgebras}\label{ch:duplicial}  In this chapter we study in details the duplicial algebras and the duplicial bialgebras. They are defined by two associative operations verifying one more relation. The coproduct is going to be coassociative. We show that there exists a good triple of operads
$$(As, Dup, Mag)\ ,$$
which is quite peculiar since the three operads are binary, quadratic, nonsymmetric, set-theoretic and Koszul. 

In order to prove that the operad $Dup$ is Koszul, we compute its dual and construct the chain complex giving rise to the homology of duplicial algebras. It turns out that it is the total complex of a certain chain complex whose horizontal (resp. vertical) components are of Hochschild type. The relationship between duplicial algebras and quantum electrodynamics can be found in \cite{BF} and \cite{Frabetti07}.

\section{Duplicial algebra} 

\subsection{Definition} A \emph{duplicial algebra}\index{duplicial algebra} (or $Dup$-algebra for short)\index{Dup-algebra} $A$  is determined by two binary operations $A\t A \to A$ called \emph {left} $(x,y)\mapsto x\g y$ and \emph{right} $(x,y)\mapsto x\d y$ respectively, satisfying the following three relations
 \begin{eqnarray*}
 (x\g y)\g z  &=& x\g (y\g z)\ , \\
(x\d y)\g z  &=& x\d (y\g z) \ ,\\
(x\d y)\d z &=& x\d (y\d z) \ .
\end{eqnarray*}
So the two operations left and right are associative. From this definition it is clear that the operad $Dup$ is binary, quadratic, nonsymmetric and set-theoretic. 

In order to describe the free duplicial algebra (or, equivalently, the operad), we need to introduce the 
planar binary rooted trees.

\subsection{Planar binary trees}\label{pbtree} By definition a \emph{planar binary rooted tree}\index{planar binary tree} (we simply say planar binary tree, or p.b.~tree for short) is a finite planar rooted tree with vertices which have two inputs and one output. The edge with no successor is called the \emph{root}. The edges with no predecessors are called the \emph{leaves}. The set of 
planar binary rooted trees with $n$ leaves is denoted $PBT_{n}$:

\begin{displaymath}
PBT_1 = \{\  \vert \ \}\quad  , \quad PBT_2=\big\{ \arbreA\big\}\quad , \quad 
{PBT_3=\Big\{ \arbreAB}\ ,\  { \arbreBA }\Big\}
\end{displaymath}
\begin{displaymath}
{PBT_4=\Bigg\{ \arbreABC}\ ,\  {{} \arbreBAC },\  {{} \arbreACA },\  {{} \arbreCAB },\  {{} \arbreCBA }\Bigg\}
\end{displaymath}
Recall that the number of elements in $PBT_{n+1}$ is the \emph{Catalan number} \index{Catalan number} $c_{n}= \frac{(2n)!}{n!(n+1)!}$ and the generating series is
$$c(t) := \sum_{n\geq 1}
c_{n-1}t^n=\frac{1}{2}(1-\sqrt{1-4t})\ .$$
  On these trees one can perform the following two kinds of grafting: the operation Over and the operation Under. By definition the operation \emph{Over}\index{Over operation}, denoted $t/s$, consists in grafting the tree $t$ on the first leaf of $s$. Similarly,  the operation \emph{Under}\index{Under operation}, denoted $t\backslash s$ consists in grafting the tree $s$ on the last leaf of $t$.
Observe that for $t\in PBT_{p+1}$ and $s\in PBT_{q+1}$ we have $t/s\in PBT_{p+q+1}, t\backslash s\in PBT_{p+q+1}$.

$${\xymatrix@R=2pt@C=2pt{
&*{}&*{}&*{}&*{}&*{}&*{}&*{}\\
&&& t &&&&\\
t/s =& &&*{}\ar@{-}[uurr]\ar@{-}[uull]&s&&&\\
&&&&*{}\ar@{-}[ul]\ar@{-}[dd]\ar@{-}[uuurrr]&&&\\
&&&&&&&\\
&&&&&&&\\
}}
\qquad
{\xymatrix@R=2pt@C=2pt{
&*{}&*{}&*{}&*{}&*{}&*{}&*{}\\
&&&&& s &&\\
t\backslash s =& &&& t &*{}\ar@{-}[uull]\ar@{-}[uurr]&&\\
&&&&*{}\ar@{-}[uuulll]\ar@{-}[dd]\ar@{-}[ur]&&&\\
&&&&&&&\\
&&&&&&&\\
}}
$$

The Over and Under operations, together with their properties appear in the work of  C.~Brouder and A.~Frabetti \cite {BF}. 
We will also use another type of grafting, denoted $t\vee s$ which consists in creating a new root and grafting the two trees to this root. Hence we have $t\vee s \in PBT_{p+q+2}$. Observe that any tree $t$ (except $\vert$) is uniquely determined by its left part $t^l$ and its right part $t^r$ so that $t = t^l\vee t^r$. The \emph{left comb}\index{left comb} is the tree $\textrm{comb}^l_{n}$ defined inductively as $\textrm{comb}^l_{1}= \arbreA$  and $\textrm{comb}^l_{n}= \arbreA / \textrm{comb}^l_{n-1}$.

\begin{prop}\label{freeDup} The free duplicial algebra on one generator is spanned by the set of planar binary trees and the left (resp. right) operation is induced by the Over operation $t/s$ (resp. Under operation $t\backslash s$). Hence the space of nonsymmetric $n$-ary operations is  $Dup_{n}= \KK[PBT_{n+1}]$. 
\end{prop}
\begin{proo} First, we verify immediately that $\bigoplus_{n\geq 1}\KK[PBT_{n+1}]$ equipped with the Over and Under operations is a duplicial algebra generated by the unique element $ \arbreA$ of $PBT_{2}$. 
Moreover we see that, for any p.b.~tree $t=t^l\vee t^r$ we have $t = t^l/\arbreA \backslash t^r$.

Let us show that $\bigoplus_{n\geq 1}\KK[PBT_{n+1}]$ satisfies the universal condition. Let $A$ be a duplicial algebra and let $a\in A$. We define a map $\phi: \bigoplus_{n\geq 1}\KK[PBT_{n+1}]\to A$ inductively by 
$\phi(\arbreA)= a$ and $\phi(t)= \phi(t^l\vee t^r) = \phi(t^l)\d a \g \phi(t^r)$. It is straightforward to check that this map is a duplicial morphism (same proof as in the dendriform case, see \cite{JLLdig} Proposition 5.7). Since we have no other choice for its value, it is the expected universal extension map.
\end{proo}

\noindent{\bf Remark.} The same structure in the category of sets, instead of vector spaces, has been investigated by Teimuraz Pirashvili  in \cite{Pira03} under the name \emph{duplex}.

\subsection{Relationship with other algebraic structures} The operad of duplicial algebras is related to several other algebraic operads by the following morphisms:
$$
{\xymatrix {
                &               &                        &Dias\ar[rr]&            &Nil \\
2mag\ar[r]&2as\ar[r]&Dup\ar[ur]\ar[dr]&                &            &\\
                 &             &                          &As^2\ar[r]&As\ar[r]&Com \\
}}
$$
\section{Duplicial bialgebra}\label{duplicialbialgebra} 

\subsection{Definition} By definition a $As^c\--Dup$-bialgebra, also called \emph{duplicial bialgebra}\index{duplicial bialgebra}, is a vector space $\HH$ equipped with a duplicial algebra structure $\g, \d$, a coassociative coalgebra structure $\dd$, and the compatibility relations are of nonunital infinitesimal type for both pairs $(\dd, \g)$ and $(\dd,\d)$:

$$\cpbDeuxDeux{\g}{*{}} = \cpbA + \cpbCC{*{}}{\g} +  \cpbEE{*{}}{\g}$$

$$\cpbDeuxDeux{\d}{*{}} = \cpbA + \cpbCC{*{}}{\d} +  \cpbEE{*{}}{\d}$$

Observe that this is a nonsymmetric bialgebra type (hence a nonsymmetric prop).

 \begin{prop}\label{Dupbialgebra} The free duplicial algebra $Dup(V)$ is a $As^c\pt Dup$-bialgebra.
 \end{prop}
 \begin{proo} 
 Let us first define the cooperation 
 $$\dd : Dup(\KK)=\bigoplus_{n} \KK[PBT_{n+1}] \to \bigoplus_{n} \KK[PBT_{n+1}]  \t \bigoplus_{n} \KK[PBT_{n+1}] .$$
  For any $t\in PBT_{n+1}$ we define 
 $$\dd(t)= \sum_{1\leq i \leq n-1}\dd_{i}(t) = \sum_{1\leq i \leq n-1} r_{i}\t s_{i}$$
 as follows. Let us number the leaves of $t$ from left to right by the integers $0, 1, \ldots , n$. For any $i=1, \ldots, n-1$ we consider the path going from the leaf number $i$ to the root. The left part of $t$ (including the dividing path) determines the tree $r_{i}$ and the right part of $t$ (including the dividing path)  determines the tree $s_{i}$. In particular $\dd(\arbreA) = 0$.
 
 Example for $i=2$:
 
$$t={\vcenter{\xymatrix@R=1pt@C=1pt{
0&&1&&2&&3&&4\\
&&&& &&&&\\
&*{}\ar@{-}[ur] &&&&&&*{}\ar@{-}[ul] & \\
&&*{}\ar@{-}[uurr] &&&&&&\\
&&&&&&&&\\
&&&&*{}\ar@{-}[uuuurrrr] \ar@{-}[uuuullll] \ar@{-}[d] &&&&\\
&&&& &&&&
}}} , \qquad  r_{2}= \arbreAB,\quad  s_{2}= \arbreBA\ .$$

 It is immediate to verify, by direct inspection, that $\dd$ is coassociative. Let us prove that 
 $$\dd(x/y) = x \t y + x_{(1)}\t x_{(2)}\g y + x \g y_{(1)}\t y_{(2)}\ ,$$
 under the notation $\dd(x) = x_{(1)}\t x_{(2)}$. Let $x\in PBT_{p+1}, y\in PBT_{q+1}$, so that $x/y \in PBT_{p+q+1}$. The element $\dd(x/y)$ is the sum of three different kinds of elements: either the dividing path starts from a leaf of $x$ not being the last one,  or starts from the last  leaf of $x$,  or starts from a leaf of $y$. In the first case we get  $x \g y_{(1)}\t y_{(2)}$, in the second case we get $ x \t y $, in the third case we get $ x_{(1)}\t x_{(2)}\g y $. The proof for $\d$ is similar.
 
 To prove that $Dup(V)$ is a $As^c\pt Dup$-bialgebra for any $V$ it suffices to extend $\dd$ to 
 $\bigoplus_{n} \KK[PBT_{n+1}]\t V^{\t n} $ by
 $$\dd(t; v_{1}\cdots v_{n})= \sum_{1\leq i\leq n-1}(r_{i}; v_{1}\cdots v_{i})\t (s_{i}; v_{i+1}\cdots v_{n})$$
 and use the property that $\To(V)$ is a n.u.i.~bialgebra (cf. \cite{LRstr} of \ref{AsAsVect}).
 \end{proo}
 
 \subsection{Remark} We could also prove this Proposition by using the inductive method described in \ref{verification} (it is a good exercise !).
 
 \begin{prop}\label{AsDupepi} The prop $As^c\pt Dup$ satisfies the hypothesis \texttt{(H2epi)}.
  \end{prop}
  
   \begin{proo} Since we are dealing with a nonsymmetric bialgebra type, it suffices to look at $Dup(\KK)=\bigoplus_{n} \KK[PBT_{n+1}]$. The map $\varphi: Dup_{n}=\KK[PBT_{n+1}]\to \KK = As^c_{n}$ is given by $\varphi(t)= \aa$ where $\aa$ is a scalar determined by the equation 
   $$\dd^{n-1}(t) = \aa \arbreA \t \cdots \t \arbreA.$$
    Here $\dd^{n-1}$ stands for the iterated comultiplication. From the explicit description of $\dd$ it comes immediately: $\dd^{n-1} (t)= \arbreA \t \cdots \t \arbreA$. Hence $\aa=1$ and the map $\varphi$ is given by $\varphi(t)=1$.
   
  Define a map $s_{n}: As_{n}=\KK \to \KK[PBT_{n+1}]$ by $s_{n}(1)={\textrm comb}^l_{n}$, where ${\textrm comb}^l_{n}$ is the left comb. It is immediate to check that $s_{.} $ induces a coalgebra map $s(V): As(V) \to Dup(V)$ which is a splitting to $\varphi(V)$. Hence hypothesis 
   \texttt{(H2epi)} is fulfilled.
   \end{proo}
   
   As a consequence  the triple $(As, Dup, \Prim_{As} Dup)$ is a good triple and it satisfies the structure Theorem over any field $\KK$ by \ref{thm:structurens}. Let us now identify the operad $ \Prim_{As} Dup$.
   
    \begin{thm}\label{primDup} The primitive operad $\Prim_{As} Dup$ of the bialgebra type $As^c\pt Dup$ is the magmatic operad $Mag$ and the functor 
    $$F: Dup\alg \to Mag\alg,\quad F(A, \g, \d) = (A, \cdot)$$ 
 is determined by
    $$ x\cdot y := x\g y - x\d y\ .$$
  \end{thm}
  \begin{cor}\label{AsDupMag} There is a good triple of operads
  $$(As, Dup, Mag)\ .$$
  \end{cor}

Before entering the proof of the Theorem and its Corollary we prove some useful technical Proposition.

\begin{prop}\label{extendingmag} Let $(R, \cdot)$ be a magmatic algebra. On $As(R)=\To(R)$ we define the operation $a\d b$ as being the concatenation (i.e. $\d = \t$) and we define the operation $a\g b$ by $a\g b = a\cdot b + a\d b$ where the operation $a\cdot b$ is defined inductively as follows:
\begin{eqnarray*}
(r\t a)\cdot b &=& r \t (a\cdot b)\\
r\cdot (s\t b) &=& (r\cdot s)\cdot b - r\cdot (s\cdot b) + (r\cdot s)\t bÊ .
\end{eqnarray*}
Then $(As(R), \g, \d)$ is a duplicial bialgebra with the deconcatenation as coproduct.
\end{prop}
\begin{proo} The last relation of duplicial algebra (associativity of $\d$) is immediate. The other two are proved by a straightforward induction argument on the degree. The compatibility relation for the pair $(\dd, \d)$ is well-known (cf.~\ref{AsAsVect}). The compatibility relation for the pair $(\dd, \g)$ is proved by induction.
\end{proo}

\subsection{Proof of Theorem \ref{primDup} and  Corollary \ref{AsDupMag}} Applying Proposition \ref{extendingmag} to the free magmatic algebra $R=Mag(V)$, we get a duplicial algebra $As(Mag(V))$. The inclusion map 
$$V= As(Mag(V))_{1} \mono As(Mag(V))$$
induces a $Dup$-map $Dup(V)\to As(Mag(V))$. From the construction of $As(Mag(V))$ it follows that this map is surjective.

From Proposition \ref{freeDup} it follows that $\dim Dup_{n} = c_{n}$.
It is also known that $\dim (As\circ Mag )_{n}= c_{n}$ because
$$f^{As}f^{Mag}(t) = \frac{c(t)}{1-c(t)}=  \frac{c(t)-t}{t}= \sum_{n\geq 1}c_{n}t^n \ .$$
(Use the identity $c(t)^2 - c(t) + t = 0$). Therefore the surjective map $Dup_{n}\to (As\circ Mag )_{n}$ is an isomorphism. Hence $Dup\to As\circ Mag$ is an isomorphism. Since the comultiplication in $As(Mag(V))$ is the deconcatenation, its primitive part is $Mag(V)$. It follows that $\Prim_{As}Dup(V)=Mag(V)$ as expected. 

Corollary \ref{AsDupMag} follows from Proposition \ref{AsDupepi} and Proposition \ref{extendingmag}.
\hfill $\square$

\begin{cor}\label{Dupfreeness} As an associative algebra for the product $\d \,= /$ the space $Dup(V)$ is free over $Mag(V)$.
\end{cor}
\begin{proo} By the structure theorem for $As^c\--Dup$-bialgebras we know that there is an isomorphism $Dup(V)\cong \To^c(\Prim_{As}Dup(V))$. Because of our choice of $s$, it turns out that the $As$-structure of the $As^c\--As$-bialgebra $ \To^c(\Prim_{As}Dup(V))$ corresponds to $/ =\ \d$ , cf.~\ref{AsAsVect}. Hence $Dup(V)$ is free for the operation $\d$.
\end{proo}

We have an extension of operads 
$$As\mono Dup\epi Mag$$
in the sense of \ref{extension}. It is even a split extension.

\subsection{Remark on the map $f: Mag\to Dup$}  Let us write $PBT_{n+1}$ as a union of two disjoint subsets $PBT_{n+1}^{a}$ and $PBT_{n+1}^{b}$, where $PBT_{n+1}^{a}$ is made of the trees of the form $\vert \vee t$ for $t\in PBT_{n}$. From the definition of $f_{n}:Mag_{n}\to Dup_{n}=\KK[PBT_{n+1}]$ and Theorem \ref{primDup} it follows that the composition of maps 
$$\KK[PBT_{n}]= Mag_{n}\to Dup_{n}=\KK[PBT_{n+1}]\epi \KK[PBT_{n+1}/PBT_{n+1}^b] =\KK[PBT_{n+1}^{a}]\cong \KK[PBT_{n}]$$
is an isomorphism. It is a nontrivial isomorphism, given in low dimension by:
\begin{eqnarray*}
\arbreA & \mapsto & - \arbreA\\
\arbreAB  & \mapsto & - \arbreBA \\
\arbreBA  & \mapsto & + \arbreBA  - \arbreAB \\
\arbreABC  & \mapsto & - \arbreCBA \\
\arbreBAC  & \mapsto & - \arbreCBA - \arbreACA \\
\arbreACA  & \mapsto & - \arbreCBA + \arbreCAB \\
\arbreCAB  & \mapsto & - \arbreCBA + \arbreACA + \arbreBAC - \arbreABC \\
\arbreCBA  & \mapsto & - \arbreCBA + \arbreACA - \arbreCAB - \arbreABC \\
\end{eqnarray*}

\section{Explicit PBW-analogue isomorphism for $Dup$} When $\HH= Dup(V)$ the isomorphism
$\HH\cong As^c(\Prim \HH)$ becomes $Dup(V)\cong \To^c(Mag(V))$. Therefore we should be able to write any linear generator of $Dup_{n}$ as a tensor of elements in $Mag_{k}, k\leq n$. Since we choose the operation $\d$ to split the map $\varphi$, we can replace the tensor by $\d$ and write an equality in $Dup(V)$ (analogous to what we did in the classical case, see \ref{explicit}). In low dimension it gives the following equalities:

\medskip

\begin{tabular}{| c c c c c c c  |}
\hline
      $Dup$      &  &  $T^1Mag $& & $T^2Mag$ & & $T^3Mag$ \\
\hline
$x$&=&$x$&&&&\\
\hline
$x\d y $&=&$ 0 $&$+ $&$ x\d y$ && \\
$x\g y $&=& $x\cdot y $&$+ $&$ x\d y $&&\\
\hline
$x\d y\d z $   &=& $0$ & $+$ & $0$                 & $+$  & $ x \d y \d z $ \\
$(x\g y)\d z $&=& $0$ & $+$ & $(x\cdot y)\d z $  & $+$  & $ x \d y \d z $ \\
$x\d y \g z $  &=& $0$ & $+$ & $x\d (y\cdot z) $  & $+$  & $ x \d y \d z $ \\
$x \g(y \d z) $&=& $(x\cdot y)\cdot z $ & $+$ & $(x\cdot y)\d z $  & $+$  & $ x \d y \d z $ \\
                     &   & $- x\cdot (y\cdot z)$ &      &                             &         &                        \\
$x\g y \g z $  &=& $(x\cdot y)\cdot z $ & $+$ & $(x\cdot y)\d z $ & $+$  & $ x \d y \d z $ \\
                     &  &                               &        & $+x\d (y\cdot z) $ &        & \\
\hline

\end{tabular}

\medskip

These formulas are consequences of Proposition \ref{extendingmag}.

\section{Koszulity of the operad $Dup$}

\subsection{Dual operad}\label{dualDup} Since the operad $Dup$ is quadratic, it admits a dual operad, denoted $Dup^{!}$, cf.~\cite{GK}. The $Dup^{!}$-algebras are duplicial algebras which satisfy the following additional relations:
$$(x\g y)\d z = 0 \qquad \textrm{and} \qquad 0= x\g (y\d z)\ .$$
This is easy to check from the conditions given in \cite{JLLdig} Appendix  for a nonsymmetric operad to be Koszul. See also \cite{LRcha}.

The free $Dup^{!}$-algebra is easy to describe (analogous to the free diassociative algebra, see \cite{JLLdig}). We have $Dup^{!}_{n}= \KK^n$, where the $i$th linear generator corresponds to 
$$ \underbrace{x\d x\d \cdots x}_{{i-1}}\d x \g\underbrace{x \cdots \g x \g x}_{n-i}\ .$$

\subsection{The total bicomplex $C^{Dup}_{**}$} Let $A$ be a duplicial algebra. We define a chain bicomplex $C^{Dup}_{**}(A)$ as follows: $C^{Dup}_{pq}(A)= A^{\t p+q+1}$ and 
$$d^h(a_{0}\cdots a_{p+q})= \sum_{i=0}^{p-1}(-1)^{i}a_{0}\cdots(a_{i}\d a_{i+1}) \cdots a_{p+q} \ , $$
$$d^v(a_{0}\cdots a_{p+q})= \sum_{j=p}^{p+q-1}(-1)^{j}a_{0}\cdots(a_{j}\g a_{j+1}) \cdots a_{p+q} \ . $$
The relation $d^h d^h =0$ follows from the associativity of the operation $\d$. The relation $d^v d^v =0$ follows from the associativity of the operation $\g$. The relation $d^h d^v + d^v d^h=0$ follows from the relation entwining $\g$ and $\d$.

\bigskip

 $C_{**}^{Dup}(A)$ \xymatrix{
\cdot \ar[d] & & & \\
A^{\t 3}\ar[d]_{\Id\t \g} & \ar[l]\ar[d]  \cdot & & \\
A^{\t 2}\ar[d]_{\g}  & \ar[l]_{\d \t \Id} A^{\t 3}\ar[d]_{-\Id\t \g}  & \ar[l] \cdot \ar[d]& \\
A & \ar[l]_{\d} A^{\t 2} &   \ar[l]_{\d\t \Id} A^{\t 3} & \ar[l] \cdot 
}

\bigskip

\begin{prop}\label{totalcplx} The (operadic) homology of a duplicial algebra $A$ is the homology of the total complex of the bicomplex
$C^{Dup}_{**}(A)$ up to a shift.
\end{prop}
\begin{proo} The operadic chain complex of a duplicial algebra is given by 
$$C^{Dup}_{n}(A) = (Dup^!_{n})^* (A)$$
and the differential $d$ is the unique coderivation which extends the $Dup$-products.

From the description of the operad $Dup^!$, cf.~\ref{dualDup}, we check immediately that $Tot\ C^{Dup}_{**}(A) = C^{Dup}_{*}(A) $. The fact that $d^h+d^v$ identifies to the operadic differential is also immediate.
\end{proo}

\begin{thm} The operad $Dup$ is a Koszul operad.
\end{thm}
\begin{proo} Let us recall some facts about Hochschild homology of non-unital algebras. Let $R$ be a non-unital algebra and let $M$ be a right $R$-module. The Hochschild complex of $R$ with coefficients in $M$ is:

$$C_{*}(R,M): \qquad \to \cdots M\t R^{\t n} \buildrel{b'}\over{\longrightarrow} M\t R^{\t n-1} \to \cdots \to M$$
where $b'(a_{0},\ldots , a_{n}) = \sum_{i=0}^{i=n-1}(-1)^{i}(a_{0},\ldots , a_{i}a_{i+1},\ldots  ,a_{n}) $ and $a_{0}\in M, a_{i}\in R$. The homology groups are denoted by $H_{*}(R,M)$. If $R$ is free over $W$, i.e.~$R=\To(W)$, then one can prove the following (cf.~for instance \cite{HC}):
\begin{eqnarray*}
H_{0}(R,M)& = & M/MR,\\
H_{n}(R,M)& = & 0\quad  \textrm{otherwise.}
\end{eqnarray*}

In order to prove the theorem it suffices to show that the Koszul complex is acyclic, or equivalently that the $Dup$ homology of the free duplicial algebra $Dup(V)$ is 
$$H_{1}^{Dup}(Dup(V))= V, \quad \textrm{and}\quad  H_{n}^{Dup}(Dup(V))= 0$$
 for $n\geq 2$.
 
Since by Proposition \ref{totalcplx} the chain complex of the duplicial algebra $A$  is the total complex of a bicomplex, we can use the spectral sequence associated to this bicomplex to compute it:
$$E^2_{pq}= H_{q}^v H_{p}^h (C^{Dup}_{**}(A) ) \Rightarrow H_{p+q+1}^{Dup}(A)\ .$$
Since $A:=Dup(V)$ is free over $Mag(V)$ as an associative algebra for $\d$ (cf. \ref{Dupfreeness}) and since the horizontal complex is the Hochschild complex (for $\d$) with coefficients in $Mag(V)$, we get
$$H_{q}^v (C^{Dup}_{p*}(A))=0, \ \textrm{for}\ q\geq 1\ \textrm{and}\ H_{0}^v (C^{Dup}_{p*}(A))=A^{\t p}\t Mag(V)\ .$$
Hence the complex $(E^1_{0*},d^1)$ is the Hochschild complex (for $\g$) of $A$ with coefficients in
$Mag(V)$. Its homology is 

\begin{eqnarray*}
E^2_{00}& = & Mag(V)/Mag(V)A= V ,\\
E^2_{0q}& = & 0\quad  \textrm{otherwise.}
\end{eqnarray*}

Hence the spectral sequence tells us that $H_{n}^{Dup}(A)=0$ for $n>1$ and that $H_{1}^{Dup}(A)=V$. So we can deduce that $Dup$ is a Koszul operad.
\end{proo}

\subsection{Alternative proof} (Bruno Vallette, private communication) Since the operad $Dup$ is set-theoretic, one can apply the poset method of Vallette \cite{Vallette05, CV} to prove its Koszulity. Here the poset is as follows. Let us fix an integer $n$. The poset $\Pi_{Dup}(n)$ is made of ordered sequences $\row t1k$ of p.b.~trees such that $\sum_{i=1}^k |t_i |= n$ and $|t_i |\geq 1$. The covering relations defining the poset structure are
$$\row t1{k+1} \to \row {t'}1k$$
if and only if the second sequence is obtained from the first by replacing two consecutive trees $t_i, t_{i+1}$ either by $t_i / t_{i+1}$ or by $t_i\backslash  t_{i+1}$. One can show that the poset is ``Cohen-Macaulay" by methods of \cite {CV}, and so, by \cite{Vallette05}, that the associated chain complex is acyclic (except in top dimension). In fact the top dimension homology group is $Dup^!_n$. This computation proves the Koszulity of the operad $Dup$.

\subsection{Question} Since $Mag^!=Nil$ and $As^!=As$ the construction proposed in \ref{Koszul} suggests the existence of a good triple of operads $(Nil, Dup^!, As)$. Does it exist  ? 

\section{On quotients of $Dup$}\label{AsDupPreLie} 

\subsection{The triple $(As, As^2,As)$}\label{asas2} An \emph{$As^2$-algebra}\index{$As^2$-algebra} is, by definition, a vector space $A$ equipped with two operations denoted $a\cdot b$ and $a*b$ such that the associativity relation
$$(a\circ_{1}b)\circ_{2} c = a\circ_{1}(b\circ_{2} c )$$
holds for any value of $\circ_{i}$ (i.e.~either equal to $\cdot$ or to $*$). In a $As^c\--As^2$-bialgebra we choose the compatibility relations to be the nonunital infinitesimal relation. It is immediate to verify that $\Prim_{As}As^2 = As$, because $x\circ y := x\cdot y - x*y$ is also associative. Therefore we get a good triple of operads 
$$(As, As^2, As),$$
which appears as a quotient of $(As, Dup, Mag)$. 
The operad $As^2$ appears in many places in the literature. For instance in \cite{Goncharov} there is a notion of $Com^c\pt As^2$-bialgebra where the compatibility relation for the pair $(\dd, \cdot)$ is Hopf (unital version) and for the pair $(\dd, *)$ it is
$$\cpbDeuxDeux{{*}}{*{}} = \cpbGG{*{}}{\cdot }{*{}}{{*}}+\cpbGG{*{}}{{*}}{*{}}{\cdot}\ .$$

\subsection{ $DupPreLie$-algebras} Let $DupPreLie$ be the operad which is  a quotient of $Dup$ by the relation
$$(x\g y)\d z - x\g (y\d z) = (x\g z)\d y - x\g (z\d y) \ .$$
This operad is still binary and quadratic, but is not nonsymmetric anymore since the added relation does not keep the variables in the same order. 

\begin{lemma} In any duplicial algebra the following equality holds:
$$(x\cdot y) \cdot z - x\cdot (y\cdot z) = (x\g y)\d z - x\g (y\d z)\ .$$
\end{lemma}
\begin{proo} Recall that $x\cdot y := x\g y - x \d y$. It is an immediate consequence of the relations.
\end{proo}
\begin{prop} There is a good triple
$$(As, DupPreLie, \Prim_{As}DupPreLie).$$
\end{prop}
\begin{proo} It follows from the Lemma that the relator defining $DupPreLie$ as a quotient of $Dup$ is a primitive operation. The claim is therefore a consequence of Proposition \ref{prop:quotient}.
\end{proo}

\noindent {\bf Question.} Do we have $ \Prim_{As}DupPreLie = PreLie$ ?

\section{Shuffle bialgebras} For duplicial algebras we used the grafting on the first leaf and the grafting on the last leaf. However there is a more subtle structure which consists in using the grafting operations on any leaf. Strictly speaking it does not give an operad because, for a given integer $i$,  the operation ``grafting on the $i$th leaf" exists only when the elements have high enough degree. This ``grafting algebra" structure has been studied in details by Mar\'\i a Ronco in \cite{Ronco06}, where she proves the analog of a structure theorem in this setting.

\section{The triple $(Dup,Dup, Vect)$}\label{DupDup} 

\subsection{Biduplicial bialgebra}  By definition an  \emph{infinitesimal biduplicial bialgebra}\index{infinitesimal biduplicial} is a $Dup^c\pt Dup$-bialgebra with the following compatibility relations:

$$\cpbDeuxDeux{\d}{\d} = \cpbA + \cpbCC{\d}{\d} +\cpbEE{\d}{\d} $$

$$\cpbDeuxDeux{\d}{\g} =  \cpbEE{\g}{\d} $$

$$\cpbDeuxDeux{\g}{\d} =  \cpbCC{\g}{\d} $$

$$\cpbDeuxDeux{\g}{\g} = \cpbA +  \cpbCC{\g}{\g} +\cpbEE{\g}{\g} $$
 
\begin{prop} Equipped with the dual basis coalgebra structure, the free duplicial algebra $Dup(V)$ is a natural infinitesimal biduplicial bialgebra.  So Hypothesis \Hone is fulfilled.
\end{prop}

\begin{proo} We use the explicit description of $Dup(V)$ in terms of planar binary trees given in \ref{freeDup}. We check the case
$$\dd_{\g}(x\g y) = x\t y + x_{(1)}^{\g}\t  x_{(2)}^{\g}\g y + x\g  y_{(1)}^{\g}\t  y_{(2)}^{\g},$$
where $\dd_{\g}(x)= x_{(1)}^{\g}\t  x_{(2)}^{\g}$. 

First we describe $\dd_{\g}(x)$ explicitly when $x$ is a p.b. tree. Along the right edge of $x$ we can cut between legs to obtain two trees denoted $x_{(1)}^{\g}$ and $   x_{(2)}^{\g}$ such that $x = x_{(1)}^{\g}\backslash   x_{(2)}^{\g}$. Then $\dd_{\g}(x)$ is the sum of these $x_{(1)}^{\g}\t   x_{(2)}^{\g}$ for all possible cuts.

Let $x$ and $y$ be two p.b. trees. Since $x\g y = x\backslash y$ in the free algebra, the cuts on the right edge of $x\d y$ are of three different types:

-- either a cut in $x$,

-- or a cut separating $x$ from $y$,

-- or a cut in $y$.

The first type of cuts gives the summands of $x_{(1)}^{\g}\t  (x_{(2)}^{\g}\g y)$; the second type of cuts gives $x\t y$; the third type of cuts gives the summands of $ (x\g  y_{(1)}^{\g})\t  y_{(2)}^{\g}$. So we are done with this case.

Since in $x\d y= x/y$ the right edge is the same as the right edge of $y$ the cuts to obtain $\dd_{\g}(y)$ are exactly the cuts of $y$. Therefore we get $\dd_{\g}(x\d y) = x \d y_{(1)}^{\g}\t  y_{(2)}^{\g}$ as expected. 

The proof of the other two cases (involving $\dd_{\d}$) is analogous.
\end{proo}

\begin{prop} The map $\varphi(V):Dup(V) \to Dup^c(V)$ identifies the basis of $Dup_{n}$ with its dual, which is a basis of $Dup^c_{n}$. So Hypothesis \Htwoiso is fulfilled.
\end{prop}

\begin{proo} Let $t$ and $s$ be p.b.~trees. We need to compute $\dd_{t(s)}$, which is of the form $\lambda\, x\t \cdots \t x$, where $\lambda$ is a coefficient and $x= \arbreA$ is the generator of $Dup(\KK)$. From the compatibility relations it is immediately seen that $\lambda = 1$ if $t=s$ and that $\lambda=0$ is $t\neq s$. 
\end{proo}

\begin{cor} The triple $(Dup, Dup, Vect)$ is a good triple of operads.
\end{cor}

\section{Towards NonCommutative Quantization}\index{quantization} There is another possible choice of compatibility relations for which the free duplicial algebra would still be a bialgebra. It consists in taking $\CC=Dup=\AA$ and the n.u.i. compatibility relation for the \emph{four} cases 
$$(\dd_{\g}, \g), (\dd_{\d}, \g), (\dd_{\g}, \d), (\dd_{\d}, \d).$$
For this new type of bialgebras the map $\varphi(V)$ is not surjective anymore because, for $Dup(V)$, we have $\dd_{\g}= \dd_{\d}=\dd$ as described in \ref{Dupbialgebra}. Hence $\varphi(V)$ factors through $As^c(V)$. This phenomenon is similar to    $\varphi(V)$ factorizing  through $Com^c(V)$ in the $As^c\--As$ case with $\between=\between_{Hopf}$.
 
 Analogously, for $\CC=As^2=\AA$ and compatibility relations as above, the free $As^2$-algebra is a bialgebra, but the map $\varphi$ is not surjective since it factors through $As^c$ (cf.~\ref{asas2}). The notion of infinitesimal associative bialgebra (with infinitesimal compatibility relation, cf.~\cite{Aguiar04}) is going to play a role in the analysis of these bialgebras.
 
 We intend to address these cases in a forthcoming paper.



\chapter{Appendix}
\section{Types of algebras mentioned in this monograph}\label{tableauCR}

References for types of non-binary operads:

\begin{itemize}
\item Sabinin: \ref{Sabinin},
\item MagFine: \ref{MagFine},
\item $MB $: \ref{2as},
\item brace: \ref{brace},
\item $Mag^{\infty}$: \ref{as2asmag},
\item $\langle k \rangle$-ary: \ref{Askary},
\item $tAs^{\langle k \rangle}$: \ref{k-ary},
\item $pAs^{\langle k \rangle}$: \ref{k-ary}.
\end{itemize}

\newpage

The following tableau contains the binary operads mentioned in this monograph.

\bigskip

\begin{tabular}{| l | l | l | l | }
\hline
\hline
name & gen. & symmetry & relations \\
\hline
\hline
associative & $xy$ & none & $(xy)z = x(yz)$ \\
\hline
commutative & $xy$ & $xy=yx$ &  $(xy)z = x(yz)$ \\
\hline
Lie & $[x,y]$ & $[x,y]=-[y,x]$ & $ [[x,y],z]+[[y,z],x]+[[z,x],y]=0$\\
\hline
parastatistics & $xy$ & none & $(xy)z = x(yz)$\\
               &          &           & $[[x,y],z]=0$\\
\hline
NLie  & $[x,y]$ & $[x,y]=-[y,x]$ &$[[x,y],z]=0$\\
\hline
magmatic & $xy$ & none & none \\
\hline
Poisson &$ x\cdot y$ & $x\cdot y = y\cdot x$ & $(x\cdot y)\cdot z = x\cdot (y\cdot z)$ \\
             &     $ [x,y]$  & $[x,y]=-[y,x]$ &$ [[x,y],z]+[[y,z],x]+[[z,x],y]=0$\\
             &                   &                  & $[x\cdot y, z] = x\cdot [y,z]+[x,z]\cdot y$\\
             \hline
2-associative & $x\cdot y $ & none & $(x\cdot y)\cdot z = x\cdot (y\cdot z)$ \\
                     & $ x*y$ &  & $(x* y)*z = x* (y* z)$ \\
                     \hline
Zinbiel & $x\cdot y$ & none  & $(x\cdot y)\cdot z = x\cdot (y\cdot z + z\cdot y)$ \\
\hline
dendriform & $x\g y$ & none & $(x\g y)\g z = x\g (y\g z+  y\d z)$ \\
                   & $x\d y$ &           & $(x\d y)\g z = x\d (y\g z) $ \\
                     &            &           & $(x\g y + x \d y)\d z = x\d (y\d z) $ \\
                     \hline
  dipterous & $x\g y$ & none & $(x\g y)\g z = x\g (y*z)$ \\
                   & $x* y$ &           & $(x*y)* z = x*(y* z) $ \\
                   \hline
tridendriform & $x\g y$ & none & 11 relations  \\
                   & $x\d y$ &           & see \ref{tridendriform}  \\
                   &$x\cdot y$&           &   \\
                     \hline
CTD              & $x\g y$ &           & $(x\g y)\g z = x\g (y\g z + z\g y)$ \\
                   & $x\cdot y$ &   $x\cdot y = y\cdot x$          & $(x\cdot y)\g z = x\cdot (y\g z)=(x\g z)\cdot y  $ \\
                     &            &           &$(x\cdot y)\cdot z = x\cdot (y\cdot z) $ \\
                     \hline
 pre-Lie &  $ \{x;y\}$         &   none        & $\{ \{x;y\};z\} - \{x;\{y;z\}\}=$\\
              &                          &                  &      $ \{ \{x;z\};y\} - \{x;\{z;y\}\}$\\
              \hline
 NAP    &   $xy$ & none & $(xy)z = (xz)y$\\
 \hline
 Nil & $xy$ & none & $(xy)z=0=x(yz)$ \\
 \hline
 interchange &$a\circ b$    & none & $(a\circ b)\bullet (c\circ d ) = (a\bullet c) \circ (b\bullet d)$\\
                      &$ a\bullet b$ &               &    \\
 \hline
 duplicial & $x\g y$ & none & $(x\g y)\g z = x\g (y\g z)$ \\
                   & $x\d y$ &           & $(x\d y)\g z = x\d (y\g z) $ \\
                     &            &           & $(x \d y)\d z = x\d (y\d z) $ \\
 \hline
 \hline
\end{tabular}

\newpage

\section{Compatibility relations $\between$ mentioned in this monograph}\label{tableauCR}

\bigskip

\noindent Hopf (unital):

$$\cpbdeuxdeux = \cpbG $$

\noindent Hopf (nonunital):

$$\cpbdeuxdeux = \cpbA +\cpbB  +\cpbC +\cpbD +\cpbE +\cpbF +\cpbG $$

\noindent Infinitesimal (unital):

$$\cpbdeuxdeux = -\cpbA  +\cpbC  +\cpbE  $$

\noindent Infinitesimal:

$$\cpbdeuxdeux = \cpbC  +\cpbE  $$

\noindent Infinitesimal (nonunital):

$$\cpbdeuxdeux = \cpbA  +\cpbC  +\cpbE  $$

\noindent Magmatic:

$$\cpbdeuxdeux = \cpbA $$

\noindent Frobenius:

$$\cpbdeuxdeux = \cpbC = \cpbE $$

\newpage

\noindent Livernet:

$$\cpbdeuxdeux = \cpbA + \cpbC + \cpbD $$

\noindent Semi-Hopf:

$$\cpbDeuxDeux{\g}{*{}} = \cpbB + \cpbCC{*{}}{\g} + \cpbDD{*{}}{{*}} + \cpbFF{*{}}{\g}+
\cpbGG{*{}}{{*}}{*{}}{\g}$$

\noindent BiLie:

$$\cpbdeuxdeux = \cpbC + \cpbD + \cpbE + \cpbF$$

\noindent Lily:

$$\cpbdeuxdeux = 2\Big(\cpbA - \cpbB\Big) + \frac{1}{2}\Big(\cpbC + \cpbD + \cpbE + \cpbF\Big)$$

\noindent Nilpotent:
$$
\cpbdeuxdeux = \cpbA -\cpblosg - \cpblosd + \cpbloslos 
$$

\noindent Infinitesimal $3$-ary (nnunital):

\medskip

$\cpbtroistrois = \cpbtroisID + \cpbtroisAB + \cpbtroisAC + \cpbtroisBA $

\hfill $+ \cpbtroisBC + \cpbtroisCA + \cpbtroisCB$
\bigskip

\newpage 
\section{Tableau of some good triples of operads}\label{tableauTriples}

\bigskip

\begin{tabular}{|  c | c | c | c | l |   }
\hline
 $\CC$ & $\AA$ & $\PP$ & $\between$ & reference \\
\hline
 $Com$  & $Com$  & $Vect$ & Hopf  & \cite{Borel}\\
 $Com$ & $Parastat$ & $NLie$ & Hopf &\cite{LPo} \\
$Com$ & $As$ & $Lie$ & Hopf & \cite{Cartier, MM} \\ 
 $Com$ & ?? &$As$&Hopf+?  &conjectural\\ 
 $Com$ &$SGV$ &$PreLie $&Hopf+ ad hoc  &\cite{LRcha}\\
  $Com$ &$ComAs$ &$SMB$&Hopf+ad hoc  &\cite{LRcha}\\
 $Com$ & $Mag$ & $Sabinin$ & Hopf &\cite{Perez06, Holtkamp05} \\
$Com$ & $Dipt$ &$ \Prim_{Com}Dipt$ & Hopf &see \ref{ss:dipterous}\\
$Com$ & ?? &$ Com$ & ? &conjectural\\
\hline
$As$ & $As$ & $Vect$ & n.u.i. & \cite{LRstr}\\
$As$ & $PreLie$ & $\Prim_{As}PreLie$ & n.u.i. & \\
$As$ & $Mag$ & $MagFine$ & n.u.i. &\cite{HLR} \\
$As$ & $2as$ & $MB$ & n.u.i.+Hopf&\cite{LRstr} \\
$As$ & $As^2$ & $As$ & n.u.i.+n.u.i.& see \ref{asas2} \\
$As$& $DupPreLie$& $\Prim_{As}DupPreLie$&n.u.i.+n.u.i. &see \ref{AsDupPreLie}\\
$As$& $Dup$& $Mag$&n.u.i.+n.u.i. &see \ref{AsDupMag}\\
 $As$ & $2as$ & $Mag^{\infty}$ & n.u.i.+n.u.i. &see \ref{as2asmag} \\
\hline
$As$ & $Zinb$ & $Vect$ &semi-Hopf & \cite{Burgunder08}\\
$As$ &$Dipt$ &$MB$&semi-Hopf  &\cite{LRnote, LRcha}\\
$As$ & $Dend$ & $Brace$ &semi-Hopf & \cite{Ronco02}, see \ref{ss:dendriform}\\
$As$ &$GV$ &$Brace$&semi-Hopf  &\cite{LRcha}\\
$As$ & $Dipt$ & $MB$ &semi-Hopf & \cite{LRnote}\\
$As$ & $CTD$ & $Com$ &semi-Hopf&\cite{JLLctd} \\
$As$ & $Tridend$ & $Brace+As$ &semi-Hopf &\cite{PR} \\
\hline
$Zinb$ & $As$ & $Vect$ &semi-Hopf & \cite{Burgunder08}\\
\hline
$Mag$& $Mag$& $Vect$&magmatic &\cite{Burgunder06} \\
\hline
$NAP$& $PreLie$& $Vect$&Liv &\cite{Livernet06} \\
$NAP$& $Mag$& $\Prim_{NAP}Mag$&Liv & see \ref{NAPA} \\
\hline
$2as$& $2as$& $Vect$& n.u.i.${}^2$+Hopf${}^2$&\cite{LRstr} \\
 $2as$& $2mag$& $\Prim_{2as}2mag$& n.u.i.${}^2$+Hopf${}^2$&\\
\hline
 $Dend$& $Dend$& $Vect$& hemisemi-Hopf&\cite{Foissy} \\
 $Dend$& $2mag$& $\Prim_{Dend}2mag$ & hemisemi-Hopf& \\
\hline
 $Lie$& $Lie$& $Vect$& Lily (not biLie)&see \ref{LieLie} \\
  $Lie$& $PreLie$& ?? & Lily+?&conjectural \\
 $Lie$& $Mag$& $\Prim_{Lie}Mag$& Lily (not biLie)&see \ref{LieLie} \\
$Lie$& $PostLie$& $\Prim_{Lie}PostLie$&Lily+? &\cite{Vallette05} \\
\hline
$Nil$&$Nil$ &$Vect$ &Nil &see \ref{NilNil} \\
$Nil$&$Mag$ &$Mag^{0,2,1}$ &Nil & see \ref{NilMag}\\
$Nil$&$Dup^!$ &$As$ &? &conjectural \\
 \hline
$IC$ & $IC$ & $Vect$  & Hopf+ad hoc &see \ref{IC}, \cite{Lafont97}\\
 \hline
$tCom^{\langle k \rangle}$ &$tAs^{\langle k \rangle}$&$tLie^{\langle k \rangle}$  &Hopf style & see \ref{k-ary} \\
  \hline
 $tAs^{\langle k \rangle}$ &$tAs^{\langle k \rangle}$& $Vect$ & n.u.i. & see \ref{k-ary} \\
  \hline
   $???$ &$PreLie$& $Lie$& ? &conjectural\\
    $???$ &$Lie\pt adm$& $Lie$& ? &conjectural\\
\hline
\end{tabular}
\vskip 1cm

The notation $\Prim_{\CC}\AA$ in the column ``$\PP$'' means that we do not know yet about a small presentation of this operad. 
\newpage

\backmatter


\printindex



\end{document}